\documentclass[11pt,a4paper,leqno]{article}
\usepackage{amsmath, ntheorem}
\usepackage{amssymb}
\usepackage{color}
\advance\topmargin by -2.2cm
\theoremstyle{change}
\newtheorem{thm}{Theorem}[section]
\newtheorem{lemma}[thm]{Lemma}
\newtheorem{assumption}[thm]{Assumption}

\newtheorem{prop}[thm]{Proposition}
\newtheorem{hyp(H)}[thm]{Hypothesis (H)}

\theorembodyfont{}
\newtheorem{defn}[thm]{Definition}

\newtheorem{rem}[thm]{Remark}

\numberwithin{thm}{subsection}

\newcommand{\bbr}{I\!\!R}

\newcommand{\bbn}{I\!\!N}

\newcommand{\cala}{{\cal A}}
\newcommand{\calb}{{\cal B}}
\newcommand{\calc}{{\cal C}}

\newcommand{\cale}{{\cal E}}
\newcommand{\calf}{{\cal F}}
\newcommand{\calg}{{\cal G}}
\newcommand{\calh}{{\cal H}}

\newcommand{\calj}{{\cal J}}

\newcommand{\call}{{\cal L}}

\newcommand{\caln}{{\cal N}}
\newcommand{\calo}{{\cal O}}

\newcommand{\cals}{{\cal S}}

\newcommand{\calu}{{\cal U}}
\newcommand{\calv}{{\cal V}}

\newcommand{\calx}{{\cal X}}
\newcommand{\caly}{{\cal Y}}
\newcommand{\calz}{{\cal Z}}

\newcommand{\barr}{\begin{array}}
\newcommand{\earr}{\end{array}}
\newcommand{\beqq}{\begin{equation}}
\newcommand{\eeqq}{\end{equation}}
\newcommand{\beao}{\begin{eqnarray*}}
\newcommand{\eeao}{\end{eqnarray*}\noindent}
\newcommand{\beam}{\begin{eqnarray}}
\newcommand{\eeam}{\end{eqnarray}\noindent}

\newcommand{\halmos}{\quad\hfill\mbox{$\Box$}}

\newcommand{\la}{\lambda}

\newcommand{\si}{\sigma}
\newcommand{\al}{\alpha}
\newcommand{\vth}{\vartheta}

\newcommand{\om}{\omega}
\newcommand{\Om}{\Omega}

\newcommand{\vep}{\varepsilon}

\newcommand{\wh}{\widehat}
\newcommand{\wt}{\widetilde}

\newcommand{\ov}{\overline}

\newcommand{\lra}{\longrightarrow}

\newcommand{\nto}{n\to\infty}

\begin{document}

\vskip1.5cm

{\Large\bf 
Ergodic branching diffusions with immigration: \\ 
properties of invariant occupation measure, \\
identification of particles under high-frequency observation,  \\
and estimation of the diffusion coefficient at nonparametric rates 
}

\vskip1.5cm 
Matthias Hammer 
({\tt matthias.hammer@tu-berlin.de})\footnote{Institut f\"ur Mathematik, Technische Universit\"at Berlin, Stra\ss e des 17. Juni 136, 10623 Berlin, Germany.}
\\
Reinhard H\"opfner ({\tt hoepfner@mathematik.uni-mainz.de})\footnote{Institut f\"ur Mathematik, Johannes Gutenberg-Universit\"at Mainz, Staudingerweg 9, 55099 Mainz, Germany.} 
\\
Tobias Berg ({\tt tobias.berg@aol.de})\footnotemark[2]

\vskip1.5cm

{\small 
{\bf Abstract: }
In branching diffusions with immigration (BDI), particles travel on independent diffusion paths in $\bbr^d$, branch at position-dependent rates and leave offspring --randomly scattered around the parent's death position-- according to position-dependent laws. We specify a set of conditions which grants ergodicity such that the invariant occupation measure is of finite total mass  and admits a continuous Lebesgue density.   

Under discrete-time observation, BDI configurations being recorded at discrete times $i\Delta$ only, $i\in\bbn_0$, we lose information about particle identities between successive observation times. We present a reconstruction algorithm which in a high-frequency setting (asymptotics $\Delta\downarrow 0$) allows to reconstruct correctly a sufficiently large proportion of particle identities, and thus allows to recover $\Delta$-increments of  unobserved diffusion paths on which particles are travelling. 
Picking some few well-chosen observations we fill regression schemes which, on cubes $A$ where the invariant occupation density is strictly positive, allow to estimate the diffusion coefficient of the one-particle motion at  nonparametric rates. 
\\

\\

high frequency observation, reconstruction algorithm, \\ estimation of the diffusion coefficient, kernel estimators, nonparametric rates
\\
}

\newpage

Branching diffusions with immigration (BDI) as strong Markov processes of finite   configurations of particles in $\bbr^d$ have been investigated since Ikeda, Nagasawa and  Watanabe (\cite{INW-66a}, \cite{INW-66b}, \cite{INW-68}, \cite{INW-69}) who study semigroups and their generators, and construct the process by killing and repasting of strong Markov processes. In our BDI model, particles travel on independent diffusion paths, branch at position-dependent rates and leave offspring --randomly scattered around the parent's death position-- according to position-dependent laws; 
with some variations, this is close to \cite{Lo-02}, \cite{Lo-02a}, \cite{Lo-04}, \cite{HHL-02}, \cite{HL-05}, \cite{Ha-12}.  
L\"ocherbach (\cite{Lo-02}, \cite{Lo-02a}) studies likelihood ratios for BDI processes and gives conditions for convergence to limit likelihood ratios, of some type (local asymptotic normality) which is important in parametric statistics, and knowledge of invariant measure or invariant occupation measure is required to check statistical assumptions. 
H\"opfner, Hoffmann and L\"ocherbach \cite{HHL-02} focus on the point process of branching times/positions and estimate nonparametrically the spatial branching rate: again the statistical assumptions require sufficiently explicit knowledge of the invariant occupation measure and its moments. 
L\"ocherbach \cite{Lo-04} and H\"opfner and L\"ocherbach \cite{HL-05} study the invariant occupation measure on the single-particle space and obtain a continuous Lebesgue density in different settings. Both papers consider the case of local branching, i.e.\ newborn particles start at their parent's death position.  \cite{Lo-04} works with binary reproduction, uniform ellipticity, uniformly bounded $\calc^\infty$-coefficients and Malliavin calculus; \cite{HL-05} with general position-dependent reproduction laws, finite order of smoothness of coefficients, but  restrictive assumptions in view of duality techniques. Assuming that newborn particles scatter randomly around the parent's death position, 
Hammer \cite{Ha-12} obtains a continuous density for the invariant measure on the configuration space by means of Fourier methods.

Our paper gives a set of relatively general conditions (such as: uniform exponential stability of the expectation semigroup, heat kernel bounds for the single-particle motion subject to position-dependent killing, \ldots; see section \ref{2theorems}) which grant ergodicity of the BDI process, provide some information on finite `moments' of the invariant measure up to some order (which is explicit from the family of position-dependent reproduction laws, see theorem \ref{1.3.7}),  
and imply that the invariant occupation measure is of finite total mass and admits a continuous Lebesgue density (see theorem \ref{1.3.10}). The `spatial subcriticality condition' of \cite{HL-05} reappears as main condition which grants ergodicity of the BDI process with finite invariant occupation measure (see assumption \ref{1.2.2} and lemma \ref{1.2.5} in section \ref{1.2}), however, with an important difference: in our model, similar to \cite{Ha-12}, we allow for essentially arbitrary non-local branching, i.e.\ offspring can be scattered randomly around the death position of the parent particle. As a consequence, the process entering our `spatial subcriticality condition' and thus determining the shape of the invariant occupation measure is no longer the single-particle motion itself but a jump diffusion whose jumps represent --in a sense of a `many-to-one'-formula-- the location of a `typical child' relative to the parent particle. 
Our approach allows us to avoid, to a large extent, restrictive smoothness assumptions on the coefficients or ellipticity conditions (of course, uniform ellipticity or smoothness `of low order' may be welcome as a sufficient condition to check e.g.\ our heat kernel bound assumption \ref{1.3.8} in section \ref{2theorems}).

We then turn to a setting which has received a lot of attention in statistics of processes: if we observe a continuous-time process at discrete time points $i\Delta$ only, $i\in\bbn_0$, how do we estimate those quantites which under continuous-time observation would induce mutual singularity of the laws of the process? 
The main example is volatility in diffusion processes: see e.g.\ Yoshida \cite{Yo-92}, Genon-Catalot and Jacod \cite{GJ-93}, Gobet \cite{Go-02}, Podolskij and Vetter \cite{PV-10}, Protter and Jacod \cite{JP-12}. In our situation of BDI processes $( \eta_t )_{t\ge 0}$ which are configuration-valued, observation at discrete times forces us to consider a new type of problem.  
If we observe at discrete time points $\,\left\{ i\Delta : i\in\mathbb{N}_0\right\}\,$ not a diffusion path but the trajectory of a BDI process $( \eta_t )_{t\ge 0}$, we will be left with pairs of configurations $(\eta_{i\Delta},\eta_{(i+1)\Delta})$ without any information on the path history in-between; these are merely pairs of random point measures on the single-particle space. 
Information on branching or immigration events inside $(i\Delta,(i{+}1)\Delta)$ is lost; even in case all particles succeeded to stay alive over the whole time interval $[i\Delta,(i{+}1)\Delta]$, we do not know which particle of the first configuration did travel to which position of the second configuration. In this context, we propose a reconstruction algorithm (see \ref{2.2.2} and theorems \ref{2.2.4} and \ref{2.2.6} in section \ref{2.1}) 
which in case of high-frequency asymptotics (i.e.\ $\Delta\downarrow 0$) will be able to recover correctly, to some large extent, particle identities in pairs $(\eta_{i\Delta},\eta_{(i+1)\Delta})$ of successive configurations.

In a next step, we make use of the reconstruction algorithm and of the setting of high-frequency asymptotics $\Delta\downarrow 0$ to fill regression schemes for estimation of the diffusion coefficient  $\,\si \si^\top$ of the single-particle motion in the BDI process, picking out of the overwhelming amount of discrete-time data $( \eta_{i\Delta} )_{i\in\bbn_0}$ some few but well-selected pairs $(\eta_{i\Delta},\eta_{(i+1)\Delta})$ of successive configurations for which we are sure --up to exceptional sets of vanishing probability as $\Delta\downarrow 0$-- that we reconstruct the particle identities correctly, for all particles involved. 
Reconstructing in this way particle identities and thus $\Delta$-increments for the trajectory of these particles, our regression scheme 
(see \ref{3.1.1NEU} and theorem \ref{3.1.4NEU} in section \ref{3.1NEU}) 
consists of particles indexed by $\al$ in some index set $\calj(\Delta)$ associated to $\Delta$, and of pairs 
$$
(\, \calx_\al \,,\, \calz_\al \,) \quad,\quad \al\in\calj(\Delta)
$$
where $\calx_\al$ represents the position of particle $\al$ at some random time $\tau_\al\Delta$ and $\calz_\al$ the reconstructed (rescaled) increment $\frac{\xi_\al((\tau_\al+1)\Delta)-\xi_\al(\tau_\al\Delta)}{\sqrt{\Delta\,}}$ for particle $\al$ on  $[\tau_\al\Delta,(\tau_\al{+}1)\Delta]$; note that the trajectory itself remains unaccessible from discrete-time data $( \eta_{i\Delta} )_{i\in\bbn_0}$. For fixed cubes $A$ in the single-particle space on which the invariant occupation density is strictly positive, we can ensure that  the `design variables' $\calx_\al$, $\al\in\calj(\Delta)$, are approximately regularly spaced over $A$.  
If we associate to particles $\al$ their driving Brownian motion $W_\al$,  results due to Jacod and Genon-Catalot \cite{GJ-93}, see also Podolskij and Vetter \cite{PV-10}, yield good approximations of type 
$$
\calz_\al := \frac{\xi_\al((\tau_\al+1)\Delta)-\xi_\al(\tau_\al\Delta)}{\sqrt{\Delta\,}}
\;\;\approx\;\;
\si(\calx_\al)\; \frac{W_\al((\tau_\al+1)\Delta)-W_\al(\tau_\al\Delta)}{\sqrt{\Delta\,}} 
\quad,\quad \calx_\al := \xi_\al(\tau_\al\Delta)  
$$
which give  
$$
\calz_\al \calz^\top_\al \;\;\approx\;\; (\,\si \si^\top)(\calx_\al) \quad+\quad \mbox{error terms with some martingale structure} \;. 
$$
All this holds on the `good sets' where our reconstruction is indeed correct, i.e.\ on the complements of exceptional sets. If the probability of the exceptional sets vanishes as $\Delta\downarrow 0$, the contribution of what we believe --falsely on the exceptional set-- to be an increment does not vanish: whereas in restriction to the `good sets' density estimation of $\,\si\si^\top(\cdot)$ works as in classical iid regression schemes for nonparametric estimation (Tsybakov \cite{Tsy-08}), we have to take care of what happens on the exceptional sets in order to reach a balance of both types of contributions to the squared pointwise risk of a nonparametric estimator.  We make this explicit in dimension $d=1$ when we discuss kernel estimation of $\si^2(\cdot)$ on intervals $A$ on which the (continuous) invariant occupation density is strictly positive (theorem \ref{3.2.2'NEU} in section \ref{3.2NEU}).

The paper is organized as follows. The setting for ergodic BDI processes is exposed in section \ref{setting}. Continuity of the Lebesgue density of the invariant occupation measure is proved in section \ref{1.3}. The reconstruction algorithm for particle identities in discretely observed BDI processes is the topic of section \ref{reconstruction}. Section \ref{regressionNEU} deals with regression schemes, filled from discrete observations, with the aim to estimate the diffusion coefficient of the single-particle motion; the example of kernel estimation of the diffusion coefficient in dimension $d=1$ appears in section \ref{3.2NEU}.


\section{ Ergodic branching diffusions with immigration: our setting}\label{setting}

This section introduces branching diffusions with immigration (BDI) and their ergodicity properties. 
In a first subsection, we introduce BDI as strong Markov processes with life time, close to \cite{Lo-02}, \cite{Lo-02a}, \cite{Lo-04}, \cite{HL-05} but more general in that we allow for quite arbitrary spatial scattering of the descendants generated at a branching event (as in \cite{Ha-12}). 
Our method is a construction by killing and repasting of strong Markov processes as in \cite{INW-66b}, \cite{INW-68}, \cite{INW-69} or \cite{Na-77}. 
In a second subsection we state a `spatial subcriticality' condition which grants positive Harris recurrence with finite invariant occupation measure. A third subsection sketches proofs as far as their techniques are relevant for the rest of the paper.

\subsection{BDI processes as strong Markov processes with life time}\label{1.1}

For $d\ge 1$, we write $E:=\bbr^d$, $\cale:=\calb(\bbr^d)$ and call $(E,\cale)$ single particle space. We call the space of (ordered) particle configurations  $S:=\bigcup_{\ell\in\bbn_0}E^\ell$ configuration space; $\delta$ denoting the void configuration, we have $E^0=\{\delta\}$. We write $\cals:=\calb(S)$ for the Borel-$\si$-field on $S$: $(S,\cals)$ is a Polish space. Lebesgue measure on $S$ is defined layer-wise (for $\ell\ge 1$, its restriction to $E^\ell$ equals Lebesgue measure on $E^\ell$). 
The length of a configuration $x\in S$ is denoted by $\ell(x)$, i.e.\ $\ell(x)=j$ iff $x\in E^j$. Sometimes we write a configuration $x\in S$ as a point measure $x(A)=\sum_{j=1}^{\ell(x)}\epsilon_{x_j}(A)$,  $A\in \cale$ (which equals $0$ if $x=\delta$). 
To measurable functions $f:E\to\bbr$ we associate functions $\ov f:S\to\bbr$ by $\ov f(x):=x(f)$, i.e.\ 
$$
\ov f(x):=\sum_{j=1}^\ell f(x_i) \quad\mbox{when $x=(x_1,\ldots,x_\ell)\in S$}\;,\quad \ov f(\delta):= 0 \;. 
$$
With these notations, BDI will be a $(S,\cals)$-valued c\`adl\`ag strong Markov process with the following properties {\bf (A1)}--{\bf (A4)}: 

{\bf (A1)} For $\ell\in\bbn$, on some random time interval which is specified through {\bf (A2)} and {\bf (A3)} below, $\ell$-particle configurations $X^\ell$ travel in $E^\ell$ as (strong) solutions to 
$$
X^\ell_t=(X^{1,\ell}_t,\ldots,X^{\ell,\ell}_t) \quad\mbox{satisfying}\quad dX^{j,\ell}_t = b(X^{j,\ell}_t)dt + \si(X^{j,\ell}_t)dW^j_t \;,\; 1\le j\le \ell \;. 
$$ 
Here $W^1,\ldots,W^\ell$ are independent $d$-dimensional Brownian motions. Drift $b:\bbr^d\to \bbr^d$ and volatility $\si:\bbr^d\to \bbr^{d\times d}$ are assumed to be Lipschitz. 

{\bf (A2)} i) Independently of each other,  particles living at the same time are killed at  position-dependent rate $\kappa : E \to (0,\infty)$; we assume that $\kappa$ is measurable and locally bounded on $E$. \\
ii) We have a transition probability  $K_1(\cdot,\cdot)$ from $E$ to $\bbn_0$ such that $p_k(y):=K_1(y,\{k\})$ gives the probability for a particle killed in position $y\in E$ to produce $k$ offspring, $k\in\bbn_0$. \\
iii) 
We have a transition probability  $K_2(\cdot,\cdot)$ from $E\times \bbn_0$ to $S$ with the property 
$$
K_2((y,k),\cdot)\text{ is concentrated on }E^k \text{ for all }y\in E,\; k\in\bbn_0
$$
which scatters offspring generated at a branching event relative to the parent's position:  $\,k$-particle offspring of a particle killed in position $y$ will be located in positions
\beqq\label{general_offspring_positions}
\quad y+v_1,\ldots,y+v_k \quad\mbox{with probability}\quad   K_2((y,k),dv_1,\ldots,dv_k) \;\;,\;\; k\ge 1 \;,\; y\in E\;, 
\eeqq
and in case $k=0$ we put $K_2((y,0),\cdot)=\epsilon_\delta(\cdot)$. \\

An important special case contained in {\bf (A2)}~iii) is given by product kernels
\begin{equation}\label{eq:K_special}
K_2((y,k),dv_1,\ldots,dv_k) \;=\; \bigotimes_{j=1}^kK(y,dv_j)
\end{equation} 
for some fixed transition probability $K(\cdot,\cdot)$ on $(E,\cale)$. 
Specializing further, 
if $K(y,dv)=q(dv)$ for some probability measure $q$ on $(E,\cale)$, independently of $y$, the distribution of newborn particles relative to their parent's position is spatially homogeneous. 
Finally, $\,q(dv):=\epsilon_0(dv)$ is the commonly considered case that particles are born exactly at the death position of their parent; we will refer to this special case as \emph{local branching}.
In this paper, we shall work under \eqref{general_offspring_positions} and shall not even assume  \eqref{eq:K_special}, i.e.\ we allow for arbitrary non-local branching mechanisms.\\

{\bf (A3)} For some probability measure $Q^{\tt i}$ on $(E,\cale)$ and some constant $0<c<\infty$, single immigrants arrive at constant rate $0<c<\infty$ and are located in $E$ according to $Q^{\tt i}(dy)$, independently of everything else. \\

Write $(\Om,\cala)$ for the canonical path space of c\`adl\`ag functions $[0,\infty)\to S$ with life time $\zeta\le\infty$, and $\eta=(\eta_t)_{t\ge 0}$ for the canonical process on $(\Om,\cala)$. Then {\bf (A1)}--{\bf (A3)} above determine uniquely a probability measure $Q$ on $(\Om,\cala)$ under which $\eta$ is a jump diffusion with life time. As long as $t<\zeta$, jumps  (finitely many on finite time intervals) arrive at rate 
$$
c + \ov\kappa(\eta_t) \;=\; c + \sum_{j=1}^\ell\kappa(\eta^i_t) 
\quad\mbox{when $\eta_t=(\eta^1_t, \ldots, \eta^\ell_t)$ belongs to $E^\ell$} \;. 
$$
By convention, the rate is $c$ when $\ell=0$. Note that by the Lipschitz assumptions on drift and diffusion coefficient in {\bf (A1)}, and by the local boundedness of $\kappa$ in {\bf (A2)}, the life time $\zeta\le\infty$ of the process is the first accumulation point of the sequence of successive jumps times $T_j$: we have  $T_j< T_{j+1}$ as long as $T_j$ is finite, and $\zeta:=\sup_j T_j \le \infty$. On events $\{T_j<\infty\}$, representing the configuration $\eta_{T_j^-}$ immediately before the jump by $x = (x^1,\ldots,x^\ell)$, $\ell\ge 1$, the new configuration $\eta_{T_j}$ at time $T_j$ is obtained from $x$ as follows:  
$$
\left\{ \begin{array}{lll}
(x^1, \ldots, x^{j-1}, x^{j+1}, \ldots, x^\ell)  & \mbox{w.\ pr.}\quad   & \frac{\kappa(x^j)\, p_0(x^j)}{c+\ov\kappa(x)} \;, \\
(x^1, \ldots, x^{j-1}, x^j+v_1, \ldots, x^j+v_k, x^{j+1}, \ldots, x^{\ell})  & \mbox{w.\ pr.}\quad    &\frac{\kappa(x^j)\, p_k(x^j)}{c+\ov\kappa(x)}\, K_2((x^j,k),dv_1,\ldots,dv_k)
\;, \\
(x^1, \ldots, x^\ell, y)   & \mbox{w.\ pr.}\quad    & \frac{c}{c+\ov\kappa(x)}\, Q^{\tt i}(dy) \;. 
\end{array}\right.
$$
With exception of  $k$ chosen equal to $1$, jumps change the length of the configuration. 
In case $\ell=0$, $\,\eta_{T_j^-}$ is the void configuration $\delta$, thus $\eta_{T_j}$ a one-particle configuration with $y$ selected by $Q^{\tt i}(dy)$. \\

When we deal with trajectories of individual particles in the BDI process, we write 
$$
d\xi_t \;=\; b(\xi_t) dt + \sigma(\xi_t) dW_t
$$
for the single-particle motion on $E$. Assumption {\bf (A1)} on drift $b(\cdot)$ and diffusion coefficient $\sigma(\cdot)$ grants that the diffusion $\xi$ has infinite life time.   \\

{\bf (A4)}\;  i) We have $\,\int_0^\infty \kappa(\xi_s)\, ds = \infty\,$ almost surely, for every choice of a starting point $y\in E$ for $\xi$. 

ii) Reproduction means $y \to \rho(y) := \sum_{k\ge 0}k\, p_k(y)$ are (finite and) locally bounded on $E$. \\

Condition {\bf (A4)}~i) grants that all $T_j$ defined above are almost surely finite stopping times.

Assumptions {\bf (A1)}--{\bf (A4)} 
and all notations of the present subsection will hold throughout the paper. So far, our construction of the BDI process is the canonical one:  $(\Om,\cala)$ is the canonical path space of c\`adl\`ag functions $[0,\infty)\to S$ with life time $\zeta\le\infty$,  $\,\eta=(\eta_t)_{t\ge 0}$ is the canonical process on $(\Om,\cala)$, and we have a unique probability law $Q$ on  $(\Om,\cala)$ such that $\,\eta\,$ is strongly Markov with the above properties: a jump diffusion with successive jump times $(T_j)_j$ which are finite stopping times and increase towards $\zeta\le\infty$.

\subsection{Ergodicity}\label{1.2}

In this subsection, we state a set of sufficient conditions which ensure that the BDI process $\eta$\\  
\quad$i)$ is positive Harris recurrent, admitting the void configuration $\delta$ as a recurrent atom (thus in par\-ti\-cu\-lar, $\eta$ will have infinite life time $\zeta=\infty$); \\
\quad$ii)$ admits a finite invariant occupation measure on $(E,\cale)$. 

Up to the general form of our kernel $K_2(\cdot,\cdot)$ in {\bf (A2)} iii), we follow the same approach as L\"ocherbach \cite{Lo-02},  \cite{Lo-02a}, \cite{Lo-04}, 
H\"opfner and L\"ocherbach \cite{HL-05} section 1.4; see also Hammer \cite{Ha-12} section 3 
where the same general form of non-local branching was allowed.
Introducing the necessary notation we state the relevant results.\\

\begin{assumption}\label{1.2.1neo} 
The functions $\,\kappa\,$ and $\,\rho\,$ are bounded on $E$. 
\end{assumption}
This assumption guarantees in particular non-explosion of the process $\eta$: by \ref{1.2.1neo}, $\,\eta\,$ has infinite life time $\zeta=\infty$ almost surely (and from here on, we will take $\,\Om$ as the usual Skorohod path space of c\`adl\`ag functions $[0,\infty)\to S$). 
Using the kernel $K_2(\cdot,\cdot)$ from assumption A2 iii), 
we define a transition probability $Q^{\tt r}(\cdot,\cdot)$ on the single-particle space $(E,\mathcal{E})$ as follows: for $y\in E$ and $f:E\to[0,\infty)$ measurable, let
\begin{equation}\label{defn:Q_r}
Q^{\tt r}(y,f):=\frac{1}{\rho(y)}\sum_{k\in\mathbb{N}}p_k(y)\int_{E^k} \big(f(y+v_1)+\cdots+f(y+v_k)\big)\, K_2\left((y,k),dv_1,\ldots, dv_k\right) \;; 
\end{equation}
if $y\in E$ is such that $\rho(y)=0$ we put $Q^{\tt r}(y,\cdot):=\nu(\cdot)$, for some fixed probability measure $\nu$ on $E$. \\

In the following we write $\call$ for the Markov generator of the single-particle-motion $\xi$ on $E=\bbr^d$
\beqq\label{generator-1}
\call f (y) \;=\; \sum_{i=1}^d b_i(y)\, \partial_i f(y)  \;+\; \frac12  \sum_{i,j=1}^d a_{i,j}(y)\, \partial_{i,j} f(y)
\eeqq
where $a=\si \si^\top$. 
Let us introduce a jump diffusion $\wt\xi$ on $E$ by defining a generator 
\beqq\label{generator-2}
\wt\call f (y) \;:=\; \call f (y) \;+\; \kappa(y)\rho(y)\int_E \left[ f(w)-f(y) \right]Q^{\tt r}(y,dw)  \;; 
\eeqq
here, with  $K_2(\cdot,\cdot)$ from {\bf (A2)}~iii), the integral contribution equals  
$$
\kappa(y) \sum_{k\in\mathbb{N}}p_k(y)\int_{E^k} \left[ \big(f(y+v_1)+\cdots+f(y+v_k)\big) - f(y) \right]\, K_2\left((y,k),dv_1,\ldots, dv_k\right) \;.
$$
These generators should be understood in the sense of the corresponding martingale problems.
The jump diffusion $\wt\xi$ can be defined probabilistically in the sense of killing and repasting of strong Markov processes,  
cf.\ Ikeda, Nagasawa and Watanabe \cite{INW-66a}, \cite{INW-66b}, \cite{INW-68}, \cite{INW-69} and Nagasawa  \cite{Na-77}:
a diffusive motion according to  $\xi$ is killed at position-dependent rate $\,\kappa\rho\,$ and restarted in a position selected by $\,Q^{\tt r}(y,dw)\,$, independently of everything else. 
Since $\xi$ has infinite life time and since $\kappa\rho$ is bounded in virtue of assumption \ref{1.2.1neo}, also the jump diffusion $\wt\xi$ has infinite life time. \\

\begin{assumption}\label{1.2.2}
We assume 
\beqq\label{erg-1}
E_y\left( \int_0^\infty e^{-\int_0^t [\kappa(1-\rho)](\wt\xi_s)\, ds}\; dt \right) \;<\; \infty \quad\mbox{for all $y\in E$} \;,   
\eeqq
\beqq\label{erg-2}
y \;\to\;  E_y\left( \int_0^\infty e^{-\int_0^t [\kappa(1-\rho)](\wt\xi_s)\, ds}\; dt \right) \quad\mbox{belongs to $L^1(Q^{\tt i})$} \;. 
\eeqq 
\end{assumption}

\vskip0.5cm
If we think of $\kappa(1-\rho)$ as a rate of annihilation/creation of mass, \eqref{erg-1} or \eqref{erg-2} deal with the total mass of the $\kappa(1-\rho)$-resolvent kernel of the jump diffusion $\wt\xi$. Assumption \ref{1.2.2} generalizes condition (6) in \cite{HL-05} to kernels satisfying {\bf (A2) }iii). It implies, see lemma \ref{1.2.4} below, `spatial subcriticality' in the sense of almost certain extinction of families starting from one ancestor located in $y\in E$. 
\\

\begin{defn}\label{1.2.3}
We shall write $\eta^{\tt r}$ for the branching diffusion $\eta$ \emph{without immigration} arising as subprocess of all direct descendants of one or several ancestors at some initial time $0$. When there is need to specify positions $y$ for one or $y_1,\ldots, y_m$ for several ancestors, we write $\eta^{{\tt r},y}$ or $\eta^{{\tt r},y_1,\ldots, y_m}$.  
\\
Recalling notation $\,\ov f( x ) = \sum_{j=1}^\ell f(x_j) = x(f)\,$ for $x=(x_1,\ldots, x_\ell)\in S$ with convention $\ov f( \delta )=0$, 
let
\beqq\label{eq:Def_H} 
H^{\tt r}(y,f) \;:=\; 
E_y\left( \int_0^\infty \ov f( \eta^{\tt r}_t )\, dt \right) \;\le\; \infty, \qquad y\in E,\;\mbox{$f:E\to[0,\infty)$ measurable} 
\eeqq
denote the expected occupation measure (finite or not) for $\eta^{\tt r}$ starting from one ancestor in $y\in E$.
\end{defn}

\vskip0.8cm
\begin{lemma}\label{1.2.4}
Under assumptions \ref{1.2.1neo} and (\ref{erg-1}) of \ref{1.2.2}, the total mass of the expected occupation measure for the progeny of an ancestor starting in position $y\in E$ is finite: We have
$$
 H^{\tt r}(y,1) \;=\; 
E_y\left( \int_0^\infty  \ell( \eta^{\tt r}_t )\, dt \right) \;\;=\;\; E_y\left( \int_0^\infty e^{-\int_0^t [\kappa(1-\rho](\wt\xi_s)\, ds}\; dt \right) \;<\; \infty \;. 
$$
\end{lemma}

\vskip0.8cm
\begin{lemma}\label{1.2.5}
Under assumptions \ref{1.2.1neo} and \eqref{erg-2} of \ref{1.2.2}, the BDI process $\eta$ is positive Harris recurrent, admits the void configuration $\delta$ as a recurrent atom, and 
has finite invariant occupation measure 
\beqq\label{explicit_form_mubar}
\ov\mu(A) \;=\; c\; [Q^{\tt i} H^{\tt r}](A) \;=\; c\, \int_E Q^{\tt i}(dy) H^{\tt r}(y,1_A) \;<\; \infty \quad,\quad A\in\cale
\eeqq 
with $\,c \,,\, Q^{\tt i}$ of {\bf (A3)} and $H^{\tt r}(\cdot,\cdot)$ given by (\ref{eq:Def_H}). The choice of the constant in \eqref{explicit_form_mubar} relates $\,\ov\mu\,$ to the invariant probability $\,\mu\,$ of the BDI process through 
$$
\ov\mu(f)\; =\; \mu(\ov f) \quad\mbox{for all $f:E\to[0,\infty)$ measurable}\;. 
$$   
\end{lemma}

\vskip1,0cm
\subsection{Sketching the proofs, and some further notation}\label{1.3neo}

This subsection will sketch proofs for lemmata \ref{1.2.4} and \ref{1.2.5} --assertions which generalize results from \cite{HL-05} to kernels according to {\bf (A2)}~iii)-- as far as the techniques which appear are of importance for the rest of the paper. \\

{\bf Proof of lemma \ref{1.2.4}: } 
Consider the process $\eta^{\tt r}$ starting from one ancestor in $y\in E$.
By {\bf (A4)}, the time $\tau$ of the first branching event in $\eta^{\tt r}$ is a.s.\ finite, thus $\,\eta^{\tt r}_{\tau^-}$ ($E$-valued) and $\,\eta^{\tt r}_\tau\,$ ($S$-valued) are well-defined random variables. 
The strong Markov property yields 
\begin{align}\begin{aligned}\label{iterate-1} 
H^{\tt r}(y,f) &\;=\; E_y\left( \int 1_{[[0,\tau[[}(t)\, f(\eta^{\tt r}_t)\, dt \;+\; E_{\eta_\tau^{\tt r}}\left(\int_0^\infty \ov f(\eta_t^{\tt r})\,dt\right) \right)\\
&\;=\; E_y\left( \int 1_{[[0,\tau[[}(t)\, f(\eta^{\tt r}_t)\, dt \;+\; \overline{H^{\tt r}(\cdot,f)}(\eta^{\tt r}_\tau) \right),
\end{aligned}\end{align}
where we combine definition \eqref{eq:Def_H} of $H^{\tt r}(\cdot,\cdot)$   with the \emph{branching property} 
(i.e.\ the fact that particles evolve independently). Writing the second contribution on the right hand side conditionally on $\eta^{\tt r}_{\tau^-}=z$ as 
$$
\sum_{k=1}^\infty p_k(z)\int_{E^k} K_2\left( (z,k) , dv_1,\ldots,dv_k \right) \left[ H^{\tt r}(z+v_1,f) + \ldots + H^{\tt r}(z+v_k,f)\right] 
$$
which by definition of $Q^{\tt r}(\cdot,\cdot)$ in \eqref{defn:Q_r} equals
$$
\rho(z)\;  \int_E Q^{\tt r}(z,dw)\; H^{\tt r}(w,f)  \;,  
$$
equation (\ref{iterate-1}) takes the form  
\beqq\label{iterate-2} 
H^{\tt r}(y,f) \;=\; E_y\left( \int1_{[[0,\tau[[}(t)\, f(\eta^{\tt r}_t)\, dt \;+\; \rho(\eta^{\tt r}_{\tau^-})  \int_E Q^{\tt r}(\eta^{\tt r}_{\tau^-},dw)\, H^{\tt r}(w,f)  \right)  \;. 
\eeqq  
Note that conditionally on $\eta^{\tt r}_0=y=\xi_0 \in E$, and up to time $\tau$ of position-dependent killing at rate $\kappa$, $\;\eta^{\tt r}$ is a single-particle diffusion $\xi$. Introducing the  $\kappa$-resolvent kernel of the diffusion $\xi$  
\begin{equation}\label{defn:R_kappa}
 R_\kappa(y,f) \;:=\; E_y\left(\int_0^\infty f(\xi_t) e^{-\int_0^t \kappa(\xi_s)\,ds}\,dt\right) 
\;=\; E_y\left(\int_0^\tau f(\xi_t)\,dt\right)  
\end{equation}
($y\in E$, $f:E\to[0,\infty)$ measurable), the law of $\,\eta_{\tau-}^{\tt r}\,$ starting from $\eta_0^{\tt r}=y\in E$ is given by  $[R_\kappa\kappa](y,\cdot)$, and we can rewrite \eqref{iterate-2} as
\begin{equation}\label{iterate-2a}
H^{\tt r}(y,f) \;=\; R_\kappa(y,f) \;+\; [R_\kappa\kappa\rho\, Q^{\tt r}H^{\tt r}](y,f) 
\end{equation}
which allows for iteration. By (\ref{iterate-2}) and \eqref{iterate-2a}, the expected occupation measure \eqref{eq:Def_H} has the following interpretation: 
At rate $\kappa$, we erase unit mass travelling along the trajectory of $\xi$, replace it by mass $\rho$ (generating $k$ particles with probability $p_k$, and then merging these $k$ particles), then shift the location $\eta^{\tt r}_{\tau^-}$ of the merged mass to a random position $w$ selected according to $Q^{\tt r}(\eta^{\tt r}_{\tau^-},dw)$. 
The underlying strongly Markovian system (again defined probabilistically by killing and repasting since the corresponding semigroup is i.g.\ not contractive) has the generator  
\begin{align}\label{iterate-2b}\begin{aligned}
&\call f (z) \;-\;  \kappa(z) f (z) \;+\; \kappa(z)\rho(z)\int_E Q^{\tt r}(z, dw)\, f(w) \\
&=\quad \call f (z) \;-\; [\kappa(1-\rho)](z) f (z) \;+\; \kappa(z)\rho(z)\int_E Q^{\tt r}(z, dw) \left[ f(w)-f(z) \right]  \\
&=\quad \wt\call f (z) \;-\; [\kappa(1-\rho)](z) f (z) 
\end{aligned}\end{align}
(notations from (\ref{generator-1}), (\ref{generator-2}), \eqref{defn:Q_r}) and is thus identified as the jump diffusion $\wt\xi$ on $E$ `killed' at position-dependent rate $z\to [\kappa(1-\rho)](z)$ 
(of course, since we do not assume $\rho\le 1$, speaking of `killing' is abuse of language). Iteration of (\ref{iterate-2a}) combined with \eqref{iterate-2b} then provides us with the following explicit solution to (\ref{eq:Def_H}): 
\begin{align}\begin{aligned}\label{iterate-3} 
H^{\tt r}(y,f) & \;\;=\; \;\sum_{n\in\mathbb{N}_0}\left[\,(R_\kappa\kappa\rho\, Q^{\tt r})^nR_\kappa\,\right](y,f)\\
&\;\;=\;\; E_y\left( \int_0^\infty dt\, f(\wt\xi_t)\,  e^{- \int_0^t [\kappa(1-\rho)](\wt\xi_s)\, ds } \right)   \;\;\le\;\; \infty \;. 
\end{aligned}\end{align}  
This is the $\kappa(1{-}\rho)$-resolvent kernel of the jump diffusion $\wt\xi$, well-defined since $f\ge 0$, and finite for bounded $f$ by (\ref{erg-1}) in assumption 1.2.2. We now take $f\equiv 1$ in (\ref{iterate-3}) and \eqref{eq:Def_H}. \halmos\\

Compare the last proof to  \cite{HL-05}, lemma 1.4 and its proof, and to \cite{Ha-12}, Prop.\ 3.2.21 and Cor.\ 3.2.30, and note the role of the kernel $K_2(\cdot,\cdot)$ from {\bf (A2)} iii) which scatters offspring produced at branching events:
Indeed, the kernel defined in \cite{Ha-12}, (3.2.23) corresponds exactly to our definition of $Q^{\tt r}$ in \eqref{defn:Q_r}.
\\

\begin{rem}\label{1.2_remark}
Following definition 4.10 in Ikeda, Nagasawa and Watanabe \cite{INW-69}, the  semigroup $(M_t)_{t\ge0}$ on the single-particle space $(E,\mathcal{E})$  
\beqq\label{expectation_semigroup_1}
M_t(y,f) \;:=\; E_y\left( {\ov f} ( \eta^{\tt r}_t) \right), \qquad t\ge 0 \;,\; y\in E\;,\;f:E\to[0,\infty) \;\mbox{measurable} 
\eeqq
is called  \emph{expectation semigroup} for the branching diffusion without immigration $\eta^{\tt r}$.  
In case $f=1_A$,  $\,M_t(y,A) = E_y\left(\eta_t^{\tt r}(A)\right)\,$ is the expected number of particles visiting $A$ at time $t$ which descend from a single ancestor in $y$ at time $0$. 
This semigroup 
was implicit in definition \ref{1.2.3}, via $H^{\tt r}(y,f)=\int_0^\infty M_t(y,f)\,dt$. 
In analogy with the derivation (\ref{iterate-1})-(\ref{iterate-3}), one can use the strong Markov property and the branching property to obtain by iteration the representation
\beqq\label{expectation_semigroup}
M_t(y,f) \;=\; E_y\left( f(\wt\xi_r)\; e^{-\int_0^t [\kappa(1-\rho)](\wt\xi_v) dv } \right) 
\eeqq 
which identifies the expectation semigroup as 
the Feynman-Kac semigroup corresponding to the jump diffusion $\wt\xi$ with generator \eqref{generator-2} and the `potential' $\kappa(1-\rho)$.
We refer to \cite{Ha-12}, Thm.\ 3.2.28 for a full proof under our present assumptions.

Identities of the form \eqref{expectation_semigroup_1}+\eqref{expectation_semigroup} expressing the expected number of particles in terms of the dynamics of a single particle have a long history, going back to (at least) Watanabe \cite{Wa-67}.
They are now commonly called `many-to-one'-formulas (see e.g.\ \cite{HH-09}) 
and have been generalized in various ways; in particular, in \eqref{expectation_semigroup_1} the function $\ov{f}(\eta_t^{\tt r})$ may be replaced by a functional depending on the whole path of the process up to time $t$.
However, most of this literature tends to focus on the case of local branching mechanisms where particles reproduce exactly at their death position. 
See \cite{BDMT-11} and \cite{Ma-18} for versions admitting non-local branching, also employing an auxiliary process as our jump diffusion $\widetilde\xi$, but still under stronger conditions on the offspring mechanism than our assumption \textbf{(A2)} iii)
(\cite{BDMT-11} assumes in addition constant rates).
\\
\end{rem}

Now we can prove lemma \ref{1.2.5}.\\

{\bf Proof of lemma \ref{1.2.5}:} 1) By lemma \ref{1.2.4} (where $\ell(x)\ge 1$ for $x\neq\delta$) and in virtue of (\ref{erg-2}) in assumption \ref{1.2.2}, the expected time to extinction of a subprocess $\eta^{{\tt r},j}$ of $\eta$ defined by all descendents of the $j$-th immigrant is finite, the $j$-th immigrant choosing its location according to $Q^{\tt i}$ by {\bf (A3)}. Since the process of immigration instants is a Poisson random measure with constant intensity $c$ on $(0,\infty)$, the BDI process $\eta$ will a.s.\ in the long run return infinitely often to the void configuration $\delta$.  

2) By 1), the BDI process $\eta$ can be rewritten in the form of a sum of i.i.d.\ excursions away from the void configuration $\delta$. 
Write $R_1, R_2, \ldots$ for the times of successive returns to $\delta$, and define a measure on the configuration space $(S,\cals)$ by 
\beqq\label{inv-1}
\mu(F) \;:=\; E_\delta\left( \int_0^{R_1} 1_F(\eta_s)\, ds \right) \quad,\quad F\in\cals \;. 
\eeqq
Sets $F\in\cals$ of positive $\mu$-measure are visited infinitely often in the long run, a.s.\ for every choice of a starting point in $S$. 
Thus $\eta$ is a Harris process (we refer to \cite{ADR-69} and \cite{Re-84}, \cite{Nu-78}, \cite{Nu-85} for Harris recurrence). A Harris process has a unique (up to constant multiples) invariant measure which for the moment we may call $\check\mu$. We know that $\check\mu$ is equivalent to $\mu$, and we have ratio limit theorems: for pairs of measurable functions $f,g:S\to[0,\infty)$, with $g$ such that $0<\check\mu(g)<\infty$, the limits
$$
\lim_{t\to\infty}\frac{\int_0^t f(\eta_s)\, ds}{\int_0^t g(\eta_s)\, ds} \;\;=\;\; \frac{\check\mu(f)}{\check\mu(g)}
$$
exist almost surely, for every choice of a starting point $x\in S$. The structure of $\eta$ as a sum of i.i.d.\ excursions away from $\delta$ then allows to identify the limits with 
$$
\lim_{\nto}\frac{\int_0^{R_n} f(\eta_s)\, ds}{\int_0^{R_n} g(\eta_s)\, ds} \;\;=\;\; \frac{\mu(f)}{\mu(g)} \;.
$$ 
Hence invariant measure $\check\mu$ equals $\mu$ defined in \eqref{inv-1}, up to constant multiples. This shows that $\mu$ defined in (\ref{inv-1}) is invariant for the BDI process $\eta$.

3) Associate to the invariant measure $\,\mu\,$  on the configuration space $(S,\cals)$ defined by \eqref{inv-1} an invariant occupation measure $\,\ov\mu\,$ on the single particle space $(E,\cale)$ via
\beqq\label{inv-2}
\ov\mu(f) \;:=\; \mu(\ov f) \;=\;
E_\delta\left( \int_0^{R_1} \ov f(\eta_s)\, ds \right)  
\quad,\quad f:E\to[0,\infty) \;\mbox{measurable}  \;. 
\eeqq
Let $(\tau_j^{\tt i})_{j\ge 1}$ denote the sequence of successive immigration times, write $\tau_j^{\tt d}$ 
for the time of extinction of the subprocess $(\eta^{\tt r,j}_t)_{t\ge 0}$ of all direct descendants of the ancestor who immigrated at time $\tau_j^{\tt i}$, then 
$$
\lim_{\nto}\frac{\int_0^{R_n} \ov f(\eta_s)\, ds}{\int_0^{R_n} \ov g(\eta_s)\, ds} \;\;=\;\; \frac{\mu(\ov f)}{\mu(\ov g)} \;\;=\;\; \frac{\ov\mu(f)}{\ov\mu(g)}  
$$ 
coincides --with $H^{\tt r}$ from \eqref{eq:Def_H} and $Q^{\tt i}$ from {\bf(A3)}-- with 
$$
\lim_{j\to\infty}\frac{\int_0^{\tau_j^{\tt d}} \ov{f}(\eta_s)ds }{\int_0^{\tau_j^{\tt d}} \ov{g}(\eta_s)ds }  \;=\;
\lim_{m\to\infty}\frac{ \sum_{j=1}^m \int_{\tau_j^{\tt i}}^{\tau_j^{\tt d}} \ov{f}(\eta^{\tt r,j}_s)ds }{ \sum_{j=1}^m \int_{\tau_j^{\tt i}}^{\tau_j^{\tt d}} \ov{g}(\eta^{\tt r,j}_s)ds }  \;=\; 
\frac{ E_{\delta}\left( \int_0^\infty \ov{f}(\eta^{\tt r,1}_s) ds \right) }{ E_{\delta}\left( \int_0^\infty \ov{g}(\eta^{\tt r,1}_s) ds \right) } \;=\; 
\frac{ [Q^{\tt i}H^{\tt r}](f) }{ [Q^{\tt i}H^{\tt r}](g) }
$$ 
in application of definition \ref{1.2.3}. 
This shows that $\ov\mu$ in \eqref{inv-2} equals $\,Q^{\tt i}H^{\tt r}\,$, up to some multiplicative constant. Combining \eqref{iterate-3} with \eqref{erg-2} in 
assumption \ref{1.2.2} we see that $\ov\mu(1)=\mu(\ov 1)=\mu(\ell)$ is finite. This implies $\,\mu(1)=E_\delta(R_1)<\infty\,$, and we have the assertion of the lemma up to choice of norming constants: $\,\mu$ on $(E,\cale)$ is a finite measure, thus we have positive Harris recurrence of the BDI process $\eta$ with finite invariant occupation measure. 

4) It remains to determine the constants. Define $\calj_n := \max\{ j : \tau^{\tt i}_j < R_n \}$ with notations of 3). 
As a consequence of {\bf (A3)} we have  almost surely as $\nto$
$$
\calj_n \;\sim\; c\, R_n \;\sim\; c\, E_\delta(R_1)\, n 
$$ 
for every choice of a starting point for the process $\eta$ (where $E_\delta(R_1)>\frac 1c$ shows that the right hand side is necessarily larger than $n$), together with 
$$ 
n\; \ov\mu(f)  \;\;\sim\;\; \int_0^{R_n} \ov f( \eta_s )\, ds   
\;=\; \sum_{j \le \calj_n} \int_{\tau_{j}^{\tt i}}^{\tau_{j}^{\tt d}} \ov f( \eta_s )\, ds 
\;\;\sim\;\; \calj_n\; [Q^{\tt i}H^{\tt r}](f)\;\;\sim\;\; n\; c\, E_\delta(R_1)\,[Q^{\tt i}H^{\tt r}](f) 
$$
almost surely as $\nto$. This establishes  
\beqq\label{the_correct_constant}  
\ov\mu(f) \;=\; c\, E_\delta(R_1)\,[Q^{\tt i}H^{\tt r}](f) 
\eeqq
when invariant measure is defined by \eqref{inv-1} and invariant occupation measure by \eqref{inv-2}. 
Now, dividing the right hand sides of \eqref{inv-1}+\eqref{erg-2} and both sides of \eqref{the_correct_constant} by $E_\delta(R_1)$ and changing notations correspondingly, we get the assertion of the lemma with respect to the invariant probability. 
\halmos\\

\section{Some properties of the invariant probability and the invariant occupation measure} 
\label{1.3}


We state and prove two theorems on the invariant measure and the invariant occupation measure. Both will be key tools in the statistical context of sections~\ref{reconstruction}
and \ref{regressionNEU}. Theorem \ref{1.3.7} deals with finite `moments' $\mu(\ell^q)$ of the invariant probability $\mu$ of the BDI process of the same order $q$ as the reproduction law in {\bf (A2)}~ii). Theorem  \ref{1.3.10} gives conditions which grant existence of a continuous Lebesgue density of the invariant occupation measure $\ov\mu$. The proofs are given in sections  \ref{5.3} and  \ref{HKB-proof_continous_densities}.

For the special case of local branching where particles reproduce exactly at their death position, the existence of a continuous invariant occupation density has been considered by H\"opfner and L\"ocherbach \cite{HL-05}; 
with different methods, L\"ocherbach \cite{Lo-04} and Hammer \cite{Ha-12} allow for interactions between  particles (see remark \ref{1.3.9} below). 
In our setting, due to the general form of the kernel in {\bf (A2)}~iii) which scatters offspring generated at a branching event relative to the parent's position, we take a different approach. 
\\

\subsection{ Two theorems }\label{2theorems}

We introduce further assumptions (not all of these will be in force at the same time) and strengthen preceding ones. 
From now on, \ref{1.2.1neo} and \ref{1.2.2} are always assumed, $\,\mu\,$ is the invariant probability of the BDI process $\eta$ on the configuration space $(S,\cals)$, and $\ov\mu$ the invariant occupation measure on the single particle space $(E,\cale)$ as specified by \eqref{explicit_form_mubar} in lemma \ref{1.2.5}. 
\\

\begin{assumption}\label{1.3.2neo} 
There is some natural number $q>1$ such that $y \to m_q(y)$ is bounded on $E$, where 
$$
m_q(y):=\sum_{k\in\bbn_0} k^q\, p_k(y) \le\infty
$$
denotes $q$-th moments of the position-dependent reproduction laws $(p_k(y))_k$ at $y\in E$ in {\bf (A2)}~ii). 
\end{assumption}

Our next assumption strengthens heavily (\ref{erg-2}) of assumption \ref{1.2.2}. Recall the expectation semigroup $(M_t)_{t\ge0}$ for the branching process without immigration $\,\eta^{\tt r}\,$ from \eqref{expectation_semigroup_1}, 
associated to the expected occupation measure \eqref{eq:Def_H}, and its 
representation as a Feynman-Kac semigroup in the `many-to-one'-formula \eqref{expectation_semigroup} in remark \ref{1.2_remark}. 
\\

\begin{assumption}\label{1.3.4}
With notation $||| M_t ||| := \sup\limits_{y\in E} M_t(y,E) = \sup\limits_{y\in E} E_y\left( e^{-\int_0^t [\kappa(1-\rho)](\wt\xi_v) dv } \right)$, we have  
\begin{equation}\label{exponential-stability}
 \,\limsup\limits_{t\to\infty}\, \frac1t\, \log\left( ||| M_t ||| \right) < 0  \;. 
\end{equation}
\end{assumption}
\vskip0.5cm

Assumption \ref{1.3.4} implies in particular that the function in \eqref{erg-1} is bounded, thus \eqref{erg-2} of assumption \ref{1.2.2} holds for any choice of an immigration measure $Q^{\tt i}$. Property \eqref{exponential-stability} is known in the general theory of semigroups as \emph{uniform exponential stability}. We refer to \cite{EN-00}, Ch.\ V, Sec.\ 1 for a number of equivalent characterizations which can be used to check our assumption \ref{1.3.4} whenever the semigroup $(M_t)_{t\ge0}$ is strongly continuous on the Banach space $\mathcal{C}_0(E)$ of continous functions vanishing at infinity.  
\\

\begin{thm}\label{1.3.7}
Under \ref{1.2.1neo}, \ref{1.3.2neo} and  \ref{1.3.4}, we have finite `moments' of the invariant measure 
$$
\mu(\ell^q) \;:=\; \int_S \ell^q(x)\;\mu(dx)  \;=\;
\sum_{\ell\in\bbn} \ell^q\, \mu(E^\ell)   \;<\; \infty
$$
where $q>1$ is specified by assumption \ref{1.3.2neo}. 
\end{thm}
\vskip0.5cm

Theorem \ref{1.3.7} will be proved in section \ref{5.3}. 
Our next assumption concerns the semigroup 
\beqq\label{1p_motion_killed_kappa}
P_t^\kappa(y,f)
\;:=\; E_y\left( f(\xi_t) \, e^
{ -\int_0^t \kappa(\xi_s)\, ds } \right) \qquad t\ge 0 \;,\; y\in E \;,\; f:E\to[0,\infty) \;\mbox{measurable}   
\eeqq
of the single-particle motion $\xi$ killed at rate $\kappa$.  
The semigroup \eqref{1p_motion_killed_kappa} was already implicit in the proof of lemma \ref{1.2.4}, see \eqref{defn:R_kappa}.
For this semigroup, 
we shall now require existence of heat kernel bounds (which Hammer \cite{Ha-12} 
used to investigate the invariant measure $\mu$ on $S$, see remark \ref{1.3.9} below).  
For sufficient conditions implying such bounds, we refer to Dynkin \cite{Dy-65} theorem 0.5  p.\ 229 appendix paragraph 6, or Friedman \cite{Fr-75} theorem 4.5 p.\ 141. 
\\

\begin{assumption}\label{1.3.8}
The semigroup in (\ref{1p_motion_killed_kappa}) admits densities 
$p_t^\kappa(y,z)\, dz = P_t^\kappa(y,dz)$ with respect to Lebesgue measure which are continuous in $z$ for fixed $y$ and admit bounds 
\begin{equation}\label{eq:HKB}
p_t^\kappa(y,z) \;\le\; C\, t^{-d/2}\, e^{ -\frac12 \frac{|z-y|^2}{C\, t}} \quad\mbox{for all $0<t\le t_0$, $\;y,z\in E$}   
\end{equation}
for some $t_0>0$ fixed and some positive constant $C$.   
\end{assumption}

Heat kernel bounds \ref{1.3.8} will be a key tool in our proof for the existence of a continuous invariant occupation density, as well as for the results in section \ref{reconstruction} below. We stress that \ref{1.3.8} is a strong assumption: even with $d=1$ and constant killing rate $\kappa\equiv 1$ it does not hold for Ornstein-Uhlenbeck one-particle motion $d\xi_t = -\vth \xi_t dt + dW_t$ when the OU parameter $\vth$ is different from $0$. On the other hand, by Dynkin \cite{Dy-65} p.\ 229,  assumption \ref{1.3.8} does hold for all choices of a H\"older continuous and bounded killing rate $\kappa$ whenever the single-particle motion $\xi$ is such that uniform ellipticity holds on $E$ and all $|b^i|$, $|\si^{i,j}|$ in {\bf (A1)} are bounded.  
Our final assumption requires that the transition probability $Q^{\tt r}(\cdot,\cdot)$ of \eqref{defn:Q_r} admits bounds of convolution type. 
\\

\begin{assumption}\label{1.3.7neo}
There exists some finite measure $\widehat Q^{\tt r}$ on the single-particle space $(E,\mathcal{E})$ such that 
\begin{equation}\label{boundedness_Q}
Q^{\tt r}(y,A)\le \widehat Q^{\tt r}(A-y), \qquad y\in E, \; A\in \cale \;. 
\end{equation}
\end{assumption}
\vskip0.5cm

Note that \eqref{boundedness_Q} is essentially a condition on the transition probability $K_2(\cdot,\cdot)$ of {\bf (A2)}~iii). 
Clearly assumption \ref{1.3.7neo} covers the case of a product structure \eqref{eq:K_special} where $K(y,dv)=q(dv)$ for some probability measure $q$ on $E$: here we take  $\widehat Q^{\tt r}:=q$ and have equality in \eqref{boundedness_Q}. It also covers the case of absolutely continuous product structures
\begin{equation}\label{abs_cont_offspring_1}
 K_2((y,k),dv_1,\ldots,dv_k) = \prod_{j=1}^k q((y,k),v_j)\,\nu(dv_j)
 \end{equation}
for $\sigma$-finite measures $\nu$ on $E$ when $\nu$-densities depend on $y$ and $k$  but are uniformly dominated by
\begin{equation}\label{abs_cont_offspring_2}
q((y,k),v)\le\widehat q(v),\qquad y\in E,\; k\in\mathbb{N},\; v\in E 
\end{equation}
where $\widehat q\in L^1(\nu)$; then \eqref{boundedness_Q} holds for $\widehat Q^{\tt r}(A):=\int_{A}\widehat q(v)\,\nu(dv)$.  
Note that we do not require the $\sigma$-finite measure $\nu$ to be Lebesgue-absolutely continuous. 
Beyond  \eqref{abs_cont_offspring_1} and \eqref{abs_cont_offspring_2}, we see from \eqref{defn:Q_r} that assumption \ref{1.3.7neo} controls in some sense the distance of a `typical' child from its parent's position.

The following is the second main probabilistic result: heat kernel bounds \ref{1.3.8} for particle motion killed at rate $\kappa$ and convolution bounds \ref{1.3.7neo} on the scattering of offspring at branching events allow to obtain a continuous Lebesgue density for the invariant occupation measure. 
Theorem \ref{1.3.10} will be proved in section \ref{HKB-proof_continous_densities}. \\

\begin{thm}\label{1.3.10}
Assume \ref{1.2.1neo}, \ref{1.3.8}, \ref{1.3.7neo}, and suppose that the immigration measure $Q^{\tt i}$ is such that  condition \eqref{erg-2} of \ref{1.2.2} is satisfied.  
If $d\ge2$, suppose in addition that $Q^{\tt i}(dx)=q^{\tt i}(x)dx$ is absolutely continuous with Lebesgue density $q^{\tt i}\in L^p(\mathbb{R}^d)$ for some $p\in(\frac{d}{2},\infty]$.
Then the invariant occupation measure $\ov\mu$ on $(E,\cale)$ (a finite measure by lemma \ref{1.2.5}) admits a continuous Lebesgue density $\,\ov\gamma \in \mathcal{C}_0(E)\,$. 
\end{thm}

\vskip0.8cm
\begin{rem}\label{1.3.9} 
\quad i) 
H\"opfner and L\"ocherbach \cite{HL-05} proved existence of a continuous Lebesgue density for $\ov\mu$ in the special case of local branching,
i.e.\ when $K_2(\cdot,\cdot)$ of {\bf (A2)}~iii) is of product type \eqref{eq:K_special} with $K(y,dv)=\epsilon_0(dv)$. 
Their approach, using stochastic flows of diffeomorphisms, is not directly applicable in our case of non-local branching where we allow for jumps in the distribution of newborn particles, reflected in the jump diffusion $\widetilde\xi$ with generator \eqref{generator-2}. 
However, it can be adapted to our setting by using duality theory for (Feller) semigroups. This approach, which will be taken up in another paper, leads to a continuous invariant occupation density under an alternative set of conditions on the single particle motion and the branching and reproduction mechanism. 
In the present work however, we restrict to the setting of assumptions \ref{1.3.8} and \ref{1.3.7neo}, since the heat kernel bounds \eqref{eq:HKB} will also be used (independently) in the proofs of our results in section \ref{reconstruction} below. \\
ii) 
For the case of local and binary branching, L\"ocherbach \cite{Lo-04} considered a generalization of the model where coexisting particles move as \emph{interacting} diffusions. 
In a $\calc^\infty_b$-setting, assuming uniform ellipticity, Malliavin calculus establishes the existence of a continuous invariant occupation density (theorem 4.2 in \cite{Lo-04}). 
\\
iii) 
Assuming existence of Lebesgue densities $q^{\tt i}$ for $Q^{\tt i}$ and of transition densities for $K_2(\cdot,\cdot)$ as in \eqref{abs_cont_offspring_1}-\eqref{abs_cont_offspring_2} such that the Fourier transforms of $q^{\tt i}$ and $\widehat q$ are integrable, Hammer \cite{Ha-12} used Fourier methods to deduce existence of a continuous Lebesgue density of the invariant measure $\mu$ on the configuration space $S$ from the heat kernel bound assumption \ref{1.3.8} 
(see assumptions 2.2.1, 2.2.5 and theorem 2.2.8 in \cite{Ha-12}), where
continuity on $S$ is understood layer-wise, i.e.\ for every $\ell$ the restiction $\mu(\cdot\cap E^\ell)$ of $\mu$ to $E^\ell$ admits a Lebesgue density $\gamma^\ell$ which belongs to $\calc_0(E^\ell)$. 
We shall not make use of this result in the present paper. 
\end{rem}
\vskip0.8cm

%
\subsection{Proof of theorem \ref{1.3.7}}\label{5.3}

This subsection is devoted to the proof of theorem \ref{1.3.7}. 
Recall that for a measurable function $f: E\to\mathbb{R}$ we write $\ov f:S\to\mathbb{R}$ for the function $\ov f(x):=\sum_{j=1}^\ell f(x_j)$, $x=(x_1,\ldots,x_\ell)\in S\,$, with $\ov f(\delta)=0$. As in definition \ref{1.2.3}, $\,\eta^{\tt r}\,$ is the branching diffusion without immigration. Let $\,(T^{\tt r}_t)_{t\ge 0}$ denote the semigroup of $\eta^{\tt r}$  
\begin{equation}\label{branching-semigroup}
T_t^{\tt r}(x,g):=E_x(g(\eta_t^{\tt r})),\qquad t\ge0, \; x\in S,\; g:S\to[0,\infty)\text{ measurable}   
\end{equation}
which is related to the expectation semigroup $(M_t)_{t\ge0}$ introduced in \eqref{expectation_semigroup_1} by 
\[
M_t(y,f)=T_t^{\tt r}(y,\ov f),\qquad t\ge 0 \;,\; y\in E\;,\;f:E\to[0,\infty) \;\mbox{measurable}.
\]
Moreover, let $(T_t)_{t\ge0}$ denote the semigroup of the BDI process $\eta$ on $S$.\\ 

We start with the branching diffusion without immigration $\,\eta^{\tt r}\,$ and study `higher moments' $\,T_t^{\tt r}(y, {\ov f}^p )\,$ for $p>1\,$ when $y$ ranges over the single-particle space $E$.
The following is Ikeda, Nagasawa and Watanabe \cite{INW-69}, (4.97) in theorem 4.15 on p.\ 144: \\

\begin{lemma}\label{2.2.1_neufassung}
(\cite{INW-69}) We have a representation 
\begin{align}\label{recursion-p-th-moment}\begin{aligned}
T_t^{\tt r}(x,\ov f^p)&=M_t(x,f^p) + \int_0^t ds\int_E M_{t-s}(x,dy)\kappa(y)\sum_{n\ge2}p_n(y)\times\\
&\qquad\times\sum_{\substack{(k_1,\ldots,k_n):\,0\le k_j<p,\\ k_1+\cdots+k_n=p}}\binom{p}{k_1,\ldots,k_n}\int_{E^n}K_2((y,n),dv_1,\ldots,dv_n)\prod_{j=1}^n T_s^{\tt r}(y+v_j,\ov f^{k_j})
\end{aligned}\end{align}
for $x\in E$, $f:E\to[0,\infty)$ bounded measurable, $p\in\bbn$. 
\end{lemma}

{\bf Sketch of Proof:} 

First, we note that the expectation semigroup \eqref{expectation_semigroup_1} has the series representation
\begin{align}\label{equation M_t 3}\begin{aligned}
M_t(x,f)=P_t^\kappa(x,f)+&\sum_{m\in\mathbb{N}}\int_0^tds_1\int_E[P_{s_1}^\kappa\kappa\varrho Q^{\tt r}](x,dy_1)\int_0^{t-s_1}ds_2\int_E[P_{s_2}^\kappa\kappa\varrho Q^{\tt r}](y_1,dy_2)\cdots\\
&\cdots\int_0^{t-s_1-\ldots-s_{m-1}}ds_m\int_E[P_{s_m}^\kappa\kappa\varrho Q^{\tt r}](y_{m-1},dy_m)P^\kappa_{t-s_1-\ldots-s_m}(y_m,f)
\end{aligned}\end{align}
where $(P_t^\kappa)_{t\ge0}$ and $Q^{\tt r}$ are defined in \eqref{1p_motion_killed_kappa} and \eqref{defn:Q_r}, respectively (see e.g.\ \cite{Ha-12} lemma 3.2.20).

Now proceeding as in \cite{INW-69}, 
we take $h\equiv1$ in their lemma 4.8, eq. (4.75) on p.\ 139 to obtain
\begin{align}
 \label{eq:0}\begin{aligned}
&T_t^{\tt r}(x,\bar f^p) = P_t^\kappa(x,f^p) + \int_0^tds\int_E [P_s^\kappa\kappa](x, dy)\sum_{n\in\mathbb{N}} p_n(y)\int_{E^n}K_2((y,n),dv_1,\ldots, dv_n)\times\\
&\qquad\qquad\qquad\qquad\qquad\qquad\times\sum_{\substack{(k_1,\ldots,k_n):\,0\le k_j\le p,\\ k_1+\cdots+k_n=p}}\binom{p}{k_1,\ldots,k_n}\prod_{j=1}^nT_{t-s}^{\tt r}(y+v_j,\bar f^{k_j})
\end{aligned}\end{align}
for each $x\in E$, where we have adjusted their notation to ours. 
(Essentially, this formula is obtained by conditioning on the first branching event and using the branching property.)

Now we decompose the sums arising in \eqref{eq:0}
\beqq\label{local_mark_aux1}
\sum_{\substack{(k_1,\ldots,k_n):\,0\le k_j\le p,\\ k_1+\cdots+k_n=p}}\binom{p}{k_1,\ldots,k_n} \prod_{j=1}^n T_{t-s}^{\tt r}(y+v_j,\ov f^{k_j})
\eeqq
into two terms. The first one collects all indices where $0\le k_j< p$ for all $j=1,\ldots,n$. 
The remaining second term, collecting indices of type $\,(k_1,\ldots,k_n)=p\, {\tt e}_j$ where ${\tt e}_j$ is the $j$-th unit vector in $\bbr^n$, shrinks to     
\beqq\label{local_mark_aux2}
\sum_{j=1}^n T_{t-s}^{\tt r}(y+v_j,\ov f^p)    
\eeqq
where the maximal power $p$ shows up. 
Both contributions have to be integrated with respect to the kernel 
$$
\sum_{n\in\mathbb{N} } p_n(y)\, \int_{E^n}K_2((y,n),dv_1,\ldots,dv_n).
$$

Using notation \eqref{defn:Q_r} and defining
\[J_p(s;y):=\sum_{n=2}^\infty p_n(y)\int_{E^n}K_2((y,n),dv_1,\ldots, dv_n)\sum_{\substack{(k_1,\ldots,k_n):\,0\le k_j<p,\\ k_1+\cdots+k_n=p}}\binom{p}{k_1,\ldots,k_n}\prod_{j=1}^nT_s^{\tt r}(y+v_j,\bar f^{k_j})\]
for $s\ge0$ and $y\in E$, this gives
\begin{align}\label{eq:2}
 \begin{aligned}
T_t^{\tt r}(x,\bar f^p)&= P_t^\kappa(x,f^p) + \int_0^t ds\int_E [P_s^\kappa\kappa](x, dy)\,J_p(t-s;y)\\
&\qquad+ \int_0^tds\int_E [P_s^\kappa\kappa\varrho Q^{\tt r}](x, dy)\,T_{t-s}^{\tt r}(y,\bar f^p).
\end{aligned}
\end{align}

The structure of the previous display (namely the occurence of $T_{t-s}^{\tt r}(y,\bar f^p)$ on the right hand side) allows for iteration: Expanding the last term on the right hand side of \eqref{eq:2} leads to
\begin{align*}\label{eq:3}
 \begin{aligned}
&T_t^{\tt r}(x,\bar f^p)\\
&= P_t^\kappa(x,f^p) + \int_0^t ds\int_E [P_s^\kappa\kappa](x, dy)\,J_p(t-s;y)\\
&+\sum_{m\in\mathbb{N}} \int_0^tds_1\int_E [P_{s_1}^\kappa\kappa\varrho Q^{\tt r}](x, dy_1)\cdots\int_0^{t-s_1-\cdots-s_{m-1}}ds_m\int_E [P^\kappa_{s_m}\kappa\varrho Q^{\tt r}](y_{m-1},dy_m)\,P_{t-s_1-\cdots-s_m}^\kappa(y_m, f^p)\\
&+ \sum_{m\in\mathbb{N}}\int_0^tds_1\int_E [P_{s_1}^\kappa\kappa\varrho Q^{\tt r}](x, dy_1)\cdots\int_0^{t-s_1-\cdots-s_{m-1}}ds_m\int_E[P_{s_m}^\kappa\kappa\varrho Q^{\tt r}](y_{m-1},dy_m) \times\\
&\qquad\qquad\qquad\times \int_0^{t-s_1-\cdots -s_m}ds_{m+1}\int_E[P_{s_{m+1}}^\kappa\kappa](y_{m},dy_{m+1})\,J_p(t-s_1-\cdots-s_{m+1};y_{m+1}).
\end{aligned}\end{align*}

Using the series representation \eqref{equation M_t 3} of the expectation semigroup, we see that the previous display can be transformed into \[\begin{aligned}
T_t^{\tt r}(x,\bar f^p)&=M_t(x,f^p) + \int_0^t ds\int_E M_s(x,dy)\,\kappa(y)\,J_p(t-s;y),
\end{aligned}\]
proving the representation \eqref{recursion-p-th-moment}.

\halmos\\

\begin{lemma}\label{5.3.6}
Assume \ref{1.2.1neo}, \ref{1.3.2neo} and \ref{1.3.4}.
Fix a natural number  $q>1$ such that \ref{1.3.2neo} holds. Then there exist $\gamma>0$ and constants $C_1,\ldots,C_q$ such that for  $f:E\to[0,\infty)$ bounded and measurable  
\beqq\label{upper-bound-semigroup-2}
\big\|T_t^{\tt r}(\ov f^p)\big|_E\big\|_\infty=\sup_{x\in E}|T_t^{\tt r}(x,\ov f^p)|\le C_p \,e^{-\gamma t}\,\|f\|_\infty^p \;,\qquad t\ge0 \;,\; p=1,\ldots,q  
\eeqq
where $\,\ldots\big|_E$ denotes the restriction of the kernel \eqref{branching-semigroup} to the single-particle space $E$. 
\end{lemma}

\vskip0.5cm
{\bf Proof:} 
The case $p=1$ is assumption \ref{1.3.4}: we know that there exist $C>0$ and $\gamma>0$ such that
\begin{equation}\label{spectral-radius-3}
|||M_t|||\le Ce^{-\gamma t},\qquad t\ge0.
\end{equation}
We proceed by induction: let 
$1<p\le q$ and assume that \eqref{upper-bound-semigroup-2} already holds for $k=1,2,\ldots,p-1$, i.e. there are constants $C_1,C_2,\ldots,C_{p-1}$ such that
\[
\big\|T_t^{\tt r}(\ov f^k)\big|_E\big\|_\infty\le C_k \,e^{-\gamma t}\,\|f\|_\infty^k,\qquad t>0,\;k=1,\ldots,p-1.
\]
For indices $(k_1,\ldots,k_n)$ appearing in the sum in \eqref{recursion-p-th-moment}, we have $0\le k_j\le p-1$ and $k_1+\cdots+k_n=p$, thus $k_j\ge1$ for at most $p$ and at least $2$ indices $j$. 
Then by induction  
\[\begin{aligned}
\prod_{j=1}^n\big\|T_s^{\tt r}(\ov f^{k_j})\big|_E\big\|_\infty&=\prod_{j:k_j\ge1}\big\|T_s^{\tt r}(\ov f^{k_j})\big|_E\big\|_\infty \;\le\;
\prod_{j:k_j\ge1}C_{k_j}\,e^{-\gamma s}\,\|f\|_\infty^{k_j}\\
&\le\quad e^{-2\gamma s}\prod_{j:k_j\ge1}C_{k_j}\prod_{j:k_j\ge1}\|f\|_\infty^{k_j} \quad\le\quad 
C'_p\,e^{-2\gamma s}\,\|f\|_\infty^{p} 
\end{aligned}\]
where we define $C'_p:=\left(\max\{1,C_1,\ldots,C_{p-1}\}\right)^p$. Substituting this into \eqref{recursion-p-th-moment} 
and making use of assumption \ref{1.3.2neo} --and again of \ref{1.3.4}-- 
we obtain
\[\begin{aligned}
&\big\|T_t^{\tt r}(\ov f^p)\big|_E\big\|_\infty \\ 
&\le |||M_t|||\,\|f\|_\infty^p+\|\kappa\|_\infty\int_0^tds\,|||M_{t-s}|||\,\bigg\|\sum_{n\ge2}p_n(\cdot)\sum_{\substack{(k_1,\ldots,k_n):\,0\le k_j<p,\\ k_1+\cdots+k_n=p}}\binom{p}{k_1,\ldots,k_n}\bigg\|_\infty\,C'_p\,e^{-2\gamma s}\,\|f\|_\infty^{p}\\
&\le Ce^{-\gamma t}\,\|f\|_\infty^p+\|\kappa\|_\infty\int_0^tds\,Ce^{-\gamma(t-s)}\,C'_p\,e^{-2\gamma s}\,\bigg\|\sum_{n\ge2}n^p\,p_n(\cdot)\bigg\|_\infty \|f\|_\infty^{p}\\
&\le Ce^{-\gamma t}\,\|f\|_\infty^p\left(1+C'_p\,\|\kappa\|_\infty\|m_p(\cdot)\|_\infty\int_0^te^{-\gamma s}\,ds\right)\\
&\le C_p\,e^{-\gamma t}\,\|f\|_\infty^p,
\end{aligned}\]
with $C_p:=C\left(1+\frac{C_p'\|\kappa\|_\infty\| m_p(\cdot)\|_\infty}{\gamma}\right)$.
Thus \eqref{upper-bound-semigroup-2} is proved. \halmos\\

We turn to the semigroup $(T_t)_{t\ge 0}$ of the branching diffusion $\,(\eta_t)_{t\ge 0}\,$ and focus on  the void configuration $\delta$ as starting point at time $t=0$.\\

\begin{lemma}\label{5.3.2} 
Fix $p\in\mathbb{N}$, $t>0$ and consider $f:E\to[0,\infty)$ measurable. Then
\begin{align}\label{recursion}\begin{aligned}
T_t(\delta,\ov f^p) \;=\; 
c\sum_{k=0}^{p-1}\binom{p}{k}\int_0^t [Q^{\tt i}T^{\tt r}_s](\ov f^{p-k})\,
T_s(\delta,\ov f^k)\, ds  \;. 
\end{aligned}\end{align}
\end{lemma}
\vskip0.5cm

{\bf Proof:} 
Write for short $\nu_t:= T_t(\delta,\cdot)$. Since immigration times are distributed according to Poisson random measure with constant intensity $c>0$, $\;\nu_t(\ov f^p)$ has the following explicit form
\[\begin{aligned}
&\nu_t(\ov f^p)\\
&=e^{-ct}\sum_{n=1}^\infty c^n\int_0^tds_1\int_0^{s_1}ds_2\cdots\int_0^{s_{n-1}}ds_n\int_S[Q^{\tt i}T^{\tt r}_{s_1}](dz_1)\cdots\int_S[Q^{\tt i}T^{\tt r}_{s_n}](dz_n)\,\left(\ov f(z_1)+\cdots+\ov f( z_n)\right)^p\\
&=e^{-ct}\sum_{n=1}^\infty c^n \sum_{\substack{(k_1,\ldots,k_n):\,k_j\ge0,\\k_1+\cdots+k_n=p}}\binom{p}{k_1,\ldots,k_n}\int_0^tds_1\int_0^{s_1}ds_2\cdots\int_0^{s_{n-1}}ds_n\prod_{j=1}^n [Q^{\tt i}T^{\tt r}_{s_j}](\ov f^{k_j}) \;. 
\end{aligned}\]
The right hand side in the previous display can be simplified: define 
\[\wt\nu_t(\ov f^p):=e^{ct}\nu_t(\ov f^p).\]
Differentiating with respect to $t$ (and sorting the terms), we get
\[\begin{aligned}
\frac{d}{dt}\wt\nu_t(\ov f^p)
&=c\,\wt\nu_t(\ov f^p) + c\sum_{k=1}^p\binom{p}{k} [Q^{\tt i}T^{\tt r}_t](\ov f^k)\,\tilde\nu_t(\ov f^{p-k})\\
&=c\,\wt\nu_t(\ov f^p) + c\sum_{k=0}^{p-1}\binom{p}{k} [Q^{\tt i}T^{\tt r}_t](\ov f^{p-k})\,\tilde\nu_t(\ov f^{k})=:c\,\wt\nu_t(\ov f^p) + h(t).
\end{aligned}\]
Solving this linear inhomogenous ODE by variation of constants yields
\[\begin{aligned}
\wt\nu_t(\ov f^p)&=e^{ct}\int_0^te^{-cs}\,h(s)\,ds=e^{ct}\,c\sum_{k=0}^{p-1}\binom{p}{k}\int_0^t [Q^{\tt i}T^{\tt r}_s](\ov f^{p-k})\,e^{-cs}\wt\nu_s(\ov f^{k})\,ds.
\end{aligned}\]
Multiplying by $e^{-ct}$ again, we obtain
\eqref{recursion}.  \halmos\\

\begin{lemma}\label{sect_2.1_neu_1} 
Assume \ref{1.2.1neo}, \ref{1.3.2neo} and \ref{1.3.4}.
Consider a natural number  $q>1$ for which \ref{1.3.2neo} holds.
Then for $f:E\to[0,\infty)$ bounded and measurable,   
$$
\sup_{0<t<\infty} T_t(\delta,\ov f ^p) \;<\; \infty \quad,\quad 1\le p\le q \;.
$$
\end{lemma}
\vskip0.5cm
{\bf Proof:} 
It is sufficient to prove the assertion in case $f\equiv 1$: then $\ov f = \ov 1 =\ell$. Lemma \ref{5.3.2} allows for recursion. 
First, in case $p=1$, we combine $T_s(\delta,\ell^0)=1$ with  \eqref{recursion}  and lemma \ref{5.3.6}: 
$$
T_t(\delta,\ell^1) \;=\; c\, \int_0^t [Q^{\tt i} T^{\tt r}_s](\ell^1)\, ds \;\le\; c\, C_1\, \int_0^\infty e^{-\gamma s}  ds \;=\; \frac{c}{\gamma}\, C_1 \;=:\; M_1 \;<\; \infty \;. 
$$
Next, consider $p< q$. If our assertion holds for all $1\le k\le p$, with suitable bounds $M_k$, then it holds for $p{+}1$: by recursion \eqref{recursion},
$$
T_t(\delta,\ell^{p+1}) \;=\; c\, \sum_{k=0}^p  \binom{p+1}{k}  \int_0^t [Q^{\tt i} T^{\tt r}_s](\ell^{p+1-k})\, T_s(\delta,\ell^k)\, ds 
$$
is bounded  (we can apply \eqref{upper-bound-semigroup-2} to every $k$-th summand since $p+1\le q$) by 
$$
c\, \sum_{k=0}^p  \binom{p+1}{k}  \int_0^\infty [C_{p+1-k}\, e^{-\gamma s\,}][M_k] ds 
\;=\; \frac{c}{\gamma}\,\sum_{k=0}^p  \binom{p+1}{k}  C_{p+1-k}\, M_k \;=:\; M_{p+1}
$$
and we are done.\halmos\\

Now we have the tools to prove theorem \ref{1.3.7}. \\

{\bf Proof of theorem \ref{1.3.7}: } 
By the ergodic theorem for Harris recurrent processes (see e.g.\ 
\cite{ADR-69}, p.\ 30) we know that for all $\mu$-integrable functions $g:S\to\mathbb{R}$ 
\begin{equation}\label{ergodic theorem}
\lim_{t\to\infty}\frac{1}{t}\int_0^t T_s(x,g)\,ds=\mu(g)\qquad\text{for }\mu\text{-a.e. }x\in S.
\end{equation}
By a simple monotone convergence argument, \eqref{ergodic theorem} clearly extends to all nonnegative measurable $g:S\to[0,\infty)$ where the limit is equal to $+\infty$ if $g$ is not $\mu$-integrable. 
Moreover, \eqref{ergodic theorem} must in particular hold for $x=\delta$ since $\mu(\delta)>0$.
Choosing $g:=\ov{f}^q$, this gives
\begin{equation}\label{eq:cesaro-limit}\mu(\ov{f}^q)=\lim_{t\to\infty}\frac{1}{t}\int_0^t T_s(\delta,\ov{f}^q)\,ds\le\infty\end{equation}
for all $f:E\to[0,\infty)$ measurable and $q\in\mathbb{N}$.
But from \eqref{recursion}, we see that $t\mapsto T_t(\delta,\ov{f}^q)$ is increasing, thus the limit $\lim_{t\to\infty}T_t(\delta,\ov{f}^q)$ exists in $[0,\infty]$ and must be equal to the Ces\`aro limit \eqref{eq:cesaro-limit} 
and so
\begin{equation}\label{eq:limit-higher-moments}
\mu(\ov{f}^q)=\lim_{t\to\infty}T_t(\delta,\ov{f}^q)\le\infty\end{equation}
for all measurable $f:E\to[0,\infty)$ and $q\in\mathbb{N}$. 
Now let $q$ satisfy  assumption \ref{1.3.2neo}. Then 
we can use lemma \ref{sect_2.1_neu_1} to conclude that the limit \eqref{eq:limit-higher-moments} is finite 
for each bounded measurable $f$. 
Now the assertion of theorem \ref{1.3.7} follows by choosing $f\equiv 1$, i.e.\ $\ov f=\ell$.
\halmos\\

\begin{rem}\label{5.3.3}
Note that for the derivation of formula \eqref{eq:limit-higher-moments}, we did not use assumption \ref{1.3.2neo} nor did we need the full force of assumption \ref{1.3.4}, but only positive Harris recurrence of the BDI process $\eta$
(for which we know from lemma \ref{1.2.5} that e.g.\ the weaker assumption 1.2.2 is already sufficient).
In fact, the recursion \eqref{recursion} can be solved to obtain the following explicit formula 
\begin{equation}\label{eq:formula-q-moment}
\mu(\ov{f}^q)=
\sum_{n=1}^{q}\frac{c^{n}}{n!}\sum_{\substack{(k_1,\ldots,k_{n}):\,k_j\ge1,\\ k_1+\cdots+k_n=q}}\binom{q}{k_1,k_2,\ldots,k_n}\prod_{j=1}^n\left(\int_0^\infty [Q^{\tt i}T^{\tt r}_{s}](\ov f^{k_j})\,ds\right)\le\infty
\end{equation}
for each $q\in\mathbb{N}$ and $f:E\to[0,\infty)$ measurable. 
Our above proof shows that assumptions \ref{1.3.2neo} and \ref{1.3.4} together are sufficient to ensure finiteness of \eqref{eq:formula-q-moment} for bounded measurable $f$, but they are probably not necessary.
\end{rem}

\subsection{Proof of theorem \ref{1.3.10}}\label{HKB-proof_continous_densities}

We recall the occupation times kernel $H^{\tt r}$ for the branching diffusion process without immigration $\eta^{\tt r}$
defined in \eqref{eq:Def_H} with series representation \eqref{iterate-3} 
$$
H^{\tt r}(x,B)=\sum_{n\in\mathbb{N}_0}[(R_\kappa\kappa\varrho\, Q^{\tt r})^nR_\kappa](x,B),\qquad x\in E,\,B\subseteq E\text{ Borel},
$$
where 
$Q^{\tt r}$ is the kernel \eqref{defn:Q_r} and
\[R_\kappa(x,dy)=\int_0^\infty P_t^\kappa(x,dy)\,dt,\qquad x\in E\]
is the $\kappa$-resolvent of the single-particle diffusion $\,\xi\,$ from 
\eqref{defn:R_kappa} and \eqref{1p_motion_killed_kappa}. Assumption \ref{1.3.8} grants existence of a Lebesgue density 
\beqq\label{gleichung_vor_explicit_form_mubar_2}
R_\kappa(x,dz)=r_\kappa(x,z)\,dz\qquad\text{with}\qquad r_\kappa(x,z):=\int_0^\infty p_t^\kappa(x,z)\,dt \; ,\qquad x,z\in E \;. 
\eeqq
Thus, by lemma \ref{1.2.5} and \eqref{iterate-3}, the  invariant occupation measure   
\begin{equation}\label{explicit_form_mubar_2}
\ov\mu(B) = c\,  [Q^{\tt i}H^{\tt r}](B) = c \sum_{n\in\mathbb{N}_0}[Q^{\tt i}(R_\kappa\kappa\varrho\, Q^{\tt r})^nR_\kappa](B) \;<\;\infty \;,\qquad B\subseteq E\text{ Borel} 
\end{equation}
admits a Lebesgue density
\begin{equation}\label{eq:series}
z\;\mapsto\; c \sum_{n\in\mathbb{N}_0}\ov\gamma_n(z) \;,\quad \ov\gamma_n(z):= \int_E[Q^{\tt i}(R_\kappa\kappa\varrho\, Q^{\tt r})^n](dx)\,r_\kappa(x,z) \;. 
\end{equation}
We will show that under the assumptions of theorem \ref{1.3.10}, we have $\ov\gamma_n\in\mathcal{C}_0(E)$ for all $n\in\mathbb{N}_0$ and that the series in \eqref{eq:series} converges uniformly.\\

We fix $\vep>0$ and observe that by the semigroup property of $(P_t^\kappa)_{t\ge0}$ we can decompose
\[\begin{aligned}
R_\kappa(x,dy)
&=\int_0^\vep P_t^\kappa(x,dy)\,dt+\int_0^\infty P_{t+\vep}^\kappa(x,dy)\,dt\\
&=:R_{\kappa,\vep}(x,dy)+R_\kappa P_\vep^\kappa(x,dy),
\end{aligned}\]
where we define $R_{\kappa,\vep}(x,dy):=\int_0^\vep P_t^\kappa(x,dy)\,dt$. For the resolvent density, this means
\begin{equation}\label{decomp_density}
r_\kappa(x,z)=r_{\kappa,\vep}(x,z)+\int_ER_\kappa(x,dy)\,p_\vep^\kappa(y,z),
\end{equation}
with notation 
\[
r_{\kappa,\vep}(x,z):=\int_0^\vep p_t^\kappa(x,z)\,dt.
\] 
Now under assumption \ref{1.3.8}, we know that if we choose $\vep\le t_0$ and define a density $\tilde p_t(\cdot)$ as the right hand side in the heat kernel bound \eqref{eq:HKB}, then we have 
\begin{equation}\label{eq:bounds_1}
 p_t^\kappa(x,z)\le\tilde p_t(z-x):=C\, t^{-d/2}\, e^{ -\frac12 \frac{|z-x|^2}{C\, t}},\qquad 0<t\le\vep,\; x,z\in E
 \end{equation}
and consequently (by symmetry)
\begin{equation}\label{eq:bounds_2}
 r_{\kappa,\vep}(x,z)\le \tilde r_\vep(x-z):=\int_0^\vep\tilde p_t(x-z)\,dt,\qquad x,z\in E.
 \end{equation}
We denote by 
\begin{equation}
 \tilde P_t(x,dy):=\tilde p_t(x-y)\,dy, \qquad \tilde R_\vep(x,dy):=\tilde r_\vep(x-y)\,dy
\end{equation}
the corresponding convolution kernels. 
Regularity of the density $\tilde r_\vep$ depends heavily on the dimension: while $\tilde r_\vep$ is bounded in $d=1$, for $d\ge2$ it has a singularity at the origin. Independently of the dimension we have  
$\,\tilde r_\vep\in L^1(E)\,$ and $\,\tilde p_t\in L^1(E)\cap\mathcal{C}_0^\infty(E)\,$ for $0<t\le\vep$, so the kernels $\tilde R_\vep$ and $\tilde P_t$ induce bounded convolution operators on $L^\infty(E)$, $\mathcal{C}_0(E)$ and on $\mathcal{C}_b(E)$, and $\tilde P_t$ induces also a bounded convolution operator $L^1(E)\to\mathcal{C}_0(E)$.  
Moreover, by assumption \ref{1.3.7neo} the `jump operator' corresponding to the kernel $Q^{\tt r}$ is bounded by
\begin{equation}\label{boundedness-Q}
\int_E Q^{\tt r}(x,dy)\,f(y)
\;\le\; \int_E f(x+v)\,\widehat Q^{\tt r}(dv) \;=:\; \int_E f(y)\,\tilde Q(x,dy) \;, 
\quad x\in E,\; f\ge0 \text{ measurable},
\end{equation}
where $\tilde Q(x,dy)$ denotes the convolution kernel corresponding to the finite measure $\widehat Q^{\tt r}$ in assumption \ref{1.3.7neo}.
Since the kernels resp.\ operators $\,\tilde P_t\,$ for $0<t\le\vep$, $\,\tilde R_\vep\,$ and $\,\tilde Q\,$ are all convolution kernels resp.\ operators, they all commute with each other, a fact which we shall exploit heavily below. \\

\begin{lemma}\label{lemma:operators}
The $n$-fold convolution $\tilde r_\vep^{*n}$ of the density $\tilde r_\vep$ with itself has the property 
\[
\tilde r_\vep^{*n}\in\mathcal{C}_0(\mathbb{R}^d)\qquad\text{for }n>\frac{d}{2} \;; 
\]
in particular $\tilde r_\vep\in\mathcal{C}_0(\mathbb{R})$ for $d=n=1$. 
The following holds for $d\ge2$: the density $\tilde r_\vep(\cdot)$ is continuous on $\mathbb{R}^d\setminus\{0\}$ but has a singularity at the origin; we have  $\,\tilde r_\vep(\cdot) \in L^{p^*}(\mathbb{R}^d)$ for all $1\le p^*<\frac{d}{d-2}$ (where we understand $\frac{d}{d-2}=\infty$ for $d=2$).
\end{lemma}
\vskip0.5cm

\textbf{Proof: } 
1) The fact that $\tilde r_\vep^{*n}\in\mathcal{C}_0(\mathbb{R}^d)$ for $n>\frac{d}{2}$ is most easily seen by Fourier inversion: with $\mathcal{F}$ denoting the Fourier transform, we have
\[
\mathcal{F}[\tilde p_t](\xi)=C^{1+d/2}(2\pi)^{d/2}e^{-\frac{1}{2}Ct\|\xi\|^2},\qquad \xi\in\mathbb{R}^d,\,t>0 
\]
with the constant $C$ from 
assumption \ref{1.3.8}.  
We obtain 
\[\mathcal{F}[\tilde r_\vep](\xi)=\int_0^\vep\mathcal{F}[\tilde p_t](\xi)\,dt=2(2\pi C)^{d/2}\,\frac{1-e^{-\frac{1}{2}C\vep\|\xi\|^2}}{\|\xi\|^2}\le \widetilde C\,\frac{1-e^{-\frac{1}{2}\|\xi\|^2}}{\|\xi\|^2}=:h(\|\xi\|),\qquad\xi\in\mathbb{R}^d\setminus\{0\}.\]
Consequently,
\[
|\mathcal{F}[\tilde r_\vep^{*n}](\xi)|=|\left(\mathcal{F}[\tilde r_\vep](\xi)\right)^n|\le h(\|\xi\|)^n. 
\]
Integration in (hyper-)spherical coordinates shows that $h(\|\cdot\|)^n$ is integrable on $\mathbb{R}^d$ if $n>d/2$. 
This gives $\,\mathcal{F}[\tilde r_\vep^{*n}](\cdot)\in L^1(\mathbb{R}^d)\,$ for all such $n$, and thus $\,\tilde r_\vep^{*n}\in\mathcal{C}_0(\mathbb{R}^d)\,$ by Fourier inversion. 
The fact that $\tilde r_\vep$ is continuous on $\mathbb{R}^d\setminus\{0\}$ in any dimension is clear by dominated convergence.  

2) The following is from Hammer \cite{Ha-12}, (3.2.63) on p.\ 103: for $d\ge 2$ and $\la>0$, consider the $\la$-resolvent 
$$
\phi_\la(x) \;:=\; \int_0^\infty e^{-\la s}\,p_s(x)\,ds   \;,\; x\in\bbr^d
$$
of the heat flow, i.e.\ $p_s(x)$ is the density of the normal law $\caln(0,sI_d)$.   Then $x\to\phi_\la(x)$ is $p^*$-integrable on $\bbr^d$ if and only if $p^*<\frac{d}{d-2}$.  
This is seen as follows. Sato \cite{Sa-99}, formulae (30.28)+(30.29) on p.\ 204, gives an explicit representation of $\phi_\la$   
$$
\phi_\la(x) \;=\; {\tt cst}\; \|x\|^{-\frac{d-2}{2}}\, K_{\frac{d-2}{2}}(\sqrt{2\la}\|x\|)
$$
where $K_\nu$ denotes a modified Bessel function (for $K_\nu$, see \cite{Sa-99}, (4.9) on p.\ 21, and \cite{Fo-92}, p.\ 159) whose asymptotics at $0$ and at $\infty$  are known: when $r\downarrow 0$ we have $K_\nu(r)\sim {\tt cst}\, r^{-\nu}$ for $\nu>0$  and $K_0(r)\sim {\tt cst}\, \log(r)$; when $r\uparrow \infty$ we have exponential decay (see Follett \cite{Fo-92}, p.\ 160). It follows that $\,x\to\phi_\la(x)\,$ is $p^*$-integrable on $\bbr^d$ if and only if $p^*<\frac{d}{d-2}$. 

3) For the $L^{p^*}$-properties of $\tilde r_\vep$, we observe that
\[\begin{aligned}
 \tilde r_\vep(x)&=\int_0^\vep\tilde p_t(x)\,dt\le e^{\vep}\int_0^\vep e^{-t}\,\tilde p_t(x)\,dt\le C e^{\vep}\int_0^\vep e^{-t}\,t^{-d/2}e^{-\frac{1}{2}\frac{\|x\|^2}{Ct}}\,dt\\
 &=(2\pi C)^{d/2} e^\vep \int_0^{C\vep}e^{-s/C}\,(2\pi s)^{-d/2}e^{-\frac{\|x\|^2}{2s}}\,ds \;\le\; {\tt cst}\; \phi_{\frac{1}{C}}(x) 
\end{aligned}\]
for all $x\in\mathbb{R}^d$,  
with $\phi_{\frac{1}{C}}$ from step 2). So the last assertion of the lemma follows from step 2). 
\halmos\\

\begin{lemma}\label{lemma:C_0}
Under the assumptions of theorem \ref{1.3.10} and with notation $\ov\gamma_n$ from \eqref{eq:series}, we have $\ov\gamma_n\in\mathcal{C}_0(\mathbb{R}^d)$ for all $n\in\mathbb{N}_0$.
\end{lemma}
\vskip0.5cm

\textbf{Proof: } We use induction on $n\in\mathbb{N}_0$.

1) For $n=0$, definition \eqref{eq:series} of $\ov\gamma_0$ combined with  decomposition \eqref{decomp_density} gives
\begin{equation}\label{proof:C_0_1}
z\mapsto \ov\gamma_0(z)=\int_EQ^{\tt i}(dx)\,r_\kappa(x,z)=\int_EQ^{\tt i}(dx)\,r_{\kappa,\vep}(x,z)+\int_E[Q^{\tt i}R_\kappa](dx)\,p^\kappa_{\vep}(x,z).
\end{equation}
Fix $x\in E$. For $t$ sufficiently small, $z \to p^\kappa_t(x,z)$ is continuous by assumption \ref{1.3.8}.  Note first that regardless of the dimension of $E=\mathbb{R}^d$, the function 
\begin{equation}\label{proof:continuity-1}
z \mapsto \int_0^\vep p_t^\kappa(x,z)\,dt = r_{\kappa,\vep}(x,z) 
\end{equation}
is continuous at $z_0\in E$ whenever $z_0\ne x$. To see this, fix $z_0\ne x$ and consider a sequence $z_n\to z_0$; we may assume that there is $\delta>0$ such that $\|z_n-x\|>\delta$ for all $n$. Then the estimate \eqref{eq:HKB} gives 

$$
0 \;\le\; p_t^\kappa(x,z_n) \;\le\; C\, t^{-\frac d 2}\, e^{ - \frac{\delta^2}{2 C t } } \quad,\quad n\in\bbn \;,\; 0<t<\vep  
$$
for $\vep$ sufficiently small. Here the right hand side is independent of $n\in\bbn$ and integrable in $0<t<\vep$, thus dominated convergence shows 
$$
\int_0^\vep p_t^\kappa(x,z_n)\,dt \;\to\; \int_0^\vep p_t^\kappa(x,z_0)\,dt
$$
which establishes \eqref{proof:continuity-1}. Based on this we can check the assertions of the lemma in case $n=0$. We start with the function 
\beqq\label{its_first_term}
z \mapsto  \int_{\mathbb{R}} Q^{\tt i}(dx)\, r_{\kappa,\vep}(x,z)   
\eeqq 
on the right hand side of \eqref{proof:C_0_1}.  

i) 
In the special case $d=1$, the function $z\mapsto r_{\kappa,\vep}(x,z)$ in \eqref{proof:continuity-1} is continuous in $z\in\mathbb{R}$ for every $x\in \mathbb{R}$ fixed, and its upper bound $\tilde r_\vep(\cdot)$ from \eqref{eq:bounds_2} is in $\mathcal{C}_0(\mathbb{R})$ by lemma \ref{lemma:operators}, thus bounded. Thus dominated convergence shows that the function \eqref{its_first_term} is continuous and bounded,  
for any probability measure $Q^{\tt i}(dx)$ on ${\mathbb{R}}$.   
Probability measures on $\mathbb{R}$ being tight, upper bounds $z \mapsto \int_{\mathbb{R}} Q^{\tt i}(dx)\, \tilde r_\vep(z-x)$ with $\tilde r_\vep(\cdot)\in\calc_0(\mathbb{R})$ vanish at $\pm\infty$. This implies that the function in \eqref{its_first_term} belongs to $\mathcal{C}_0(\mathbb{R})$.

ii) Let $d\ge2$ and assume that $Q^{\tt i}(dx)=q^{\tt i}(x)dx$ is absolutely continuous with density in $L^p(\mathbb{R}^d)$ for some $p\in(\frac{d}{2},\infty]$.
Here, without loss of generality we may assume that $p<\infty$, since if $q^{\tt i}$ is bounded, then (being a probability density) it is in $L^1(\mathbb{R}^d)\cap L^\infty(\mathbb{R}^d)=\bigcap_{1\le p\le\infty}L^p(\mathbb{R}^d)$.
Then the dual exponent $p^*$ satisfies $1< p^*<\frac{d}{d-2}$. By lemma \ref{lemma:operators}, the function $\tilde r_\vep(\cdot)$ is in $L^{p^*}(\mathbb{R}^d)$, and is continuous on $\mathbb{R}^d\setminus\{0\}$. We have to consider
\begin{equation}\label{proof:continuity-2}
 z\mapsto \int_{\mathbb{R}^d}q^{\tt i}(x)\,r_{\kappa,\vep}(x,z)\,dx \le \int_{\mathbb{R}^d}q^{\tt i}(x)\,\tilde r_{\vep}(x-z)\,dx = [q^{\tt i}*\tilde r_\vep](z). 
\end{equation}
Since the convolution of two functions from dual $L^p$-spaces is 
in $\mathcal{C}_0({\mathbb{R}^d})$ for $p\in(1,\infty)$ (see e.g.\ \cite{LL-01}, Lemma 2.20, \cite{HS-65}, p.\ 398), the right hand side in \eqref{proof:continuity-2} is a $\mathcal{C}_0$-function of $z$.
For convergent sequences $z_n\to z_0$ and $x\neq z_0$ we have    
$$
r_{\kappa,\vep}(x,z_n) \to r_{\kappa,\vep}(x,z_0) \;\;,\;\;  \tilde r_{\vep}(x-z_n) \to \tilde r_{\vep}(x-z_0)  
$$
as $\nto$, using \eqref{proof:continuity-1} and continuity of $\tilde r_{\vep}(\cdot)$ on $\mathbb{R}^d\setminus\{0\}$ by lemma \ref{lemma:operators}.  Integrals on the right hand side of \eqref{proof:continuity-2} being continuous in $z$, Pratt's lemma applies (see e.g.\ \cite{Els-11}, theorem\ VI.5.1, \cite{Sc-05}, p.\ 101) and shows     
\[
\int_{\mathbb{R}^d}q^{\tt i}(x)\,r_{\kappa,\vep}(x,z_n)\,dx \;\lra\; \int_{\mathbb{R}^d}q^{\tt i}(x)\,r_{\kappa,\vep}(x,z_0)\, dx
\]
as $\nto$. We have proved that the function \eqref{its_first_term}
is continuous. Since its upper bounds in \eqref{proof:continuity-2} are in $\mathcal{C}_0({\mathbb{R}^d})$, the function \eqref{its_first_term} is in $\mathcal{C}_0({\mathbb{R}^d})$. 

iii)
So far we have shown that the first term \eqref{its_first_term} on the right hand side of \eqref{proof:C_0_1} is a $\mathcal{C}_0$-function of $z$. The second term on the right hand side of \eqref{proof:C_0_1}, as a function of $z$, is always in $\mathcal{C}_0(\mathbb{R}^d)$: 
by the regularity properties of $p_\vep^\kappa$ in assumption \ref{1.3.8} and since $[Q^{\tt i}R_\kappa](dx)$ is a finite measure, the argument is analogous to step i) above, and no regularity of the immigration measure is needed for this term. 

We have proved that in case $n=0$,  $\,z \to \ov\gamma_0(z)$ in \eqref{proof:C_0_1} has the property stated in the lemma.

2) To prove that all $\ov\gamma_n(\cdot)$ have the property stated in the lemma, we proceed by induction on $n$. Suppose we know already that $\ov\gamma_n\in\mathcal{C}_0(E)$. 
Then by \eqref{decomp_density} and \eqref{eq:series} 
\[\begin{aligned}
 \ov\gamma_{n+1}(z)&=\int_E[Q^{\tt i}(R_\kappa\kappa\varrho\, Q^{\tt r})^{n+1}](dx)\,r_\kappa(x,z)\\
 &=\int_E dx\,\ov\gamma_n(x)\,\kappa(x)\rho(x)\int_E Q^{\tt r}(x,dy)\,r_{\kappa,\vep}(y,z)+\int_Edx\,\ov\gamma_{n+1}(x)\,p_\vep^\kappa(x,z).
\end{aligned}\]
Again the second term on the right hand side is a $\calc_0$-function of $z$ because of the regularity of $p_\vep^\kappa$ and since $\ov\gamma_{n+1}\in L^1(E)$.
We consider the first term and put it using notations \eqref{boundedness-Q} in the form 
\beqq\label{xxx1_neo}
\begin{aligned}
z \;\lra\; &\int_E dx\,\ov\gamma_n(x)\,\kappa(x)\rho(x)\int_E\tilde Q(x,dy)\,q^{\tt r}(x,y)\,r_{\kappa,\vep}(y,z) \\
&=\int_Edx\int_E\widehat Q^{\tt r}(dv)\,\ov\gamma_n(x)\,\kappa(x)\rho(x)\,q^{\tt r}(x,x+v)\,r_{\kappa,\vep}(x+v,z) 
\end{aligned}
\eeqq
where $q^{\tt r}(x,y)\le 1$ denotes the density of $Q^{\tt r}(x,dy)$ with respect to $\tilde Q(x,dy)$, see \eqref{boundedness-Q}. 
Here, for each $z_0\in E$ fixed,  the mapping $\,z \to r_{\kappa,\vep}(x+v,z)\,$ is  continuous at $z_0$ whenever $x+v\neq z_0$, by \eqref{proof:continuity-1}, and  $\,z \to \tilde r_\vep(x+v,z)\,$ is  continuous at $z_0$ whenever $x+v\neq z_0$, by lemma \ref{lemma:operators}. 
Note that the set $\{(x,v)\in E^2: x+v=z_0\}$ is a null set under the measure $dx \otimes\widehat Q^{\tt r}(dv)$. 
The image of the measure $dx \otimes\widehat Q^{\tt r}(dv)$ under the mapping $(x,v)\to x+v$ coincides again with Lebesgue measure on $E=\mathbb{R}^d$.  Let us write $M_1$ for the image of the measure 
$$ 
q^{\tt r}(x,x+v) \left(\, \ov\gamma_n(x)\,\kappa(x) \rho(x)\, dx \otimes\widehat Q^{\tt r}(dv) \,\right)  \quad\mbox{on}\quad E\times E
$$
under the mapping $(x,v)\to x+v$. 
By assumption \ref{1.2.1neo}, induction assumption which implies $\ov\gamma_n\in L^1(E)\cap L^\infty(E)$ and since $q^{\tt r}(x,y)\le 1$, the measure $M_1$ on $E$ is finite and admits a Lebesgue density which is bounded: 
$M_1(du)=:m_1(u)du$ for some $m_1\in L^1(\mathbb{R}^d)\cap L^\infty(\mathbb{R}^d)$. Now we rewrite the function on the right hand side in \eqref{xxx1_neo} in the form 
\beqq\label{xxx2_neo}
z \;\lra\; \int_E du\; m_1(u)\; r_{\kappa,\vep}(u,z) 
\eeqq 
where for every $z_0\in E$ fixed, the mapping $z \to  r_{\kappa,\vep}(u,z)$ is continuous at $z_0$  for $M_1$-almost all $u\in E$. 
The rest of the argument is analogous to step 1) above with $Q^{\tt i}(dx)$ replaced by $m_1(u)du$. We only give the details for the case $d\ge2$:
From \eqref{eq:bounds_2} we have bounds of convolution type 
\beqq\label{xxx3_neo}
z \;\lra\; \int_E du\; m_1(u)\; \tilde r_\vep(z-u) 
\eeqq
where $m_1\in L^1(\mathbb{R}^d)\cap L^\infty(\mathbb{R}^d)=\bigcap_{1\le p\le\infty}L^p(\mathbb{R}^d)$ and $\tilde r_\vep(\cdot)\in L^{p^*}(\mathbb{R}^d)$ for $p^*\in(1,\frac{d}{d-2})$ by lemma \ref{lemma:operators}.
Choosing such a $p^*$ and setting $p:=\frac{p^*}{p^*-1}\in(\frac{d}{2},\infty)$ as the dual exponent, we see that the upper bound \eqref{xxx3_neo} is in $\calc_0(\mathbb{R}^d)$, from \cite{HS-65} p.\ 398. 
Again Pratt's lemma applies and shows that the function \eqref{xxx2_neo} is in $\calc_0(\mathbb{R}^d)$. 
We have proved that $\ov\gamma_{n+1}$ is in $\calc_0(\mathbb{R}^d)$ whenever $\ov\gamma_{n}$ is bounded, for all $n\ge 1$. 
This concludes the proof of the lemma. \halmos\\

Now we can finish the \\

{\bf Proof of theorem \ref{1.3.10}:}
In view of lemma \ref{lemma:C_0}, it remains only to show that the series \eqref{eq:series} converges uniformly. By \eqref{eq:series} and the decomposition \eqref{decomp_density}, we have
\[
\begin{aligned}
\ov\gamma_n(z)&=\int_E[Q^{\tt i}(R_\kappa\kappa\varrho\, Q^{\tt r})^n](dx)\,r_\kappa(x,z)\\
&=\int_E[Q^{\tt i}(R_\kappa\kappa\varrho\, Q^{\tt r})^n](dx)\,r_{\kappa,\vep}(x,z)+\int_E[Q^{\tt i}(R_\kappa\kappa\varrho\, Q^{\tt r})^nR_\kappa](dx)\,p_\vep^\kappa(x,z) \;. 
\end{aligned}
\]
Moreover, it is easy to show by induction that
\begin{equation}\label{eq:expansion}
(\kappa\rho\,Q^{\tt r}R_\kappa)^n(x,dy)=(\kappa\rho\,Q^{\tt r}R_{\kappa,\vep})^n+\sum_{k=0}^{n-1}(\kappa\rho\,Q^{\tt r}R_\kappa)^{n-k}P_\vep^\kappa(\kappa\rho\,Q^{\tt r}R_{\kappa,\vep})^k,\qquad n\in\mathbb{N}_0.
\end{equation}
Now we fix $n_0>\frac{d}{2}$. Then we have for all $n>n_0$
\[
\begin{aligned}
\ov\gamma_n(z)
&=\int_E[Q^{\tt i}(R_\kappa\kappa\varrho\,Q^{\tt r})^{n-n_0}(R_\kappa\kappa\varrho\,Q^{\tt r})^{n_0}](dx)\,r_{\kappa,\vep}(x,z)+\int_E[Q^{\tt i}(R_\kappa\kappa\varrho\, Q^{\tt r})^{n}R_\kappa](dx)\,p^\kappa_\vep(x,z)\\
&=\int_Edx\,\ov\gamma_{n-n_0}(x)\int_E \left[ (\kappa\varrho\,Q^{\tt r}R_{\kappa})^{n_0-1}\kappa\rho\, Q^{\tt r} \right] (x,dy)\,r_{\kappa,\vep}(y,z)+\int_Edx\,\ov\gamma_{n}(x)\,p^\kappa_\vep(x,z)
\end{aligned}
\]
where we rewrite the term $(\kappa\varrho\,Q^{\tt r}R_{\kappa})^{n_0-1}$ using \eqref{eq:expansion} with $n_0-1$ in place of $n$:
\[
\begin{aligned}
&=\int_Edx\,\ov\gamma_{n-n_0}(x)\int_E \left[ (\kappa\varrho\,Q^{\tt r}R_{\kappa,\vep})^{n_0-1}\kappa\rho\,Q^{\tt r} \right] (x,dy)\,r_{\kappa,\vep}(y,z)\\
&\qquad+\int_Edx\,\ov\gamma_{n-n_0}(x)\sum_{k=0}^{n_0-2}\int_E[(\kappa\rho\,Q^{\tt r}R_\kappa)^{n_0-1-k}P_\vep^\kappa(\kappa\rho\,Q^{\tt r}R_{\kappa,\vep})^{k}\kappa\rho\,Q^{\tt r}](x,dy)\, r_{\kappa,\vep}(y,z)\\
&\qquad+\int_Edx\,\ov\gamma_{n}(x)\,p^\kappa_\vep(x,z) \;. 
\end{aligned}
\]
Rearranging terms, this last equation takes the form 
\[
\begin{aligned}
\ov\gamma_n(z) 
&=\int_Edx\,\ov\gamma_{n-n_0}(x)\int_E \left[ (\kappa\varrho\,Q^{\tt r}R_{\kappa,\vep})^{n_0-1}\kappa\rho\,Q^{\tt r} \right] (x,dy)\,r_{\kappa,\vep}(y,z)\\
&\qquad+  \sum_{k=0}^{n_0-2}\int_E dx\, \gamma_{n-1-k}(x)\, P_\vep^\kappa(\kappa\rho\,Q^{\tt r}R_{\kappa,\vep})^{k}\kappa\rho\,Q^{\tt r}](x,dy)\, r_{\kappa,\vep}(y,z)\\
&\qquad+\int_Edx\,\ov\gamma_{n}(x)\,p^\kappa_\vep(x,z) \;. 
\end{aligned}
\]
Using the bounds \eqref{eq:bounds_1}--\eqref{boundedness-Q} and the fact that the operators $\tilde P_\vep$, $\tilde R_\vep$ and $\tilde Q$ induce convolutions and thus all commute, the last display is bounded by
\[
\begin{aligned}
&\le\quad \int_Edx\,\ov\gamma_{n-n_0}(x)\|\kappa\varrho\|_\infty^{n_0}\int_E[(\tilde Q\tilde R_{\vep})^{n_0-1}\tilde Q](x,dy)\,\tilde r_{\vep}(y-z)\\
&\qquad+\sum_{k=0}^{n_0-2}\|\kappa\rho\|_\infty^{k+1}\int_Edx\,\ov\gamma_{n-k-1}(x)\int_E[\tilde P_\vep(\tilde Q\tilde R_{\vep})^{k}\tilde Q](x,dy)\,\tilde r_{\vep}(y-z)\\
&\qquad+\int_Edx\,\ov\gamma_{n}(x)\,\tilde p_\vep(x-z)\\
&=\quad \|\kappa\rho\|_\infty^{n_0}\int_Edx\,\ov\gamma_{n-n_0}(x)\int_E\tilde Q^{n_0}(x,dy)\,\tilde r_\vep^{*n_0}(y-z)\\
&\qquad+\sum_{k=0}^{n_0-2}\|\kappa\rho\|_\infty^{k+1}\int_Edx\,\ov\gamma_{n-k-1}(x)\int_E \tilde Q^{k+1}(x,dy)\, [\tilde r_{\vep}^{*(k+1)}*\tilde p_\vep](y-z) \\
&\qquad+[\ov\gamma_n*\tilde p_\vep](z) \;. 
\end{aligned}
\]
Now  $n_0>\frac{d}{2}$ implies that the function $\tilde r_\vep^{*n_0}$ is bounded,  by lemma \ref{lemma:operators}. Thus we obtain a bound
\[
\begin{aligned}
&\le\quad  \widehat Q^{\tt r}(E)^{n_0}\|\kappa\rho\|_\infty^{n_0}\,\|\tilde r_\vep^{*n_0}\|_\infty\,\|\ov\gamma_{n-n_0}\|_1\\
&\qquad+\sum_{k=0}^{n_0-2} \widehat Q^{\tt r}(E)^{k+1}\|\kappa\rho\|_\infty^{k+1}\|\tilde r_{\vep}\|^{k+1}_1\,\|\tilde p_\vep\|_\infty\,\|\ov\gamma_{n-k-1}\|_1 \\
&\qquad+\; \|\ov\gamma_n\|_1\,\|\tilde p_\vep\|_\infty 
\end{aligned}
\]
for $\ov\gamma_n(\cdot)$, again using notations \eqref{boundedness-Q}. Thus we have shown that for $n>n_0$
\[
\|\ov\gamma_n\|_\infty\;\le\; C\,\sum_{k=0}^{n_0}\|\ov\gamma_{n-k}\|_1 
 \;\;=\; C\sum_{m=n-n_0}^n  \|\ov\gamma_{m}\|_1 
\]
where $C>0$ is some constant that does not depend on $n$.
The total mass $\ov\mu(E)=\sum_{n\in\mathbb{N}_0}\|\ov\gamma_n\|_1$  being finite, the last expression is summable in $n>n_0$. 

We have shown that the series \eqref{eq:series} converges uniformly in $z$, and since each term is in $\mathcal{C}_0(E)$ resp.\ $\mathcal{C}_b(E)$ by lemma \ref{lemma:C_0}, the same holds for the limit. This finishes the proof of theorem \ref{1.3.10}. 
\halmos\\

%
\section{Reconstruction of increments for particle trajectories \\when the BDI process is observed discretely in time}\label{reconstruction}


Discretely observed diffusions have received a lot of attention, from Yoshida \cite{Yo-92}, Genon-Catalot and Jacod \cite{GJ-93}, Bibby and S{\o}rensen \cite{BS-95}, Kessler \cite{Ke-97} via Gobet \cite{Go-02} to Podolskij and Vetter \cite{PV-10} or Protter and Jacod \cite{JP-12}.  
Financial data have been a main motivation, and a main issue is estimation of the unknown volatility or of functionals thereof. 
If we observe at discrete time points $\, t_i:=i\Delta\,,\, i\in\mathbb{N}_0,\,$ not a diffusion path but the trajectory of a BDI process $( \eta_t )_{t\ge 0}$, a new type of  problem arises: 
we will be left with pairs of configurations 
$\,(\eta_{i\Delta},\eta_{(i+1)\Delta})\,$, i.e.\ pairs of random point measures on the single-particle space, 
without any information on the path history of the continuous-time process in-between. 
Segments $\,\eta_{[i\Delta,(i+1)\Delta]}\,$ of the trajectory of a BDI process will in general contain branching or immigration events, and even for those particles which succeeded to stay alive over the time interval $[i\Delta,(i+1)\Delta]$ --and thus did travel on diffusion paths-- any information which particle in $\eta_{i\Delta}$ did travel to which position in the configuration $\eta_{(i+1)\Delta}$ will be lost. 
In this section we propose an identification algorithm which  asymptotically as $\Delta\downarrow 0$ will be able to recover correctly, to some large extent, the particle identities in pairs of successive configurations $(\eta_{i\Delta},\eta_{(i+1)\Delta})$.  The algorithm appears in Brandt \cite{Br-05} and was investigated by Berg \cite{Be-15} in dimension $d=1$. 
In the present paper, we show that the result holds in arbitrary dimension $d\ge 1$; heat kernel bounds according to assumption \ref{1.3.8} play a key role.  
The reconstruction algorithm is presented in definition \ref{2.2.2}; the main results are theorems \ref{2.1.2}, \ref{2.2.4} and \ref{2.2.6}. 
\vspace{0.5cm}

\subsection{$\vep$-wellspread configurations, identifiable pairs of configurations, the \\ reconstruction algorithm and the problem of correct identification}\label{2.1}

Recall that we write $x=(x_1,\ldots,x_\ell)$ for configurations $x\in S$ and $x_i=(x_{i,1},\ldots,x_{i,d})$ for particle positions in $E=\bbr^d$. 
In this section, our assumptions will always include \ref{1.2.1neo} and \ref{1.2.2} for arbitrary choice of an immigration measure, and invariant probability measure $\mu$ on $S$ and invariant occupation measure $\ov\mu$ on $E$ are as in lemma \ref{1.2.5}.

\begin{defn}\label{2.1.1}
We call a configuration $x=(x_1,\ldots,x_\ell)$ with $\ell\ge 2$ {\em $\,\vep$-wellspread} if all two-particle subconfigurations $(x_{i_1},x_{i_2})$, $\,1\le i_1< i_2\le \ell$, are such that 
$$
\min_{1\le j\le d}\; |x_{i_1,j} - x_{i_2,j}| \;\ge\; \vep \quad\mbox{for all components $1\le j\le d$} \;. 
$$
We extend the definition to $\ell\in\{0,1\}$ by adopting the convention that one-particle configurations and the void configuration are $\vep$-wellspread.   
\end{defn}

Write ${\tt D}(\vep)$ for the set of $\vep$-wellspread configurations in $S$, and ${\tt N}(\vep) := S\setminus {\tt D}(\vep)$ for its complement:     
\beqq\label{def_neps}
{\tt N}(\vep) \;=\; \left\{ x=(x_1,\ldots,x_\ell) \in S : \;\mbox{$\ell\ge 2$, there is $i_1\neq i_2$ and $j$ such that $|x_{i_1,j}-x_{i_2,j}|<\vep$}\; \right\} \;. 
\eeqq
Then ${\tt N}(\vep)$ is the set of all configurations in $S$ for which at least one pair of particles presents $\vep$-close components. 
The following generalizes \cite{Be-15}, theorem 2.11. \\

\begin{thm}\label{2.1.2}
Assuming  \ref{1.2.1neo}, \,\ref{1.3.2neo} with $q:=3$, \,\ref{1.3.4} and the heat kernel bounds \ref{1.3.8}, 
we have the following asymptotics for ${\tt N}(\vep)$ in  (\ref{def_neps}): 
$$ 
\mu(\, {\tt N}(\vep) \,) \;\;\le\;\; \calo(\vep) \quad\mbox{as $\vep$ tends to $0$} \;. 
$$ 
\end{thm}
\vskip0.5cm

The proof will be given in subsection \ref{proofs_for_2.1}. We turn to discrete observation of the continuous-time BDI process $\eta=(\eta_t)_{t\ge 0}$. Fix $\Delta>0$ and let $\{t_i:=i\Delta : i\in\mathbb{N}_0 \}$ denote a scheme of equidistant observation times. 
Observing discretely in time, pairs of successive observations $(\eta_{i\Delta},\eta_{(i+1)\Delta})$ are merely pairs of finite point measures, possibly of different total mass, without any indication whether or not particles may have died, reproduced or immigrated between times $i\Delta$ and $(i+1)\Delta$, and without any information on particle trajectories in-between. 
Thus discrete observation raises the problem of `particle identification', to be solved prior to all statistical issues. 
In the following, we consider configurations $x=(x_1,\ldots,x_{\ell(x)})$ 
and $y=(y_1,\ldots,y_{\ell(y)})$ in $S$, and let $y\circ\pi$ denote the rearrangement $(y_{\pi(1)},\ldots,y_{\pi(\ell(y))})$ of particles in $y$ by any permutation $\pi$ of $(1,\ldots,\ell(y))$. \\

\begin{defn}\label{2.2.1}
For  $\Delta>0$ and $0 <\la<\frac12$, a pair $(x,y)$ in $S\times S$ is called {\em $(\Delta,\la)$-identifiable\,} if \\
i) $x$ is $4\Delta^\la$-wellspread, \\
ii) $y$ is $2\Delta^\la$-wellspread, \\
iii) $\ell(x)=\ell(y)=:\ell$ for some $\ell\ge 1$, \\
iv) there is some permutation $\pi$ of $(1,\ldots,\ell)$ (in case $\ell=1$, $\pi(1)=1$)
which achieves 
$$
|y_{\pi(i),j} - x_{i,j}| < \Delta^\la \quad\mbox{
for all $i=1,\ldots,\ell$, $j=1,\ldots,d$} \;.  
$$
Note that by ii), the permutation $\pi$ in iv) --if it exists-- is necessarily unique.  We write ${\tt ID}(\Delta,\la)$ for the Borel subset of $(\Delta,\la)$-identifiable pairs $(x,y)$ in $S\times S$.  
\end{defn}
\vskip0.5cm

We call a pair $(\eta_{i\Delta},\eta_{(i+1)\Delta})$ of successive discrete observations  $(\Delta,\la)$-identifiable when $(\eta_{i\Delta},\eta_{(i+1)\Delta})$ takes its value in the set $ {\tt ID}(\Delta,\la)$. Based on \ref{2.2.1} and \ref{2.1.1} we propose a reconstruction algorithm. \\

\begin{defn}\label{2.2.2} (Reconstruction algorithm): 
For $(\Delta,\la)$-identifiable pairs  $(\eta_{i\Delta},\eta_{(i+1)\Delta})$, for $\pi$ the permutation which achieves \ref{2.2.1}~iv) for $x:=\eta_{i\Delta}$ and $y:=\eta_{(i+1)\Delta}$, we decide to view 
$$
\mbox{$y_{\pi(k)}$ as the position at time $(i+1)\Delta$ of the particle which was in position $x_k$ at time $i\Delta$} 
$$ 
for $k=1,\ldots,\ell(x)$. 
\end{defn}
\vskip0.5cm

A decision proposed by algorithm \ref{2.2.2} may be correct or incorrect. $(\Delta,\la)$-identifiability of  a pair of successive observations $(\eta_{i\Delta},\eta_{(i+1)\Delta})$ is defined in terms of the $\si$-field
$$
\calh^\Delta_{i+1} \quad\mbox{with notation}\quad \calh^\Delta_r := \si\left( \eta_{j\Delta} : 0\le j\le r\right) \;,\; r\in\bbn_0  
$$
and the algorithm \ref{2.2.2} proposes a decision on the basis of this information. In order to judge whether the proposed decision identifies particles correctly or fails to do so, we need the continuous-time filtration generated by the process $(\eta_t)_{t\ge 0}$, i.e.\ 
$$
\calf_{(i+1)\Delta}   \quad\mbox{with notation}\quad \calf_t = \bigcap_{r>t} \si\left( \eta_s : 0\le s\le r \right) \;,\; t\ge 0  
$$
and have to consider path segments 
\beqq\label{path_segments}
\eta_{[i\Delta,(i+1)\Delta]} : \quad [i\Delta,(i+1)\Delta] \;\ni\; t \; \lra\; \eta_t  \;\in\; S  
\eeqq
as  $\calf_{(i+1)\Delta}$-measurable random variables taking values in $D([i\Delta,(i+1)\Delta],S)$, the path space of c\`adl\`ag functions $[i\Delta,(i+1)\Delta] \to S$. 
The notion introduced now refers to the larger $\si$-field $\calf_{(i+1)\Delta}$.\\

\begin{defn}\label{2.2.3}
For  $\Delta>0$ and $0 <\la<\frac12$, elements $f$ in $D([i\Delta,(i+1)\Delta],S)$ are   {\em $\;(\Delta,\la)$-good\;}~if  \\
i) $f(i\Delta)$ is $4\Delta^\la$-wellspread, \\
ii) there is $\ell\ge 1$ such that $f$ belongs to the space of continuous functions $C([i\Delta,(i+1)\Delta],E^\ell)$,\\
iii) writing for short $x:=f(i\Delta)$ and $y:=f((i+1)\Delta)$, we have  
$$|y_{k,j} - x_{k,j}| < \Delta^\la  \quad\mbox{for all}\quad k=1,\ldots,\ell \;,\;  j=1,\ldots,d  \;. 
$$  
We write ${\tt CI}(\Delta,\la)$ for the Borel set of $(\Delta,\la)$-good elements in the path spaces  $D([i\Delta,(i+1)\Delta],S)$, irrespectively of $i\in\bbn_0$. 
\end{defn}

Paths segments $\,\eta_{[i\Delta,(i+1)\Delta]}\,$  in (\ref{path_segments}) are called {\em $\;(\Delta,\la)$-good\;}~when $\,\eta_{[i\Delta,(i+1)\Delta]}$ takes its value in the set ${\tt CI}(\Delta,\la)\,$. Then definitions \ref{2.2.3} and \ref{2.2.1} imply the following 
assertions \eqref{fondamental_property_bis} and \eqref{fondamental_property}:  
\beqq\label{fondamental_property_bis}
\left\{\, \eta_{[i\Delta,(i+1)\Delta]} \in {\tt CI}(\Delta,\la) \right\} \;\subset\; \left\{\, 
(\eta_{i\Delta},\eta_{(i+1)\Delta}) 
\in {\tt ID}(\Delta,\la) \right\} \;;
\eeqq
\beqq\label{fondamental_property}
\mbox{on $\,\left\{\, \eta_{[i\Delta,(i+1)\Delta]} \in {\tt CI}(\Delta,\la) \right\}$\,, algorithm \ref{2.2.2} identifies the particles correctly} \;. 
\eeqq
The proportion of $(\Delta,\la)$-identifiable pairs observed up to time $T$ is  
\beqq\label{proportion_ID}
\frac{1}{\lfloor \frac{T}{\Delta} \rfloor } 
\sum_{i=0}^{\lfloor \frac{T}{\Delta} \rfloor - 1} 
1_{{\tt ID}(\Delta,\la)}( \eta_{i\Delta} , \eta_{(i+1)\Delta})  \;, 
\eeqq 
and \eqref{fondamental_property_bis} allows to lower-bound this by 
\beqq\label{proportion_CI}
\frac{1}{\lfloor \frac{T}{\Delta} \rfloor } 
\sum_{i=0}^{\lfloor \frac{T}{\Delta} \rfloor - 1} 
1_{{\tt CI}(\Delta,\la)}(\eta_{[i\Delta,(i+1)\Delta]}) \; 
\eeqq
which is $\calf_T$-measurable.  As a consequence of (\ref{fondamental_property}), ratio (\ref{proportion_CI}) provides a lower bound for the proportion of pairs of successive discrete observations to which algorithm \ref{2.2.2} first applies and second proposes a correct reconstruction of particle identities. 
We aim at lower bounds for (\ref{proportion_CI}) in stationary regime when $0<\la<\frac12$ is fixed and $\Delta>0$ is small enough. 
Below, $\,Q_\mu$ is the law on the canonical path space of the BDI process $\eta=(\eta_t)_{t\ge 0}$ running stationary, i.e.\ with initial law  $\mu$. 
\\

\begin{thm}\label{2.2.4}
Fix $0 <\la<\frac12$.  
\quad a)\; Under the assumptions of theorem \ref{2.1.2}, we have as $\Delta\downarrow 0$
\beqq\label{fixed_step_conv_01}
Q_\mu\left(\;  \ell(\eta_0)\ge 1 \,,\, \eta_{[0,\Delta]} \notin {\tt CI}(\Delta,\la)  \;\right) \;\;\le\;\; \calo(\Delta^\la)  \;, 
\eeqq
\beqq\label{fixed_step_conv_02}
Q_\mu\left(\;   \ell(\eta_0)\ge 1 \,,\, (\eta_0, \eta_\Delta) \in {\tt ID}(\Delta,\la) \; , \; \eta_{[0,\Delta]} \notin {\tt CI}(\Delta,\la) \;\right) \;\;\le\;\; \calo(\Delta)  \;. 
\eeqq
b)\; There is some $\Delta_0>0$ and some constant $D<\infty$ such that for all $\,0<\Delta<\Delta_0$ the following inequalities hold (note that $Q_\mu\left( \ell(\eta_0)\ge 1 \right) $ is strictly smaller than $1$): 
\beqq\label{fixed_step_conv_3}
Q_\mu\left( \ell(\eta_0)\ge 1 \right) - D\, \Delta^\la 
\;<\; Q_\mu\left(\, \eta_{[0,\Delta]} \in {\tt CI}(\Delta,\la) \,\right)    \;<\;
Q_\mu\left( \ell(\eta_0)\ge 1 \right)   
\eeqq 
\beqq\label{fixed_step_conv_4}
Q_\mu\left(\, \eta_{[0,\Delta]} \in {\tt CI}(\Delta,\la) \,\right)  \;<\;
Q_\mu\left(\, (\eta_0,\eta_\Delta)\in {\tt ID}(\Delta,\la) \,\right) \;<\; Q_\mu\left(\, \eta_{[0,\Delta]} \in {\tt CI}(\Delta,\la) \,\right) + D \Delta \;. 
\eeqq
\end{thm}
\vskip0.5cm

Theorem \ref{2.2.4} will be proved in section \ref{proofs_for_2.1}. The main consequence of the theorem is the following; in the language of discretely observed semimartingales, it deals with high frequency observation schemes.  \\

\begin{thm}\label{2.2.6}
Write $T$ for deterministic time horizons.
Fix $\,0<\la<\frac12$ and let $\Delta$ decrease to $0$. Then we have the following convergences in $Q_\mu$-probability:\\
a) When $T \uparrow\infty$, 
$$
\frac{1}{\lfloor \frac{T}{\Delta} \rfloor } 
\sum_{i=0}^{\lfloor \frac{T}{\Delta} \rfloor - 1} 
1_{{\tt CI}(\Delta,\la)}(\eta_{[i\Delta,(i+1)\Delta]}) \;\;\lra\;\; Q_\mu\left( \ell(\eta_0)\ge 1 \right) \;=\; 1-\mu(E^0)  \;.  
$$
Here the limit is deterministic and strictly between $0$ and $1$. \\
b) When $T<\infty$ remains fixed,  
$$ 
\frac{1}{\lfloor \frac{T}{\Delta} \rfloor } 
\sum_{i=0}^{\lfloor \frac{T}{\Delta} \rfloor - 1} 
1_{{\tt CI}(\Delta,\la)}(\eta_{[i\Delta,(i+1)\Delta]})
\;\;\lra\;\; \frac1T\int_0^T 1_{\{\, \ell(\eta_s)\ge 1 \,\}}\, ds  \;.
$$
Here the limit is a random variable taking values in $[0,1]$. \\
c) In both cases a) and b) above, we have 
$$
\frac{1}{\lfloor \frac{T}{\Delta} \rfloor } 
\sum_{i=0}^{\lfloor \frac{T}{\Delta} \rfloor - 1} 
1_{\, \{  \ell(\eta_{i\Delta})\ge 1 \,,\, \eta_{[i\Delta,(i+1)\Delta]}  \,\notin\; {\tt CI}(\Delta,\la) \,\} }  
\;\;=\;\; 
o_{(Q_\mu)}(1) \quad\mbox{as $\;\Delta\downarrow 0$} \;. 
$$
\end{thm}
\vskip0.8cm

On the basis of definitions \ref{2.2.3} and 
\ref{2.2.1} and of  \eqref{fondamental_property_bis}, \eqref{fondamental_property}, \eqref{proportion_ID} and \eqref{proportion_CI}, we can resume theorem \ref{2.2.6} as follows. For high frequency observation schemes, in the sense of asymptotics  a) or~b),  the reconstruction algorithm \ref{2.2.2} first applies to pairs of successive observations $\,(\eta_{i\Delta},\eta_{(i+1)\Delta})\,$ and second  
proposes a correct answer to the problem of particle identification   
$$\mbox{
`for eventually all pairs $\,(\eta_{i\Delta},\eta_{(i+1)\Delta})\,$  which have $\,\ell(\eta_{i\Delta}) > 0\,$'  
}$$
asymptotically as $\Delta\downarrow 0$. It is clear that pairs $\,(\eta_{i\Delta},\eta_{(i+1)\Delta})\,$ with $\eta_{i\Delta}=\delta$ are of no use in view of reconstruction of particle identities. 
The proof of theorem \ref{2.2.6} will be given in section \ref{proofs_for_2.1}. \\

\subsection{\bf Proofs for subsection \ref{2.1}  }\label{proofs_for_2.1}


We prove the results of the preceding section. \\

{\bf Proof of theorem \ref{2.1.2}: } 1) Recall that assumption \ref{1.3.4} implies \ref{1.2.2} for arbitrary choice of an immigration measure.
Under \ref{1.2.1neo},  \ref{1.3.2neo} with $q:=3$ and \ref{1.3.4} we know from theorem \ref{1.3.7} 
\beqq\label{third_moments_mubar}
\mu(\ell^3) \;=\; \sum_{\ell\in\bbn} \ell^3\; \mu( E^\ell ) \;\;<\;\; \infty \;. 
\eeqq
 
2) 
With respect to $\vep>0$ arbitrary but fixed --which for a while we suppress from notations-- write for two-particle configurations  $(x_1, x_2)\in E^2$ with $x_i=(x_{i,1},\ldots,x_{i,d})$ 
$$ 
h(x_1,x_2) := \sum_{j=1}^d 1_{ \{ |x_{1,j}-x_{2,j}|<\vep \} } 
$$
and define a function $g:S\to[0,\infty)$ by  
$$
g(x_1,\ldots,x_\ell) \;:=\; \sum_{ 1\le i_1<i_2\le \ell}  h(x_{i_1},x_{i_2}) \quad\mbox{for}\;\; \ell\ge 2 \quad,\quad g\equiv 0 \;\;\mbox{on}\; E^0\cup E^1 \;. 
$$
Then we can write for short 
\beqq\label{repr-neps}
1_{ {\tt N}(\vep) }(x) \;\;\le\;\; g(x)\; \quad\mbox{for all}\quad x\in S  
\eeqq
and have for the invariant probability $\mu$ on $S$ (use \eqref{inv-1} in the proof of lemma \ref{1.2.5} plus norming) 
\beqq\label{repr-neps_neu2}
\mu(\,{\tt N}(\vep)\,) \;\;\le\;\; \int_S g\, d\mu \;\;=\;\; \frac{1}{E_\delta(R_1)}\, 
E_{\delta} \left( \int_0^{R_1} g(\eta_s)\, ds \right)   
\eeqq
where --writing  $(T_j)_j$ for the sequence of jump times in the BDI process, and using the Markov property-- 
\eqref{inv-1} implies that  
\beqq\label{repr-neps_neu3}
E_{\delta} \left( \int_0^{R_1} g(\eta_s)\, ds \right)   \;\;=\;\;  
\sum_{n=0}^\infty E_{\delta} \left( \,1_{\{T_n<R_1\}}\, \int_{T_n}^{T_{n+1}} g(\eta_s)\, ds\right)   
\eeqq
where we can write 
\beqq\label{repr-neps_neu4}
E_{\delta} \left( \,1_{\{T_n<R_1\}}\, \int_{T_n}^{T_{n+1}} g(\eta_s)\, ds\right)  \;\;=\;\;
E_{\delta} \left( \,1_{\{T_n<R_1\}}\, E_{\eta_{T_n}}\left( \int_0^{T_1} g(\eta_s)\, ds \right) \right) \;. 
\eeqq
\vskip0.5cm

3) We insert an auxiliary step and prove the bound  
\beqq\label{updated_proof_bound1}
\int_{E^2} [ P^\kappa_t(x_1,dy_1){\otimes}P^\kappa_t(x_2,dy_2) ]\; h(y_1,y_2) \;\;\le\;\; 
d\; C^2 (2\pi C)^d\; 2\vep\; \frac{1}{\sqrt{2\pi C t}}
\eeqq
valid for  $0< t\le t_0$, where
we recall that $(P^\kappa_t)_{t\ge 0}$ is the semigroup of the one-particle motion killed at rate $\kappa$ and $t_0$ is from assumption \ref{1.3.8} (note that the right hand side of \eqref{updated_proof_bound1} is free from $x_1$ and $x_2$).  
Indeed, for such $t$ 
the left hand side of \eqref{updated_proof_bound1} is by \ref{1.3.8} smaller than  
\beqq\label{updated_proof_bound2}
C^2 (2\pi C)^d \int_{E^2} p_{Ct}(y_1-x_1)\, p_{Ct}(y_2-x_2)\; h(y_1,y_2)\;  dy_1\, dy_2  
\eeqq
where $p_s$ denotes the density (on $\mathbb{R}^d$) at time $s$ for $d$-dimensional standard Brownian motion. Since 
$$ 
h(y_1,y_2) = \sum_{j=1}^d 1_{ \{ |y_{1,j}-y_{2,j}|<\vep \} }  = \sum_{j=1}^d  1_{B_\vep(y_{1,j})}(y_{2,j})
$$
depends on $\vep$ (here $B_\vep$ denotes a ball in $\mathbb{R}^1$), it is sufficient to calculate in \eqref{updated_proof_bound2} the maximum of a one-dimensional normal density to prove \eqref{updated_proof_bound1}.

4) Still keeping $\vep>0$ fixed and suppressed from notations, we evaluate  typical terms in the expressions \eqref{repr-neps_neu3} and \eqref{repr-neps_neu4} in order to establish a bound 
\beqq\label{updated_proof_bound3_bis}
E_{(x_1,\ldots,x_\ell)}\left( \int_0^{T_1} g(\eta_s)\, ds \right) \;\;\le\;\;   {\tt cst}\; \vep\; d\; \frac{\ell(\ell-1)}{2}\; (t_0^{-1/2} +t_0^{1/2})   \;; 
\eeqq
here and below, `$\tt cst$' collects constants which are not of interest (and which may change from line to line). 
To prove \eqref{updated_proof_bound3_bis}, we start from 
\beqq\label{updated_proof_bound3}
E_{(x_1,\ldots,x_\ell)}\left( \int_0^{T_1} g(\eta_s)\, ds \right) \;=\; \sum_{ 1\le i_1<i_2\le \ell} E_{(x_1,\ldots,x_\ell)}\left( \int_0^{T_1}  h(\eta^{i_1}_s,\eta^{i_2}_s) \, ds \right)
\eeqq
for $(x_1,\ldots,x_\ell)\in S$ and $\ell\ge 1$, with $g$ and $h$ as above. We have by definition of the BDI process 
\beao
&&E_{(x_1,\ldots,x_\ell)}\left( \int_0^{T_1}  h(\eta^{i_1}_s,\eta^{i_2}_s) \, ds \right) \;=\; 
E_{(x_1,\ldots,x_\ell)}\left( \int_0^\infty dt\; e^{- \int_0^t (c+\ov\kappa(\eta_v)) dv }\; h(\eta^{i_1}_t,\eta^{i_2}_t) \,\right) \\
&&\quad=\quad  \int_0^\infty dt\; e^{-ct}  \int_{E^\ell} \prod_{j=1}^\ell P^\kappa_t(x_j,dy_j)\;\; h(y_{i_1},y_{i_2}) \\
&&\quad\le\quad  \int_0^\infty dt\; e^{-ct}  \int_{E^2} [ P^\kappa_t(x_{i_1},dy_{i_1}){\otimes}P^\kappa_t(x_{i_2},dy_{i_2}) ]\; h(y_{i_1},y_{i_2}) .
\eeao
Let us consider the initial part   $\int_0^{t_0} dt \ldots$  of the last integral first: using \eqref{updated_proof_bound1},   
$$
\int_0^{t_0} dt\; e^{-ct}  \int_{E^2} [ P^\kappa_t(x_{i_1},dy_{i_1}){\otimes}P^\kappa_t(x_{i_2},dy_{i_2}) ]\; h(y_{i_1},y_{i_2}) 
$$
is bounded above by
\beqq\label{updated_proof_bound4}
\int_0^{t_0} dt\; e^{-ct}\;  d\; C^2 (2\pi C)^d\; 2\vep\; \frac{1}{\sqrt{2\pi C t}} 
\quad\le\;\;  {\tt cst}\; \vep\; d\; t_0^{1/2} \;. 
\eeqq
Turning to the remaining part  $\int_{t_0}^\infty dt \ldots$  of the  integral above, we shall prove a bound 
\beqq\label{updated_proof_bound5}
\int_{t_0}^\infty dt\; e^{-ct}  \int_{E^2} [ P^\kappa_t(x_{i_1},dy_{i_1}){\otimes}P^\kappa_t(x_{i_2},dy_{i_2}) ]\; h(y_{i_1},y_{i_2}) \;\;\le\;\;  {\tt cst}\; \vep\; d\; t_0^{-1/2} \;. 
\eeqq
To see this, write the left hand side of \eqref{updated_proof_bound5} as 
$$
e^{- c t_0} \int_0^\infty dv\; e^{-cv}\; F(v,t_0,x_{i_1},x_{i_2}) 
$$  
with short notation 
$$
F(v,t_0,x_{i_1},x_{i_2}) \;:=\;   \int_{E^2} \left[ [ P^\kappa_v P^\kappa_{t_0} ](x_{i_1},dy_{i_1}){\otimes}[ P^\kappa_v P^\kappa_{t_0} ](x_{i_2},dy_{i_2}) \right]\; h(y_{i_1},y_{i_2}) \;.   
$$
Rearranging terms and applying again  \eqref{updated_proof_bound1}, this is smaller than 
\beao
&&\int_{E^2} [ P^\kappa_v(x_{i_1},dz_{i_1}){\otimes}P^\kappa_v(x_{i_2},dz_{i_2}) ]   \int_{E^2} [ P^\kappa_{t_0}(z_{i_1},dy_{i_1}){\otimes}P^\kappa_{t_0}(z_{i_2},dy_{i_2}) ]\; h(y_{i_1},y_{i_2}) \\
&&\le\quad  \int_{E^2} [ P^\kappa_v(x_{i_1},dz_{i_1}){\otimes}P^\kappa_v(x_{i_2},dz_{i_2}) ]\; d\; C^2 (2\pi C)^d\; 2\vep\; \frac{1}{\sqrt{2\pi C t_0}}  \quad\le\quad d\; C^2 (2\pi C)^d\; 2\vep\; \frac{1}{\sqrt{2\pi C t_0}}
\eeao
which gives \eqref{updated_proof_bound5}. Combining \eqref{updated_proof_bound5}, \eqref{updated_proof_bound4} and \eqref{updated_proof_bound3} we have proved \eqref{updated_proof_bound3_bis}. 

5) Now we insert the bound \eqref{updated_proof_bound3_bis} into the three equations \eqref{repr-neps_neu2}, \eqref{repr-neps_neu3} and \eqref{repr-neps_neu4} which gives
\beao
\mu(\,{\tt N}(\vep)\,) &\le& {\tt cst}\; \sum_{n=0}^\infty E_\delta\left(\, 1_{\{T_n<R_1\}}\; E_{\eta_{T_n}}\left( \int_0^{T_1} g(\eta_s)\, ds \right) \right) \\
&\le&  {\tt cst}\;\vep\;d\;(t_0^{-1/2}+t_0^{1/2})\;   \sum_{n=0}^\infty E_\delta\left(\, 1_{\{T_n<R_1\}}\; \left(\ell(\eta_{T_n})\right)^2  \right) \;. 
\eeao
Now from the definition of the BDI process, for all $x\in S$ with $\ell(x)\ge 1$, 
$$
1 = E_x\left( \int_0^\infty dt\; (c+\ov\kappa(\eta_t)) e^{-\int_0^t (c+\ov\kappa(\eta_v)) dv} \right) \le {\tt cst}\; \ell(x)\; E_x( T_1 ) \;. 
$$
Then, absorbing also the dimension $d$ and the term $(t_0^{-1/2}+t_0^{1/2})$ into the constants, we obtain 
\beao
\mu(\,{\tt N}(\vep)\,) &\le& {\tt cst}\; \vep\; \sum_{n=0}^\infty E_\delta\left(\, 1_{\{T_n<R_1\}}\; \left(\ell(\eta_{T_n})\right)^3\, E_{\eta_{T_n}}( T_1 ) \right) \\ 
&\le& {\tt cst}\;\vep\;  \sum_{n=0}^\infty E_\delta\left(\, 1_{\{T_n<R_1\}}\;E_{\eta_{T_n}}\left(\int_0^{T_1} (\ell(\eta_s))^3\, ds\right) \right) \\ 
&\le& {\tt cst}\;\vep\;  E_\delta\left(\int_0^{R_1}(\ell(\eta_s))^3\, ds \right) 
\eeao
which gives 
$$
\mu(\,{\tt N}(\vep)\,) \;\;\le\;\;  {\tt cst}\; \vep\; \mu(\, \ell^3\,)
$$
and finishes the proof of theorem \ref{2.1.2}. \halmos\\

{\bf Proof of theorem \ref{2.2.4}: } 
1) In a first step, for $0<\la<\frac12$ fixed and $\Delta>0$ small enough, define 
$$
h_{\Delta,\la}(x,y) := \sum_{j=1}^d 1_{\{ |x_j-y_j|>\Delta^\la \}} \;\;,\;\; x,y\in E=\mathbb{R}^d
$$
for one-particle configurations. We shall show that asymptotically as  $\Delta\downarrow 0$, 
\beqq\label{collect_current_result_step_1}
\int_E P_\Delta^\kappa(x,dy)\, h_{\Delta,\la}(x,y)  \;\;=\;\; o(\Delta)   
\eeqq
where $o(\Delta)$ denotes bounds which do not depend on $x\in E$. 
Indeed, denoting by $\Phi$ resp.\ $\varphi$ the distribution function resp.\ density of the standard normal law on $\mathbb{R}$, as in the proof of theorem \ref{2.1.2} the heat kernel bounds \ref{1.3.8} yield  
\beao
\int_E P_\Delta^\kappa(x,dy)\, h_{\Delta,\la}(x,y) &\le& C \sqrt{2\pi C}^d \int_E dy\, p_{C\Delta}(y-x)\, h_{\Delta,\la}(x,y) \\
&=& d\; \,C \sqrt{2\pi C}^d\; 2\left(\, 1-\Phi( \Delta^\la / \sqrt{C\Delta}  ) \,\right).
\eeao
Using $0<\la<\frac12$ and the elementary inequality $\,1-\Phi(v) < \frac1v 
\varphi(v)\,$ which is true for all $v>0$, this in turn  is bounded by 
$$
\le\;\; {\tt cst}\; \Delta^{\frac12-\la}\; e^{-\frac{1}{2C} \Delta^{2\la-1} } 
\;=\; {\tt cst}\; \Delta \;\; (\frac1\Delta)^{\frac12+\la}\; e^{-\frac{1}{2C} (\frac1\Delta)^{1-2\la} } \;=\; o(\Delta)
$$ 
as $\Delta\downarrow 0$. Here and below, `$\tt cst$' denotes constants which may change from line to line.

2) Consider a one-particle motion $\xi$ killed at rate $\kappa$, starting at time $0$ from $x\in E$. Over a time interval of length $\Delta$ we will have three possibilities. 
Either killing occurs before time $\Delta$, i.e.\ $\xi_\Delta=\wh\delta$ where $\wh\delta$ represents some cemetery point for processes with life time, 
or $\xi_\Delta$ takes values in a cube $U_{\Delta,\la}(x):=  \mathop{\sf X}\limits_{j=1}^d (x_j{-}\Delta^\la,x_j{+}\Delta^\la)$ centred at $x=(x_1,\ldots,x_d)$, or 
$\xi_\Delta$ takes values in $E\setminus U_{\Delta,\la}(x)$. If we write $\wh P^\kappa_t(\cdot,\cdot)$ for the semigroup of the one-particle motion including possible jumps to the cemetery point  on the extended state space $\wh E := \{\wh\delta\}\cup E$, and as before $P^\kappa_\Delta(\cdot,\cdot)$ for the killed semigroup on $E$ in assumption \ref{1.3.8}, this means that 
asymptotically as $\Delta\downarrow 0$, 
$$
\wh P^\kappa_\Delta\left(x,\left( \{\wh\delta\}\cup E\setminus U_{\Delta,\la}(x) \right) \right)  \;\le\;  \|\kappa\|_{\infty}\, \Delta  +  P^\kappa_\Delta(x, E\setminus U_{\Delta,\la}(x) )  
\;\le\;  \|\kappa\|_{\infty}\, \Delta  \;+\;   o(\Delta) 
$$
where we have used  \eqref{collect_current_result_step_1}. The bounds do not depend on $x\in E$, for all $\Delta$ sufficiently small. 
We turn to $\ell$-particle motions, $\ell\ge 1$,  killed and jumping to some cemetery point $\wh\delta$ at rate $c+\ov\kappa(\cdot)$, and write $\wh P^{c+\ov\kappa}_t(\cdot,\cdot)$ for its semigroup on the extended state space $\wh E^\ell := \{\wh\delta\}\cup E^\ell$. 
Independence of the killed motions of individual particles shows that  
$$
\wh P^{c+\ov\kappa}_\Delta\left( x , \left( \{\wh\delta\}\cup E^\ell\setminus \mathop{\sf X}\limits_{i=1}^\ell U_{\Delta,\la}(x_i) \right)  \right)  
\;\le\;  ( c + \|\kappa\|_{\infty} \ell ) \Delta \;+\; \sum_{i=1}^\ell P^\kappa_\Delta\left( x_i,  E\setminus U_{\Delta,\la}(x_i) \right)
$$
for starting positions $x=(x_1,\ldots,x_\ell)$ in $E^\ell$, where again by \eqref{collect_current_result_step_1}  
\beqq\label{bound_multiparticle_motion_killed}
\wh P^{c+\ov\kappa}_\Delta\left( x , \left( \{\wh\delta\}\cup E^\ell\setminus \mathop{\sf X}\limits_{i=1}^\ell U_{\Delta,\la}(x_i) \right)  \right)  
\;\;\le\;\;  {\tt cst}\; \ell\;  \Delta   \quad\mbox{as $\Delta\downarrow 0$}    
\eeqq
for all $\ell\ge 1$, with some global constant not depending on $\ell\ge 1$ or  $x\in E^\ell$.

3) So far we have exploited assumptions  \ref{1.2.1neo} and  \ref{1.3.8}. The following argument will exploit  \ref{1.3.2neo} (with $q:=3$) together with \ref{1.3.4}  (which implies \ref{1.2.2} for arbitrary choice of an immigration measure), and will conclude the proof.  Consider a path segment $\eta_{[i\Delta,(i+1)\Delta]}$ under $Q_\mu$. 
By stationarity it  is sufficient to consider $i=0$, with random initial position $\eta_0\in S$ selected according to invariant measure $\mu$. Directly from definitions \ref{2.2.3}, \ref{2.2.1} and \eqref{def_neps} we have the inclusion 
$$
\left\{   \ell(\eta_0)\ge 1 \,,\, \eta_{[0,\Delta]} \notin {\tt CI}(\Delta,\la)  \right\} 
\;\;\subset\;\; \left\{ \eta_0 \in {\tt N}(4\Delta^\la) \right\} \;\cup\; \left\{  \ell(\eta_0)\ge 1 \,,\, \eta_0 \in {\tt D}(4\Delta^\la) \,,\,\eta_{[0,\Delta]}\notin {\tt CI}(\Delta,\la)\right\}   \;. 
$$
Whenever an initial configuration $x\in S$ with $\ell=\ell(x)\ge 1$ is $4\Delta^\la$-wellspread, we consider as in step~2) the $\ell$-particle motion starting from $x$, killed and jumping to a cemetery point $\wh\delta$ at rate $c+\ov{\kappa}(\cdot)$; then 
\beqq\label{bound1_segment_notCI}
Q_x\left( \eta_{[0,\Delta]} \notin {\tt CI}(\Delta,\la) \right) \;\;=\;\;  \wh P^{c+\ov\kappa}_t\left( x , \left( \{\wh\delta\}\cup E^\ell\setminus \mathop{\sf X}\limits_{i=1}^\ell U_{\Delta^\la}(x_i) \right)  \right)  \quad\le\quad 
 {\tt cst}\; \ell\; \Delta 
\eeqq
from (\ref{bound_multiparticle_motion_killed}), with $x=(x_1,\ldots,x_\ell)$.
As a consequence, $\,\ov\mu (1)\, = \sum_\ell \ell\, \mu(E^\ell)\,$ being finite in virtue of lemma \ref{1.2.5}, we arrive at  
\beqq\label{bound2_segment_notCI}
Q_\mu\left(\, \ell(\eta_0)\ge 1  \,,\, \eta_0 \in {\tt D}(4\Delta^\la) \,,\,  \eta_{[0,\Delta]} \notin {\tt CI}(\Delta,\la)  \, \right)  \;\;\le\;\;   {\tt cst} \; \Delta \; \sum_{\ell=1}^\infty \ell\, \mu(E^\ell)  \;\;=\;\; \calo(\Delta)   
\eeqq
as $\Delta\downarrow 0$. 
So far, from the above inclusion and the bound (\ref{bound2_segment_notCI}), 
$$
Q_\mu\left(\, \ell(\eta_0)\ge 1 \,,\, \eta_{[0,\Delta]} \notin {\tt CI}(\Delta,\la) \,\right) \;\le\; \mu\left( {\tt N}(4\Delta^\la) \right) + \calo(\Delta) 
$$ 
as $\Delta\downarrow 0$.   
Now we make use of assumption \ref{1.3.2neo} with $q:=3$: applying theorem \ref{2.1.2} to the first term on the right hand side, we have proved assertion \eqref{fixed_step_conv_01} in theorem \ref{2.2.4}. Note that the rate in \eqref{fixed_step_conv_01} comes from the exceptional set ${\tt N}(4\Delta^\la)$ in the above inclusion. \\
Assertion \eqref{fixed_step_conv_02} follows from the bound \eqref{bound2_segment_notCI} since $\,\{ (\eta_{0}, \eta_{\Delta}) \in {\tt ID}(\Delta,\la)\} \subset \{ \eta_{0}\in {\tt D}(4\Delta^\la) \}\,$ by definition in \ref{2.2.1}~i). We have proved part a) of theorem \ref{2.2.4}.

4) We prove part b) of the theorem.  
By definition in \ref{2.2.3}, path segments with $\eta_0=\delta$ never belong to ${\tt CI}(\Delta,\la)$, thus 
\beqq\label{coupling_CI2ell_1}
\left\{ \eta_{[0,\Delta]} \in {\tt CI}(\Delta,\la) \right\} \dot\cup \left\{ \ell(\eta_0)\ge 1 \,,\, \eta_{[0,\Delta]} \notin {\tt CI}(\Delta,\la) \right\} \;=\;  \left\{ \ell(\eta_0)\ge 1 \right\} \;. 
\eeqq
By \eqref{fixed_step_conv_01} in part~a), the $Q_\mu$-probability of the second event on the left hand side  is $\calo(\Delta^\la)$ as $\Delta\downarrow 0$. 
This establishes \eqref{fixed_step_conv_3}. Similarly, \eqref{fondamental_property_bis} allows to write 
$$
\left\{ 
\eta_{[0,\Delta]} \in {\tt CI}(\Delta,\la)
\right\} \dot\cup \left\{ 
(\eta_{0},\eta_{\Delta}) \in {\tt ID}(\Delta,\la)
\,,\, \eta_{[0,\Delta]} \notin {\tt CI}(\Delta,\la) \right\} \;=\;  \left\{ 
(\eta_{0},\eta_{\Delta}) \in {\tt ID}(\Delta,\la)
 \right\}  
$$
where the $Q_\mu$-probability of the second event on the left hand side is $\calo(\Delta)$ as $\Delta\downarrow 0$, by \eqref{fixed_step_conv_02} in part~a).  
Note that $(\eta_{0},\eta_{\Delta}) \in {\tt ID}(\Delta,\la)$ implies $\,\ell(\eta_0)\ge 1\,$ by definition \ref{2.2.1}. This establishes \eqref{fixed_step_conv_4}. The proof of theorem \ref{2.2.4} is finished. \halmos\\

We mention a consequence of theorem \ref{2.2.4}, of minor importance. \\

\begin{prop}\label{2.2.5}
Fix  $0<\la<\frac12$ and $\Delta\in(0,\Delta_0)$, with $\Delta_0>0$ from part b) of theorem \ref{2.2.4}. Then for arbitrary choice of a starting point $x\in S$, we have $Q_x$-almost sure convergence as $N\to\infty$ 
\beqq\label{fixed_step_conv_1}
\frac1N \sum_{i= 0}^{N-1} 1_{{\tt CI}(\Delta,\la)}(\eta_{[i\Delta,(i+1)\Delta]})   \;\lra\; Q_\mu\left(\, \eta_{[0,\Delta]} \in {\tt CI}(\Delta,\la) \,\right)  
\eeqq
\beqq\label{fixed_step_conv_2}
\frac1N \sum_{i=0}^{N-1} 1_{{\tt ID}(\Delta,\la)}(\eta_{i\Delta},\eta_{(i+1)\Delta})  \;\lra\; Q_\mu\left(\, (\eta_0,\eta_\Delta)\in {\tt ID}(\Delta,\la) \,\right)
\eeqq
where the limits \eqref{fixed_step_conv_1} and \eqref{fixed_step_conv_2} are such that 
inequalities \eqref{fixed_step_conv_3} and \eqref{fixed_step_conv_4} hold for $0<\Delta<\Delta_0$. 
\end{prop}
\vskip0.5cm

{\bf Proof: }
By lemma \ref{1.2.5}, the continuous-time process $\eta=(\eta_t)_{t\ge 0}$ is positive Harris. A particular feature of the BDI model is that $\eta$ returning infinitely often to the void configuration $\delta$ will remain there during an independent exponential time with parameter $c>0$. 
As a consequence, for $\Delta>0$ arbitrary but fixed, there will be an infinite number of intervals $[j\Delta,(j+1)\Delta]$ on which $(\eta_t)_{t\ge 0}$ remains visiting the void configuration. Thus the Markov chain of path segments $\left( \eta_{[i\Delta,(i+1)\Delta]} \right)_{i\in\bbn_0}$, taking values in $D([0,\Delta],S)$, will admit an infinite number of visits in state $\equiv \delta$ viewed as a path in $D([0,\Delta],S)$). 
As consequences of this fact, both the path segment chain and the chain of successive pairs $(\eta_{i\Delta},\eta_{(i+1)\Delta} )_{i\in\mathbb{N}_0}$ are positive Harris chains. We thus have strong laws of large numbers: the 
rescaled additive functionals on the left hand side in \eqref{fixed_step_conv_1} and \eqref{fixed_step_conv_2} converge $Q_x$-almost surely, for every choice of a starting point $x\in S$, to the limits on the right. \halmos\\

We explain why proposition \ref{2.2.5} is of minor importance. For $0<\la<\frac12$ fixed and for $0<\Delta<\Delta_0$ small but fixed, drawing a large number of discrete-time observations $\,\eta_{i\Delta}$, $0\le i\le N$  (the asymptotics is in $N\to\infty$), the reconstruction algorithm \ref{2.2.2} applies to a proportion of roughly 
$$
Q_\mu\left(\, (\eta_0,\eta_\Delta)\in {\tt ID}(\Delta,\la) \,\right)   
\quad>\quad 
Q_\mu\left( \ell(\eta_0)\ge 1 \right) - D\, \Delta^\la 
$$
pairs of successive observations $(\eta_{i\Delta},\eta_{(i+1)\Delta})$, using theorem \ref{2.2.4}; it is clear that observed pairs with $\ell(\eta_{i\Delta})=0$ are of no use in view of the algorithm. 
If the decision proposed by algorithm \ref{2.2.2} will be correct for a large amount of the data to which the algorithm applies, there will remain some small proportion of approximately 
$$
0 \quad<\quad Q_\mu\left(\, (\eta_0,\eta_\Delta)\in {\tt ID}(\Delta,\la) \;,\;  \eta_{[0,\Delta]} \notin {\tt CI}(\Delta,\la) \,\right)   \quad<\quad   D\; \Delta  
$$
per cent of the data on which the algorithm \ref{2.2.2} may take a decision which fails to identify the underlying (unobserved) travelling particles correctly. Note that $\Delta$ is small but fixed. 
In the language of discretely observed semimartingales (\cite{JP-12}, \cite{PV-10}), proposition \ref{2.2.5} belongs to the frame\-work of `low frequency' asymptotics. A fully satisfactory answer to the problem of particle identities requires a setting of `high frequency' observation, i.e.\  $\Delta$ tending to $0$. We prove theorem \ref{2.2.6}.\\

{\bf Proof of theorem \ref{2.2.6}: }
 1) For deterministic $T$, stationarity allows to write 
\beao
&&Q_\mu\left( \frac{1}{\lfloor \frac{T}{\Delta} \rfloor } 
\sum_{i=0}^{\lfloor \frac{T}{\Delta} \rfloor - 1}
1_{ \{  \ell(\eta_{i\Delta})\ge 1 \,,\, \eta_{[i\Delta,(i+1)\Delta]}  \,\notin\; {\tt CI}(\Delta,\la) \} }  \;>\; \vep \right) \\
&&\le\quad \frac{1}{\vep}\;  E_\mu \left( \frac{1}{\lfloor \frac{T}{\Delta} \rfloor } 
\sum_{i=0}^{\lfloor \frac{T}{\Delta} \rfloor - 1}
1_{ \{   \ell(\eta_{i\Delta})\ge 1 \,,\, \eta_{[i\Delta,(i+1)\Delta]}  \,\notin\; {\tt CI}(\Delta,\la) \} } \right)  \\[2mm]
&&=\quad \frac{1}{\vep}\;  Q_\mu\left(\,  \ell(\eta_{0})\ge 1 \,,\, \eta_{[0,\Delta]} \notin {\tt CI}(\Delta,\la) \,\right) \quad\le\quad \frac{1}{\vep}\, \calo(\Delta^\la)     
\eeao
as $\Delta \downarrow 0$, using theorem \ref{2.2.4}.  As a consequence, irrespectively of the behaviour of $T$, 
$$
\frac{1}{\lfloor \frac{T}{\Delta} \rfloor } 
\sum_{i=0}^{\lfloor \frac{T}{\Delta} \rfloor - 1}
1_{\, \{  \ell(\eta_{i\Delta})\ge 1 \,,\, \eta_{[i\Delta,(i+1)\Delta]}  \,\notin\; {\tt CI}(\Delta,\la) \,\} }  
\;\;=\;\; 
o_{(Q_\mu)}(1) \quad\mbox{as $\;\Delta\downarrow 0$} \;. 
$$
This is c). Together with a decomposition as in (\ref{coupling_CI2ell_1}) we obtain the  $Q_\mu$-stochastic equivalence  
\beqq\label{coupling_CI2ell_2}
\frac{1}{\lfloor \frac{T}{\Delta} \rfloor } 
\sum_{i=0}^{\lfloor \frac{T}{\Delta} \rfloor - 1}
1_{{\tt CI}(\Delta,\la)}(\eta_{[i\Delta,(i+1)\Delta]}) 
\;\;=\;\; 
\frac{1}{\lfloor \frac{T}{\Delta} \rfloor } 
\sum_{i=0}^{\lfloor \frac{T}{\Delta} \rfloor - 1} 
1_{ \{ \ell(\eta_{i\Delta}) \ge 1   \} }  \;+\; o_{(Q_\mu)}(1) \quad\mbox{as $\;\Delta\downarrow 0$} 
\eeqq
which holds in both cases under consideration: the case where time horizon $T$ is fixed and finite, and the case where $T$ is increasing to $\infty$. 
We underline that high-frequency asymptotics $\Delta\downarrow 0$ is a necessary condition for \eqref{coupling_CI2ell_2}. 
\\
2) When $T<\infty$ is fixed and $\Delta\downarrow 0$, the following convergence 
\beqq\label{pathwise_limit_ell_1}
\frac{1}{\lfloor \frac{T}{\Delta} \rfloor } 
\sum_{i=0}^{\lfloor \frac{T}{\Delta} \rfloor - 1}
1_{ \{ \ell(\eta_{i\Delta}) \ge 1   \} }  
\;\;\lra\;\; 
\frac1T\int_0^T 1_{\{ \ell(\eta_s)\ge 1 \}}\, ds    
\eeqq
holds pathwise since every path of $\,(\eta_t)_{t\ge 0}\,$ is a c\`adl\`ag function. Here the limit is a random variable taking values in $[0,1]$. Combining (\ref{pathwise_limit_ell_1}) and (\ref{coupling_CI2ell_2}) we have proved b). 
\\
3) When $T\uparrow\infty$ and $\Delta\downarrow 0$, the following convergence   
\beqq\label{pathwise_limit_ell_2}
\lim\limits_{ T\uparrow\infty }\; \lim\limits_{ \Delta\downarrow 0  }\; \frac{1}{\lfloor \frac{T}{\Delta} \rfloor } 
\sum_{i=0}^{\lfloor \frac{T}{\Delta} \rfloor - 1}
1_{ \{ \ell(\eta_{i\Delta}) \ge 1   \} }  
\;\;=\;\; \lim\limits_{ T\uparrow\infty }\; \frac1T\int_0^T 1_{\{ \ell(\eta_s)\ge 1 \}}\, ds 
\;=\; Q_\mu\left( \ell(\eta_0)\ge 1 \right) 
\eeqq
holds for $Q_\mu$-almost all paths of the BDI process $(\eta_t)_{t\ge 0}$, as a consequence of positive Harris recurrence by lemma \ref{1.2.5}. Combining (\ref{pathwise_limit_ell_2}) and (\ref{coupling_CI2ell_2}) we have proved assertion a) of theorem \ref{2.2.6}. 
\halmos \\

\section{Regression schemes for estimation of the diffusion coefficient \\from discrete BDI observations }\label{regressionNEU}


This section needs all assumptions of sections \ref{setting} and \ref{1.3}, and in particular relies heavily on theorems \ref{1.3.10} and \ref{2.2.4}. For the one-particle motion in {\bf (A1)}, the additional assumption 
\beqq\label{invertible_diff_coeffNEU}
a(y) := (\si \si^\top)(y) \;\;\mbox{invertible for all}\;\; y\in E=\bbr^d
\eeqq  
will be in force. Throughout, the continuous-time BDI process $(\eta_t)_{t\ge 0}$ is stationary,  $\,Q_\mu$ is the stationary law on the canonical path space as in section \ref{reconstruction}, and we deal with discrete-time observation at step size $\Delta$ when $\Delta$ is small. \\

\subsection{The regression scheme: construction and properties}\label{3.1NEU}


By theorem \ref{1.3.10}, the invariant occupation measure $\ov\mu$ on the single particle space  $E=\bbr^d$ admits a Lebesgue density  $\ov \gamma \in \calc_0(E)$. 
Fix any cube $A$ in $E$ such that 
\beqq\label{M1NEU}
M_1 \;:=\; \inf \{ \ov\gamma(x) : x\in A \}  \;\;>\;\; 0 
\eeqq  
(in virtue of theorem \ref{1.3.10}, such cubes do exist).
Fix $0<\la<\frac 1 2$ and consider asymptotics $\Delta \downarrow 0$ as in theorem \ref{2.2.4}; introduce the sequence of integers 
\beqq\label{n_delta_lambda_NEU}
n(\Delta) \;:=\; \lfloor\, L(A)\; \Delta^{-\frac{1}{2d}} \rfloor   
\eeqq
increasing to $\infty$ as $\Delta \downarrow 0$, with $L(A)$ the edge length of $A$ selected in \eqref{M1NEU}.

Writing for short $n=n(\Delta)$, we partition the cube $A$ in \eqref{M1NEU} into $\,n^d \sim {\rm vol}(A)  \Delta^{-\frac12}\,$ cells of equal volume $\,\sim  \Delta^{\frac12}\,$ and of equal edge length $\,\sim \Delta^{\frac{1}{2d}}\,$  
in every component. 
In the special case where $A$ is the unit cube  $A := [0,1]^d$, every such cell is identified through its upper right corner $(\frac{j_1}{n},\ldots,\frac{j_d}{n})$. Then we write $\calj(\Delta)$ for the set of all multiindices $\al=(j_1,\ldots,j_d)$ arising in this way, and $A_\al$ for the cell whose upper right corner $(\frac{j_1}{n},\ldots,\frac{j_d}{n})$ makes appear $\al\in\calj(\Delta)$. 
For a general cube $A$ selected in \eqref{M1NEU}, a linear transformation component by component maps $A$ to $[0,1]^d$, allows to identify cells $A_\al$ by the image $(\frac{j_1}{n},\ldots,\frac{j_d}{n})$ of their upper right corner, and thus again yields a represention of $A$ as a union of cells $A_\al$, $\al\in \calj(\Delta)$. 
From \eqref{n_delta_lambda_NEU} we have a one-to-one correspondence between $\calj(\Delta)$ and $\{1,\ldots,n\}^d$ where $n=n(\Delta)$. 
We shall make use of the decomposition of  $A$ meeting \eqref{M1NEU} into cells 
$$
A_\al  \quad,\quad  \al\in \calj(\Delta)
$$
to fill from discrete BDI observations $(\eta_{i\Delta})_{i\in\bbn_0}$ regression schemes, see \ref{3.1.1NEU} below. 
Upon careful choice of $0<\la<\frac 12$ in definitions  \ref{2.2.1} and \ref{2.2.3} and thus in the reconstruction algorithm \ref{2.2.2}, such schemes will allow  to estimate the diffusion coefficient \eqref{invertible_diff_coeffNEU} of the one-particle motion. 
\\

\begin{defn}\label{3.1.1NEU} (Regression scheme)\,
Fix a cube $A$ meeting \eqref{M1NEU}. Fix $0<\la<\frac12$  and let $\Delta$ decrease to $0$. For  $\Delta$ small enough, define pairs   
$$
( \calx_{\al} , \calz_{\al} ) \;:\; \al\in \calj(\Delta)
$$
as follows: \\
i) For  $\al\in \calj(\Delta)$, define
$$
\tau_\al=\tau_\al (\Delta)\;:=\; \inf\left\{ i\in\bbn_0 : \mbox{$\;(\eta_{ i\Delta }, \eta_{(i{+}1)\Delta })\,$  is  $(\Delta,\la)$-identifiable  and  satisfies $\eta_{i\Delta}(A_\al)\ge 1$ } \right\} \;. 
$$
At time $\tau_\al\Delta$, writing for short 
$$
x:=\eta_{ \tau_\al\Delta } \;,\; 
y:=\eta_{(\tau_\al{+}1)\Delta }\;,\; 
\ell:=\ell(x)=\ell(y)\ge 1\;,\; 
x=(x_1,\ldots,x_\ell)\;,\; 
y=(y_1,\ldots,y_\ell)
$$
and $\pi$ for the unique permutation of $(1,2,\ldots,\ell)$ such that  $y_{\pi(i)}$ is close to $x_i$ for all $1\le i\le\ell$  
$$
\left| y_{\pi(i),j} - x_{i,j}\right| \;<\; \Delta^\la \quad\mbox{for all $1\le i\le \ell$, $j=1,\ldots, d$ }  
$$
in the sense of definition \ref{2.2.1} (the permutation is trivial in case $\ell=1$). Then $\eta_{i\Delta}(A_\al)\ge 1$ allows to pick  $m=m(\al)$ such that particle $x_{m(\al)}$ is located in the cell $A_\al$ at time $\tau_\al\Delta$. For this $m(\al)$ we define  
$$
\calx_{\al} := x_{m(\al)}  \quad\mbox{together with}\quad \calz_{\al} := \frac{ y_{\pi(m(\al))}- x_{m(\al)} }{ \sqrt{\Delta} }  \;.  
$$
ii) In order to do so for all cells $A_\al$, $\al\in \calj(\Delta)$, we define 
$$
\tau^*=\tau^*(\Delta)  \;=\;  \max \{ \tau_\al : \al\in \calj(\Delta)\} \;. 
$$ 
\end{defn}

\vskip0.5cm
We make some comments. First, $\,\tau_\al$ is a stopping time with respect to the filtration $(\calh^\Delta_{i+1})_{i\in\bbn_0}$ defined after definition \ref{2.2.2} of the reconstruction algorithm; the same holds for $\tau^*$. 
We have to show that these stopping times are almost surely finite (at least), then $\calx_\al$ and $\calz_\al$ are well-defined random variables taking values in $E=\bbr^d$. 
Second, it may occur that we fill disjoint cells $A_\al \neq A_{\al'}$ simultaneously at the same time.
When $\,\om\in\{\tau_\al=\tau_{\al'}\}\,$, pairs $(x,y)$ defined by 
$$
x=\eta_{ \tau_\al\Delta}(\om)=\eta_{ \tau_{\al'}\Delta}(\om) \;,\; y=\eta_{(\tau_\al{+}1)\Delta}(\om) =\eta_{(\tau_{\al'}{+}1)\Delta}(\om)
$$
being $(\Delta,\la)$-identifiable in the sense of definition \ref{2.2.1}, we must have the following: 
$m(\al)\neq m(\al')$ since $A_\al \neq A_{\al'}$, thus $\pi(m(\al)) \neq \pi(m(\al'))$ since $\pi$ in  \ref{2.2.1} is a permutation,  
thus $x_{m(\al)}\neq x_{m(\al')}$ and $y_{\pi(m(\al))}\neq y_{\pi(m(\al'))}$ also in this case.\\

\begin{prop}\label{3.1.3''NEU}
In the framework of definition \ref{3.1.1NEU}, the stopping times $\tau^*=\tau^*(\Delta)$ have finite expectation. 
\end{prop}

The proof given in section \ref{3.3NEU} below will also show that the expected time  $\,E_\mu\left( \Delta\,  \tau^*(\Delta) \right)\,$ which we need to fill the scheme \ref{3.1.1NEU} remains bounded as $\Delta\downarrow 0$.  \\

\begin{thm}\label{3.1.4NEU}
In the framework of  \ref{3.1.1NEU}, the scheme 
$$
( \calx_\al , \calz_\al ) \;:\; \al\in\calj(\Delta)
$$
has the following  properties:

i) We have $\calx_\al\in A_\al$ for all $\al\in\calj(\Delta)$, thus  
$$
\mbox{ variables  $\calx_\al$, $\al\in\calj(\Delta)$, are approximately regularly spaced in the cube $A$ } \;.  
$$
ii) There are `good events' $\,{\tt G}(\Delta)$  of probability $\ge 1-\calo( \Delta^{\frac12} )$ on which the full collection  $\calz_\al$, $\al\in\calj(\Delta)$, represents true (rescaled) $\Delta$-increments of the underlying one-particle motion in  
{\bf (A1)}. 

More precisely, on the `good event' $\,{\tt G}(\Delta)$, there is a collection   $W_\al$ of independent $d$-dimensional standard Brownian motions and a collection $\xi_\al$ of path segments which are solutions to       
$$
d\xi_\al(s) \;=\; b(\xi_\al(s)) ds + \si(\xi_\al(s))dW_\al(s)  
\quad,\quad  0 \le s\le \Delta
\quad,\quad \al\in\calj(\Delta)
$$
(strong solutions driven by $W_\al$)   such that all pairs $(\calx_\al,\calz_\al)$, $\al\in\calj(\Delta)$, admit the representation 
\beqq\label{rescaled_increments_GnDeltaNEU}
\frac{ \xi_\al(\Delta) - \xi_\al(0) }{ \sqrt{\Delta} }  \;=\; \calz_\al 
\quad,\quad 
\xi_\al(0) \;=\; \calx_\al   \quad,\quad \al\in\calj(\Delta)  \;. 
\eeqq

iii) There are exceptional events  $\,{\tt F}(\Delta)$  of probability  
$$
Q_\mu\left(\,{\tt F}(\Delta)\,\right) \;\le\; \calo(\Delta^{\frac12})
$$
on which at least one entry $\calz_\al$ to the scheme \ref{3.1.1NEU}, $\al\in\calj(\Delta)$, fails to represent a true (rescaled) increment of the underlying one-particle motion. 

iv) By construction in \ref{3.1.1NEU} we have   
\beqq\label{deterministic_bounds_znalNEU}
\calz_\al \;\in\; \mathop{\sf X}\limits_{i=1}^d \left( -\Delta^{\la-1/2} \,,\, \Delta^{\la-1/2} \right) 
\quad\mbox{for all}\quad \al\in\calj(\Delta)  \;.
\eeqq
\end{thm}

\vskip0.8cm
The proof will be given in section \ref{3.3NEU}. \\

We intend to use the scheme  \ref{3.1.1NEU} as a regression scheme for nonparametric estimation of the diffusion coefficient 
of the one-particle motion in  {\bf (A1)},  based on `high-frequency'  discrete BDI observations. The motivation is the following (see proposition 6 in Genon-Catalot and Jacod \cite{GJ-93}, approximation (16) in Podolskij and Vetter \cite{PV-10}, or the book Jacod and Protter \cite{JP-12}): on the good events $\,{\tt G}(\Delta)$ in theorem \ref{3.1.4NEU}, good approximations  
\beqq\label{good_approximationNEU}
\calz_\al  \;=\; \frac{ \xi_\al(\Delta) - \xi_\al(0) }{ \sqrt{\Delta} }  \;\;\approx\;\; 
\si(\xi_\al(0))\; \frac{ W_\al(\Delta) - W_\al(0) }{ \sqrt{\Delta} }
\;=\; \si(\calx_\al)\; \calu_\al(1)
\eeqq
are available for the terms in \eqref{rescaled_increments_GnDeltaNEU} where by construction 
$$
\calu_\al(s) := \frac{ W_\al(s\Delta) - W_\al(0) }{ \sqrt{\Delta} } \;,\; 0\le s\le 1 \;,\;\al\in\calj(\Delta)
$$
are independent $d$-dimensional standard Brownian motions on the time interval $[0,1]$.   
Clearly $W_\al$ is independent of $\xi_\al(0)=\calx_\al$. Since increments of the same Brownian motion over time intervals $[i\Delta,(i+1)\Delta]$, $[i'\Delta,(i'+1)\Delta]$, $i'\neq i$, are independent, and since different particles have independent driving Brownian motions, the construction in theorem \ref{3.1.4NEU} grants that --in restriction to the good events $\,{\tt G}(\Delta)$-- the collection of Brownian motions $\{\calu_\al:\al\in\calj(\Delta)\}$ is independent of the collection of design variables $\{\calx_{\al'}:\al'\in\calj(\Delta)\}$. 

From It\^{o}'s formula, using superscript $i$ for the components of $\calu_\al$ and $\delta^{(i,j)}$ for Kronecker's $\delta$, 
$$
\calu_\al(1)\calu_\al^\top(1) \;=\; \left( 
\delta^{(i,j)} + \int_0^1 \calu_\al^{(i)}(s)\, d\calu_\al^{(j)}(s) + \int_0^1 \calu_\al^{(j)}(s)\, d\calu_\al^{(i)}(s) 
\right)_{1\le i,j\le d}
$$
which means that, always on the `good event' $\,{\tt G}(\Delta)$ in theorem \ref{3.1.4NEU},   
\beqq\label{interpretation_schemeNEU}
\calz_\al\calz_\al^\top \;=\; (\si\si^\top)(\calx_\al) \;+\;\;\mbox{error terms of martingale structure} \;. 
\eeqq
Thanks to theorem \ref{3.1.4NEU} and \eqref{interpretation_schemeNEU},  \ref{3.1.1NEU} provides us --in restriction to the good sets $\,{\tt G}(\Delta)\,$-- with 
a regression scheme --in a sense analogous to sections 1.5.1 or 1.6.1 of Tsybakov \cite{Tsy-08}-- for nonparametric estimation of the diffusion coefficient $a(\cdot)=(\si\si^\top)(\cdot)$  on ${\rm int}(A)$. 
The design variables $\calx_\al$, $\al\in\calj(\Delta)$ are approximately regularly spaced over the cube $A$ selected in the beginning.  

In contrast to the good sets however, the picture is less pleasant on the exceptional sets $\,{\tt F}(\Delta)\,$: here we have nothing except the trivial bounds from theorem \ref{3.1.4NEU}~iv). 
We shall illustrate the effect of the exceptional sets by an example (kernel estimation in dimension $d=1$) in section \ref{3.2NEU}.\\

It remains to make precise what `good approximation' in \eqref{good_approximationNEU} means.  We quote a result from Genon-Catalot and Jacod \cite{GJ-93}; by theorem \ref{3.1.4NEU}, their result applies in restriction to the good set $\,{\tt G}(\Delta)$ where all $\calz_\al$ in \eqref{rescaled_increments_GnDeltaNEU}, $\al\in\calj(\Delta)$, are increments of rescaled one-particle motions. 
Below, $g$ denotes a polynomial on $\bbr^d$ of arbitrary finite degree $\gamma$. 
\\

\begin{lemma}\label{3.1.5NEU}
(\cite{GJ-93}, proposition 6):  Assume that  the coefficients of the one-particle motion in {\bf (A1)} are $\,\calc^2$ on $E=\bbr^d$ and satisfy \eqref{invertible_diff_coeffNEU}  together with 
\beqq\label{additional_assumption_J+GCNEU} 
|b_l| \;,\;  | \sigma_{l,j} |  \;,\; | \partial_i \sigma_{l,m} |   \;,\; | \partial_i \partial_j \sigma_{l,m} |    
\;\;\mbox{are bounded by some constant $L$} \;.   
\eeqq 
Then, using notations of  theorem \ref{3.1.4NEU}, \eqref{rescaled_increments_GnDeltaNEU}  and \eqref{good_approximationNEU}, there is some constant $C$  
(which depends on $L$ and on the degree $\gamma$ of the polynomial $g$, but does not depend on $\Delta$ as $\Delta\downarrow 0$) 
such that the following deterministic bound holds:   
\beqq\label{GJ-93_proposition_6NEU}
E\left(\; 1_{{\tt G}(\Delta)}\! \left[\, g\left( \calz_\al \right) - g\left( \,\si(\calx_\al)\, U_\al(1)\, \right) \,\right]^2 
 \,\middle|\, 
\{\calx_{\al'}:\al'\in\calj(\Delta)\} 
\right) \;\;\le\;\; C\, \Delta
\quad,\quad \al\in\calj(\Delta)   \;. 
\eeqq
\end{lemma}
\vskip0.8cm

\subsection{Proofs for section \ref{3.1NEU}}\label{3.3NEU}

We start with some remarks motivating the construction in \ref{3.1.1NEU}. 
Positive Harris recurrence of the continuous-time process grants that $(\eta_t)_{t\ge 0}$ visits infinitely often the void configuration $\delta$, spending an exponentially distibuted time in $\delta$ at every visit. 
Thus for fixed $\Delta>0$, discrete observations $\,(\eta_{i\Delta})_{i\in\bbn_0}$ form a positive recurrent Markov chain, with recurrent atom $\{\delta\}$ and invariant measure $\mu$ on $S$, and pairs of successive observations $(\eta_{i\Delta}, \eta_{(i+1)\Delta})_{i\in\bbn_0}$ form a positive recurrent Harris chain, with recurrent atom $\{(\delta,\delta)\}$ and invariant measure $\mu(dx) P_{\Delta}(x,dx')$  on $S\times S$. 
Positive recurrence ensures that in the long run we will collect an overwhelming amount of data: 
out of these we pick few but well-selected ones --using the reconstruction algorithm \ref{2.2.2}-- to fill the  scheme \ref{3.1.1NEU}. \\

\begin{prop}\label{prop_4.2.1NEU}
For the cube $A$ meeting \eqref{M1NEU} decomposed into disjoint cells $A_\al$, $a\in\calj(\Delta)$, 
$$
\al \;\lra\; \frac{1}{ \ov\mu(A) } \int_S \mu(dx)\, 1_{\{ x(A_\al)\ge 1\}} 
$$
defines a probability law on the finite set $\calj(\Delta)$ which is equivalent to the uniform law on  $\calj(\Delta)$. 
\end{prop}

{\bf Proof: } From theorem \ref{1.3.10}, the invariant occupation measure $\ov\mu$ on the single particle space  $E=\bbr^d$ admits a Lebesgue density $\ov \gamma \in \calc_0(E)$. 
By choice of the cube $A$ in \eqref{M1NEU} we have $\ov \gamma >0$ on $A$, hence 
$$
\al \;\lra\;  \frac{1}{ \ov\mu(A) } \int_S \mu(dx)\, x(A_\al)  \;=\; 
\frac{ \ov\mu(A_\al) }{ \ov\mu(A) } \;=\; 
\frac{ 1 }{ \ov\mu(A) } \int_E dy\; \ov\gamma(y)\, 1_{A_{\al}}(y)     \quad>\quad 0
$$
defines a probability on  $\calj(\Delta)$ under which every $\al\in\calj(\Delta)$ carries strictly positive mass. In the sense of equivalence of measures, this probability is equivalent to the uniform law on  $\calj(\Delta)$.  Also 
$$
\al \;\lra\; \int_S \mu(dx)\, x(A_\al)  \quad,\quad   \al \;\lra\; \int_S \mu(dx)\, 1_{\{ x(A_\al)\ge 1\}}  
$$
are equivalent measures on  $\calj(\Delta)$. \halmos\\

As an application of theorem \ref{2.2.4}, positive Harris recurrence of $(\eta_{i\Delta})_{i\in\bbn_0}$ combined with the last proposition yields a proof that the stopping times $\tau^*(\Delta)$ in  proposition \ref{3.1.3''NEU} have finite expectation:  \\

{\bf Proof of proposition \ref{3.1.3''NEU}: } 
For $\Delta>0$ fixed, \eqref{M1NEU} implies that we need at most a geometric number of life cycles to observe the first  occurrence of $\{\eta_{i\Delta}(A)\ge 1\}$, $i\in\bbn_0$. The expected length of a life cycle is finite. For $\al\in\calj(\Delta)$ fixed, \ref{prop_4.2.1NEU} grants that we need at most a geometric number of occurrences of $\{\eta_{i\Delta}(A)\ge 1\}$ to record the first occurrence of $\{\eta_{{i'}\Delta}(A_\al)\ge 1\}$, ${i'}\in\bbn_0$. 
Here the success probability of the geometric number of trials is $\calo(n^{-d})=\calo(\Delta^{\frac12})$, as a consequence of proposition \ref{prop_4.2.1NEU}. For $\Delta$ small enough, theorem \ref{2.2.4} ensures that after at most a geometric number of occurrences of $\{\eta_{{i'}\Delta}(A_\al)\ge 1\}$ we will record the first occurrence of $\{\eta_{ {i''}\Delta}(A_\al)\ge 1 \,,\, \eta_{[{i''}\Delta , ({i''}+1)\Delta]}\in{\tt CI}(\Delta,\la) \}$, ${i''}\in\bbn_0$. By \eqref{fondamental_property_bis}, this is an occurrence (not necessarily the first 
one) of the desired event 
$$
\left\{\, \eta_{ {i''}\Delta}(A_\al)\ge 1 \,,\, \left(\eta_{{i''}\Delta} , \eta_{({i''}+1)\Delta}\right)\in{\tt ID}(\Delta,\la) \,\right\} \;. 
$$
Since $\calj(\Delta)$ is a finite set, proposition  \ref{3.1.3''NEU}  is proved. \halmos\\

We remark that this proof indicates the following: Since $\calj(\Delta)$ has $n^d$ elements, since for every element $\al$ of $\calj(\Delta)$ we need in $Q_\mu$-expection $\calo(n^d)$ trials to hit $A_\al$, the expected time $\,E_\mu( \Delta\, \tau^*(\Delta))\,$ which we need to fill the scheme \ref{3.1.1NEU} will be of order 
$$
\calo( n^d n^d \Delta) \;=\; \calo(1)
$$
as $\Delta\downarrow 0$. This is the motivation for the choice of $n=n(\Delta)$ in \eqref{n_delta_lambda_NEU}. \\

{\bf Proof of theorem \ref{3.1.4NEU}: } 
For every  $\Delta$ as $\Delta\downarrow 0$, the construction in \ref{3.1.1NEU} uses a fixed cube $A$ meeting \eqref{M1NEU}, partitioned into  $n^d=\calo(\Delta^{-\frac 12})$ cells $A_\al$ of equal size, $\,\al\in\calj(\Delta)$, and pairs 
\beqq\label{scheme_final_representationNEU}
\left(\eta_{\tau_\al\Delta} , \eta_{(\tau_\al+1)\Delta}\right)\in{\tt ID}(\Delta,\la) \quad\mbox{such that}\quad \eta_{\tau_\al\Delta}(A_\al)\ge 1 \quad,\quad  \al\in \calj(\Delta) 
\eeqq
where $\,n=n(\Delta)$ is given by  \eqref{n_delta_lambda_NEU}.

1) We prove assertions iv) and i) of the theorem. By definition in \ref{2.2.1}, both configurations in \eqref{scheme_final_representationNEU} have equal length $\ell\ge 1$. 
For every $\al\in\calj(\Delta)$,  
a unique permutation $\pi$ of $(1,\ldots,\ell)$ maps 
particles in the configuration $\eta_{(\tau_\al+1)\Delta}=:y=(y_1,\ldots,y_\ell)$ to particles in the configuration $\eta_{\tau_\al\Delta}=:x=(x_1,\ldots,x_\ell)$ via  
$$
| y_{\pi(i),j} - x_{i,j} | <  \Delta^\la \quad\mbox{for all $1\le i\le \ell$ and all $1\le j\le d$} \;. 
$$
Since $\eta_{\tau_\al\Delta}(A_\al)\ge 1$, some particle $m=m(\al)$ in the configuration $x$ is visiting $A_\al$ at time $\tau_\al\Delta$: with
$$
\calx_\al := x_m \quad,\quad \calz_\al := \frac{ y_{\pi(m)} - x_m }{\sqrt{\Delta}}
$$
we achieve $\calx_\al\in A_\al$ and $\calz_\al \in \mathop{\sf X}\limits_{j=1}^d \left( -\Delta^{\la-1/2} , \Delta^{\la-1/2} \right) $.

2) To prove assertions ii) and iii), we define the `good events' in  ii) as follows: 
$$
{\tt G}(\Delta) \;:=\;\left\{\,
\mbox{all  pairs} \;\; 
(\eta_{\tau_\al\Delta},\eta_{(\tau_\al{+}1)\Delta}) \;\;\mbox{in \eqref{scheme_final_representationNEU} are such that}\;\; \eta_{[\tau_\al\Delta,(\tau_\al{+}1)\Delta]}
 \in {\tt CI(\Delta,\la)}   
\,\right\} \;. 
$$
Then definition \ref{2.2.3} and \eqref{fondamental_property_bis}+\eqref{fondamental_property} show the following:  in restriction to  ${\tt G}(\Delta)$, the reconstruction algorithm \ref{2.2.2} applied to data \eqref{scheme_final_representationNEU} will reconstruct all particle identities correctly.
In particular it reconstructs correctly --on ${\tt G}(\Delta)$, on $\{\tau_\al=i\}$ for  $\al\in\calj(\Delta)$-- the increments $\,\xi_\al((i{+}1)\Delta)-\xi_\al(i\Delta)\,$ in the trajectory of the selected particle over the time interval $[i\Delta, (i{+}1)\Delta]$. Uniquely associated to this particle and this time interval is the increment $\,W_\al((i{+}1)\Delta)-W_\al(i\Delta)\,$ of the driving Brownian motion in {\bf (A1)}. Driving Brownian motions for different particles are independent 
(we might select the same particle twice: then this happens over disjoint time intervals $[i\Delta, (i{+}1)\Delta]$, $[i'\Delta, (i'{+}1)\Delta]$, $i'\neq i$, and we have again independence of the increments of the Brownian motion). This is \eqref{rescaled_increments_GnDeltaNEU}, up to a change of time. 

3) It remains to prove the bound in assertion iii). 
Note that  ${\tt G}(\Delta)$ belongs to the $\si$-field $\calf_{(\tau^*+1)\Delta}$ associated to continuous-time observation  up to time $(\tau^*+1)\Delta$. Define  ${\tt F}(\Delta)$ as the event that the scheme \eqref{scheme_final_representationNEU} will involve some pair $(\eta_{i\Delta},\eta_{(i{+}1)\Delta})$ for which the segment $\eta_{[i\Delta,(i{+}1)\Delta]}$ lacks to be $(\Delta,\la)$-good. 
For $i$ fixed, as a consequence of \ref{2.2.4} and stationarity, the probability under $Q_\mu$ to have 
$$
\left(\eta_{i\Delta},\eta_{(i{+}1)\Delta}\right) \in {\tt ID}(\Delta,\la) \quad,\quad \eta_{[i\Delta,(i{+}1)\Delta]} \notin {\tt CI}(\Delta,\la)
$$ 
is bounded by $D\Delta$ for all $0<\Delta<\Delta_0$. Here $D$ is some constant $D$ which does not depend on $\Delta$. 
The cube $A$ being partitioned into  $n^d=\calo(\Delta^{-\frac12})$ cells $A_\al$, $\,\al\in\calj(\Delta)$, $\,n=n(\Delta)$ as in  \eqref{n_delta_lambda_NEU}, we need $n^d$ pairs $(\eta_{i\Delta},\eta_{(i{+}1)\Delta})$ to fill the regression scheme \ref{3.1.1NEU}. This gives 
$$
Q_\mu\left( {\tt F}(\Delta) \right) \;\le\; 
\calo(n^d \Delta) \;=\; \calo(\Delta^{\frac12}) \;. 
$$
Thus $\,{\tt F}(\Delta)$ is an exceptional event in a sense of vanishing probability as $\Delta\downarrow 0$  . \halmos\\

\subsection{ Example: Kernel estimators for the diffusion coefficient in dimension $d=1$
}\label{3.2NEU} 


We restrict the setting to dimension $d=1$. Based on data
$\,(\eta_{i\Delta})_{i\in\bbn_0}\,$ from discrete observation 
at step size $\Delta$ and on a regression scheme \ref{3.1.1NEU} filled from these data, 
we wish to construct a kernel estimator for the diffusion coefficient $\si^2(\cdot)$. Calculating its pointwise risk under squared loss, we have to make sure that the influence of the exceptional sets in theorem \ref{3.1.4NEU}~iv) does not dominate the (classical) bounds which do hold on the good sets: 
this will oblige us to select $\la$ quite close to $\frac 12$. Recall that $0<\la<\frac{1}{2}$ remains fixed in \ref{2.2.1}, \ref{2.2.3} and \ref{3.1.1NEU} while asymptotics is in $\Delta\downarrow 0$. 

We make use of the following notations. For $\wt\beta\in\bbn$, a kernel of order $\wt\beta$  (\cite{Tsy-08} pp.\ 5, 10) is a function $K:\bbr\to\bbr$ supported by $[-1,1]$ and Lipschitz continuous on $(-1,1)$ such that   
$$
\int K(v)\, dv \;=\; 1 \;,\;  
\int v^r\, K(v)\, dv \;=\; 0 \quad\mbox{for $r=1,\ldots, \wt\beta$}  \;. 
$$
For any choice $h>0$ of a bandwidth we write $\,K_h(v):=\frac1h K(\frac v h)\,$.

For $\beta>1$ let $\beta'$ denote the largest natural number which is strictly smaller than $\beta$. The H\"older class $\calh(\beta,L)$  of order $\beta$ (\cite{Tsy-08} p.\ 5) is the class of all functions $f$ in $\calc^{\beta'}$ with the following property: the derivative $f^{(\beta')}$ of order $\beta'$ is H\"older with index  $\,\beta{-}\beta' \in (0,1]$ and with H\"older constant $L$: 
$$
\left| f^{(\beta')}(x) - f^{(\beta')}(x') \right| \;\le\; L\, |x-x'|^{\beta-\beta'} \quad,\quad x,x'\in\bbr \;. 
$$
An example: whenever in dimension $d=1$ the one-particle motion in {\bf (A1)} satisfies \eqref{invertible_diff_coeffNEU} together with \eqref{additional_assumption_J+GCNEU}, the diffusion coefficient $\si^2:\bbr\to (0,\infty)$ belongs to the H\"older class $\calh(2,L)$. When we model BDI in dimension $d=1$ we might have reasons to work with diffusion coefficients which have `more smoothness': our statistical model below will assume that $\si^2(\cdot)$ belongs to a H\"older class  $\,\calh(\beta,L)$ of order $\beta\ge 2$, with $\beta$ fixed and known. To every  $\beta\ge 2$ we associate a critical value    
\beqq\label{critical_valueNEU}
\la_0(\beta) \;:=\; \frac 12 - \frac{1}{8(2\beta+1)} \;=\; \frac{8\beta + 3}{16\beta + 8} \quad\in\;\;  (\,0\,,\,\tfrac 12 \,)   \;. 
\eeqq

\vskip0.5cm
In dimension $d=1$, the cube $A$ meeting \eqref{M1NEU} in the beginning of section \ref{3.1NEU} is an interval, decomposed according to \eqref{n_delta_lambda_NEU} into   
\beqq\label{n_delta_lambda_d1NEU}
n(\Delta) \;=\; \lfloor {\rm length}(A)\; \Delta^{-\frac12} \rfloor   
\eeqq
disjoint subintervals $\,A_\al\,$ of equal length $\,\sim \Delta^{\frac 12}\,$, $\al\in\calj(\Delta)$.  For given order $\beta\ge 2$, it is well known that one needs kernels of order depending on $\beta$ to estimate $\si^2(\cdot)$ within class $\,\calh(\beta,L)$ at the optimal nonparametric rate. In our case, exceptional sets being present in theorem \ref{3.1.4NEU}, we need more: 
we also have to relate $\,\la\,$ in \ref{2.2.1}, \ref{2.2.3} and \ref{3.1.1NEU} to the given order $\beta$, and have to work with 
\beqq\label{critical_value_bedingungNEU}
\la_0(\beta) \;\le\; \la\;<\; \frac12  
\eeqq
with $\la_0(\beta)$ from \eqref{critical_valueNEU}. Then, from discrete observation of the BDI process at step size $\Delta$, we fill regression schemes \ref{3.1.1NEU} where $\la$ is fixed according to \eqref{critical_value_bedingungNEU} 
and where $\Delta$ tends to $0$: at stage $\Delta$ of the asymptotics, we work with  
$$
\left( \calx_\al , \calz_\al \right) \;,\; \al\in\calj(\Delta)    
$$
(and shall sometimes write $(\calx_\al,\calz_\al)_{1\le \al\le n}$). Next, we 
fix a kernel $K$ of order $\beta'$, the largest natural number strictly smaller than $\beta$, define the bandwidth $h=h(\Delta)$ by 
\beqq\label{bandwidthNEU}
h(\Delta) \;:=\; (n(\Delta))^{-\frac{1}{2\beta+1}} \;=\; \calo( \Delta^{ \frac12 \frac{1}{2\beta+1} } )   
\eeqq
and introduce the estimator  
\beqq\label{def_kernelestNEU}
\wh{\si^2_{\Delta}}(a)  \;:=\; \sum\limits_{\al\in\calj(\Delta)} {\rm length}(A_\al)\,  \,\caly_\al\;\, K_h(\calx_\al-a)  
\quad,\quad \caly_\al := \calz_\al^2
  \quad,\quad a \in {\rm int}(A) \;. 
\eeqq
Recall that all cells $A_\al$ are of equal length $\,\sim \Delta^{\frac 12}\;$ (which of course is $\,\calo(\frac 1n)\,$ as $\,\Delta\downarrow 0\,$, but note:   $\,{\rm length}(A_\al)$ asymptotically does not depend  on the size of the interval $A$ selected in \eqref{M1NEU}, whereas $\,\frac{1}{n(\Delta)}=\frac 1n\,$ does).
In this setting, we give a rate as $\,\Delta\downarrow 0\,$ for the risk of this estimator, writing $E_{\si^2}(\ldots)$, $Q_{\si^2}(\ldots)$ instead of $E_\mu(\ldots)$, $Q_\mu(\ldots)$ in order to stress the dependence on $\si^2$.\\

\begin{thm}\label{3.2.2'NEU}
For $\beta\ge2$ fixed, let the diffusion coefficient in {\bf (A1)} belong to class $\calh(\beta,L)$. Choose $0<\la<\frac12$ sufficiently close to $\frac12$ to satisfy  \eqref{critical_value_bedingungNEU}. Then asymptotically as $\Delta\downarrow 0$, with $n=n(\Delta)$ from  \eqref{n_delta_lambda_d1NEU},  $h=h(\Delta)$ from \eqref{bandwidthNEU}, and $K$ of order $\beta'$, 
the pointwise risk of the kernel estimator (\ref{def_kernelestNEU}) under squared loss satisfies 
\beqq\label{rate_squared_errorsNEU}
\limsup_{\Delta\downarrow 0}\;\; n^{\frac{2\beta}{2\beta+1}}\; 
E_{\si^2}\left( \left| \wh{\si^2_\Delta}(a) - \si^2(a) \right|^2  \right)  
\;\;<\;\;  \infty  
\eeqq
at every point $a\in{\rm int}(A)$. 
\end{thm}

\vskip0.5cm
The theorem will be proved in subsection \ref{3.4NEU} below. Note that the rate which appears in (\ref{rate_squared_errorsNEU}) is the nonparametric rate which is known to be optimal (see Tsybakov \cite{Tsy-08} section 2.5) under squared loss in standard regression schemes 
$ 
\;(\calx'_\al,\caly_\al)  \;,\;  \caly_\al=f(\calx'_\al)+\vep_\al  \;,\; 1\le \al\le n\;
$
with unknown $f\in\calh(\beta,L)$, in the classical setting of i.i.d.\ square-integrable errors $\vep_\al$ and equispaced deterministic design points~$\calx'_\al$.

\subsection{Proving theorem \ref{3.2.2'NEU}}\label{3.4NEU}


We prepare the proof of theorem \ref{3.2.2'NEU} by a series of auxiliary steps. The proof will need the diffusion coefficient $\si^2$ only in restriction to compacts like $\{ y\in\bbr : d(y,A)\le 1 \}$, and we shall  write again $L$ --as before in assumption \eqref{additional_assumption_J+GCNEU}-- to denote a bound on this compact for derivatives of the function $\,\si^2 \in \calh(\beta,L)$ up to order $\beta'$. Assumptions and nota\-tions are those of \ref{3.1.1NEU}, of theorem \ref{3.1.4NEU} together with lemma \ref{3.1.5NEU}, and of \ref{3.2.2'NEU}. Often we write for short $n$, $h$ instead of $n(\Delta)$, $h(\Delta)$ etc.  \\

\begin{lemma}\label{3.3.2NEU}
As $\Delta\downarrow 0$, we have deterministic bounds 
$$
\sup_{a \in {\rm int}(A)}\;\;  \limsup_{\Delta\downarrow 0}\quad
nh \left|\; \sum\limits_{\al\in\calj(\Delta)} {\rm length}(A_\al)\, K_h( \calx_\al-a ) \;-\; 1 \;\right|  
\quad\le\;\; M \;\;<\quad  \infty   
$$  
at rate 
\beqq\label{nhrateNEU}
n h  \;=\;  n^{ \frac{2\beta}{2\beta + 1} }  \;=\;   \calo( \Delta^{ \frac{-\beta}{2\beta+1} } ) \quad\lra\quad\infty \;. 
\eeqq
Recall that all cells $A_\al$ have equal length $\,\sim \Delta^{\frac 12} = \calo(\frac 1n)\,$ for all $\al\in\calj(\Delta)$.    
\end{lemma}

\vskip0.5cm
{\bf Proof: }  Fix $a\in{\rm int}(A)$. By continuity of $K$ on $(-1,1)$, whenever $\Delta$ is small enough and a cell $A_\al$  fully contained in $B_h(a)$, the ball of radius $h$ around $a$, we can select a point $\zeta_\al \in {\rm cl}(A_\al)$ such that 
$$
\int_{A_\al} K_h(x-a)\, dx \;=\; {\rm length}(A_\al)\; K_h(\zeta_\al -a)   \;\sim\; \Delta^{\frac 12}\, K_h(\zeta_\al -a) \;. 
$$
Now  $K$ is Lipschitz on $(-1,1)$, thus 
$$
\left|\, K\Big( \frac{\zeta_\al-a}{h} \Big) - K\Big( \frac{\calx_\al-a}{h} \Big) \,\right| \;\;\le\;\; \calo\left( {\rm length}(A_\al) \frac 1 h \right)\; \;=\;\calo\left( \frac{1}{nh} \right)    
$$
since $\calx_\al\in A_\al$ and by definition of $n=n(\Delta)$ in \eqref{n_delta_lambda_d1NEU}, thus 
$$
\int_{A_\al} K_h(x-a)\, dx \;=\; {\rm length}(A_\al)\; K_h(\calx_\al -a) \;+\; \calo\left(\frac{1}{(nh)^2}\right)  \;. 
$$
For $\calo(nh)$ indices $\al$ in $\calj(\Delta)$ the cell $A_\al$ will be fully contained in $B_h(a)$; for at most two additional values of $\al$, a cell $A_\al$ may intersect $B_h(a)$. Since $a\in{\rm int}(A)$, we have for $\Delta$ small enough 
$$
1 \;=\;  \int_{B_h(a)} K_h(x-a)\, dx \;=\; \sum_{\al\in\calj(\Delta)} \int_{A_\al} K_h(x-a)\, dx
\;=\; \sum_{\al\in\calj(\Delta)}{\rm length}(A_\al)\; K_h(\calx_\al -a) 
\;+\; \calo\left(\frac{1}{nh}\right)  
$$
which is the assertion. 
By  \eqref{n_delta_lambda_d1NEU} and \eqref{bandwidthNEU}, $\;n h\,$  
tends to $\infty$ as $\Delta\downarrow 0$.  \halmos\\

\begin{lemma}\label{3.3.1NEU}
i) For all powers $m\in\bbn_0$ we have deterministic bounds 
\beqq\label{*NEU}
\sup_{a \in {\rm int}(A)}\;\;  \limsup_{\Delta\downarrow 0}\quad
\sum\limits_{\al\in\calj(\Delta)} {\rm length}(A_\al)\;  \frac{1}{h}\,|K|^m(\tfrac{\calx_{\al}-a}{h})  \quad\le\;\; M\;\; <\quad \infty 
\eeqq
with suitable constants $M=M(K,m)$. \\ 
ii) For continuous functions $f:\bbr\to\bbr$ we have  
$$
\sum\limits_{\al\in\calj(\Delta)} {\rm length}(A_\al)\; f(\calx_\al)\; K_h(\calx_\al-a)  \;\;\lra\;\;  f(a)  \quad,\quad a \in {\rm int}(A)  
$$  
almost surely and in $L^q$, $q\ge 1$ arbitrary.\\
iii) For continuous functions $g:\bbr\to\bbr$ we have  
$$
\sum\limits_{\al\in\calj(\Delta)} {\rm length}(A_\al)\; g(\tfrac{\calx_\al-a}{h})\; K_{h}(\calx_\al-a)  \;\;\lra\;\;  \int_{-1}^1 g(v)\, K(v)\, dv  \quad,\quad a \in {\rm int}(A) 
$$ 
almost surely and in $L^q$, $q\ge 1$ arbitrary. 
\end{lemma}

{\bf Proof: } For (\ref{*NEU}) it is sufficient to note that $K$ is bounded, 
that $\calx_\al$ belongs to cell $A_\al$ by construction of the regression scheme in \ref{3.1.1NEU}, and that the number of cells $A_\al$ which inter\-sect the support of $K( \frac{\cdot{-}a}{h} )$ is  $\,\calo( n h )\,$ which tends to $\infty$: so the left hand side of   (\ref{*NEU})  is of type 
$$
\frac1n\sum\limits_{\al=1}^n  \frac{1}{h}\,|K|^m(\tfrac{\calx_{\al}-a}{h})
$$
and assertion i) is proved. To prove ii) and iii), the same argument shows that for $\Delta$ small enough,  
at points $a\in {\rm int}(A)$, the random variables 
$$
\frac{\calx_\al-a}{h} \quad\mbox{where $\al$ is such that}\quad  \frac{\calx_\al-a}{h}  \;\in\; [-1,1]
$$
are approximately equispaced over $[-1,1]$, the spacing of the design variables $\calx_\al$ being of order $\;{\rm length}(A_\al) \sim \Delta^{\frac 12} = \calo(\frac 1n)\,$. As a consequence, left hand sides in ii) and iii) are Riemann sums 
and converge almost surely as $\Delta\downarrow 0$. Since (\ref{*NEU}) provides constants $M = M(K,1)$ such that the convergence is dominated, we have also convergence in $L^q$ for $q\ge 1$ arbitrary.\halmos\\

\begin{lemma}\label{3.3.xxNEU}
For the kernel $K$ of order $\beta'$, we have deterministic bounds 
$$
\sup_{a \in {\rm int}(A)}\;\;  \limsup_{\Delta\downarrow 0}\quad
nh \left|\; \sum\limits_{\al\in\calj(\Delta)} {\rm length}(A_\al)\, ( \calx_\al-a )^j\, K_h( \calx_\al-a ) \;\;-\; 0 \;\right|  
\quad\le\;\; M\;\; <\quad  \infty   
$$    
for every $1\le j \le \beta'$. Recall that $\beta'$ is the greatest integer strictly smaller than $\beta$. 
\end{lemma}

\vskip0.5cm
{\bf Proof: } Since $K$ is a kernel of order $\beta'$, we know from lemma \ref{3.3.1NEU}~iii) that for every $1\le j\le \beta'$
$$
 \sum\limits_{\al\in\calj(\Delta)} {\rm length}(A_\al)\, \left( \frac{\calx_\al-a}{h} \right)^j K_h(\calx_\al-a)   \quad\lra\quad  \int_{-1}^1 v^j\, K(v)\, dv  \;=\;  0 \quad,\quad 1\le j\le \beta'
$$ 
holds almost surely as $\,\Delta\downarrow 0\,$, and in $L^q$ for arbitrary $q\ge 1$. We shall combine this with 
the steps of the proof of lemma \ref{3.3.2NEU}. Fix $1\le j\le \beta'$ and select $\zeta_\al=\zeta_\al(j) \in {\rm cl}(A_\al)$ such that 
$$
\int_{A_\al} (x-a)^j\,K_h(x-a)\, dx \;=\; {\rm length}(A_\al)\; (\zeta_\al -a)^j\, K_h(\zeta_\al -a)   \;. 
$$
Since $x \to (x-a)^j\,K_h(x-a)$ is Lipschitz on $B_h(a)$ and $\calx_\al\in A_\al$, we have 
$$
\left|\, (\zeta_\al -a)^j\, K\Big( \frac{\zeta_\al-a}{h} \Big) - (\calx_\al -a)^j\, K\Big( \frac{\calx_\al-a}{h} \Big) \,\right| \;\;\le\;\; \calo\left( {\rm length}(A_\al) \frac 1 h \right)\; \;=\;\calo\left( \frac{1}{nh} \right)    
$$
and thus 
$$
\int_{A_\al} (x-a)^j\,K_h(x-a)\, dx \;=\; {\rm length}(A_\al)\; (\calx_\al -a)^j\, K_h(\calx_\al -a) \;+\; \calo\left(\frac{1}{(nh)^2}\right)  \;. 
$$
For $\calo(nh)$ indices $\al$ in $\calj(\Delta)$ the cell $A_\al$ will be fully contained in $B_h(a)$; for at most two additional values of $\al$, a cell $A_\al$ may intersect $B_h(a)$. Since $a\in{\rm int}(A)$, we have for $\Delta$ small enough 
\beao
0 &=&  \int_{B_h(a)} (x-a)^j\,K_h(x-a)\, dx \;\;=\; \sum_{\al\in\calj(\Delta)} \int_{A_\al} (x-a)^j\,K_h(x-a)\, dx \\
&=& \sum_{\al\in\calj(\Delta)}{\rm length}(A_\al)\; (\calx_\al -a)^j\, K_h(\calx_\al -a) 
\;\;+\;\; \calo\left(\frac{1}{nh}\right)  
\eeao
for every $1\le j\le \beta'$ since $K$ is a kernel of order $\beta'$. This is the assertion.  \halmos\\

\begin{lemma}\label{3.3.4NEU}
Our assumptions on the kernel $K$ combined with the H\"older property for $\si^2$ grant  
$$
\sup_{a \in {\rm int}(A)}\;\;  \limsup_{\Delta\downarrow 0}\quad
E_{\si^2}\left(\, nh \left[\, \sum\limits_{\al\in\calj(\Delta)} {\rm length}(A_\al)\;  \left( \si^2(\calx_\al) - \si^2(a) \right)\;  K_h(\calx_\al-a) \,\right]^2 \,\right) \quad<\quad \infty   \;. 
$$
\end{lemma}

\vskip0.8cm
{\bf Proof: } Fix $a\in {\rm int}(A)$. We start from a Taylor expansion of $\,f:=\si^2\in\calh(\beta,L)\,$ at $\,a\,$  as in  \cite{Tsy-08}~p.~14   
$$
f(a+h) - f(a) \;=\; \sum_{j=1}^{\beta'}  \frac{ f^{(j)}(a) }{j!} h^j \;+\; \frac{h^{\beta'}}{(\beta'-1)!} \int_0^1 (1-\tau)^{(\beta'-1)} [ f^{(\beta')}(a+\tau h) - f^{(\beta')}(a) ] \, d\tau
$$
and apply the preceding lemmata. First, derivatives of order $1\le j\le \beta'$ in the above square brackets produce terms 
$$
\sum\limits_{\al\in\calj(\Delta)} {\rm length}(A_\al)\;  (\, \calx_\al-a \,)^j\; K_h(\calx_\al-a) 
$$
to be multiplied with deterministic factors $\frac{(\si^2)^{(j)}(a)}{j!}$. By lemma \ref{3.3.xxNEU}, such terms admit deterministic bounds $\calo(\frac{1}{nh})$ as $\Delta\downarrow 0$. Second, writing $\psi$ for the derivative of order $\beta'$ of $\si^2$, remainder terms in the above square brackets take the form  
$$ 
\sum\limits_{\al\in\calj(\Delta)} {\rm length}(A_\al)\;  K_h(\calx_\al-a)\; \frac{ (\calx_\al-a)^{\beta'} }{(\beta' - 1)!}\; \int_0^1 (1-\tau)^{\beta' -1} \left(\, \psi\left( a + \tau (\calx_\al-a) \right) - \psi( a ) \,\right) d\tau \;.
$$
Now $\psi$ being H\"older of order $\beta-\beta'$ with H\"older constant $L$, this is bounded in absolute value by 
$$
\sum\limits_{\al\in\calj(\Delta)} {\rm length}(A_\al)\;  \frac{ L\, |\calx_\al-a |^{\beta} }{ \beta' !} \; \frac{1}{h} |K|( \tfrac{\calx_\al-a}h ) 
$$
which using \eqref{*NEU} with $M=M(K,1)$, \eqref{n_delta_lambda_d1NEU}, \eqref{bandwidthNEU} and \eqref{nhrateNEU} admits deterministic bounds 
$$
\frac{L\, M}{\beta' !}\; h^{\beta} \;=\; \calo(h^{\beta}) \;=\; \calo( n^{\frac{ - \beta}{2\beta+1}}) \;=\; \calo( \frac{1}{\sqrt{nh\,}} )
$$ 
as $\Delta\downarrow 0$. Combining these bounds concludes the proof of the lemma.\halmos\\

\begin{lemma}\label{3.3.yyNEU}
In restriction to the good events $\,{\tt G}(\Delta)\,$ of theorem \ref{3.1.4NEU} we have for $\caly_\al=\calz^2_\al$
$$
\sup_{a \in {\rm int}(A)}\;\;  \limsup_{\Delta\downarrow 0}\quad
E_{\si^2}\left(\, nh \left[\, \sum\limits_{\al\in\calj(\Delta)} {\rm length}(A_\al)\;  1_{{\tt G}(\Delta)}\!\!\left(\, \caly_\al - \si^2(\calx_\al) \right)  K_h(\calx_\al-a) \,\right]^2  \,\right) \quad<\quad \infty \;.   
$$ 
\end{lemma}

\vskip0.8cm
{\bf Proof: } We use the notations of \ref{3.1.4NEU}, \eqref{good_approximationNEU}, \eqref{interpretation_schemeNEU} and \ref{3.1.5NEU}: in particular we have in dimension  $d=1$ on the good event $\,{\tt G}(\Delta)\,$
$$
\caly_\al \;=\; \calz^2_\al \;\;\approx\;\;  [\si(\calx_\al)\,\calu_\al(1)]^2 \;=\;
\si^2(\calx_\al)\left\{ 1 + 2 \int_0^1\calu_\al(s)d\calu_\al(s) \right\} 
\quad,\quad \al\in\calj(\Delta)
$$
with one-dimensional standard Brownian motions $\,\calu_\al$, $\al\in\calj(\Delta)$, which by construction in theorem \ref{3.1.4NEU} are independent of each other and independent of the $\calx_\al$, $\al\in\calj(\Delta)$. Thus as $\Delta\downarrow 0$, the expression in square brackets in the assertion is the sum in $L^2(Q_{\si^2})$ of two terms: first, 
\beqq\label{have_to_prove_contibution3NEU}
\cals_1(\Delta) \;:=\; 
\sum\limits_{\al\in\calj(\Delta)} {\rm length}(A_\al)\;  1_{{\tt G}(\Delta)}\!\!\left[\, \caly_\al - [\si(\calx_\al)\,\calu_\al(1)]^2 \right]  K_h(\calx_\al-a)   \; 
\eeqq
second, since $\;\calu^2_\al(1) = 1 + 2 \int_0^1\calu_\al(s)d\calu_\al(s)\,$ on {\tt G}$(\Delta)$, 
\beqq\label{have_to_prove_contibution2NEU}
\cals_2(\Delta) \;:=\; 
2 \sum\limits_{\al\in\calj(\Delta)} {\rm length}(A_\al)\; \si^2(\calx_\al)\;  1_{{\tt G}(\Delta)}\!\! 
\int_0^1\calu_\al(s)d\calu_\al(s)\;   K_h(\calx_\al-a)   \;. 
\eeqq

1) Consider \eqref{have_to_prove_contibution2NEU} first: the $\,\calu_\al$ being independent of each other and independent of the  $\calx_\al$,  $\al\in\calj(\Delta)$,  
$$
E_{\si^2}\left(\, \left[\, \cals_2(\Delta)\,\right]^2  \,\right) 
\;=\; 
4\; E_{\si^2}\left(\, \left[\, \sum\limits_{\al\in\calj(\Delta)} {\rm length}(A_\al)\;  \si^2(\calx_\al)\; 
1_{{\tt G}(\Delta)}\!\! \int_0^1\calu_\al(s)d\calu_\al(s)\;    
K_h(\calx_\al-a) \,\right]^2  \,\right) 
$$
reduces to the expectation of the sum of squares of diagonal terms which admits bounds 
$$
\le\quad 4\;\; \frac{L^4}{2}\;\;  \frac{1}{nh}\;\; E_{\si^2}\left(\,  \sum\limits_{\al\in\calj(\Delta)} {\rm length}(A_\al)\;    \frac 1h 
|K|^2(\tfrac{\calx_\al-a}{h}) \,\right) \quad\le\;\; 2\, L^4 M\;\;  \frac{1}{nh}
$$
where $\,\frac 12\,$ is the expectation  $\,E\left( (\int_0^1 \calu_1 d\calu_1)^2 \right)= E\left( \int_0^1 \calu^2_1(s) ds \right)\,$, 
$\,L$ the bound for $|\si|$ from \eqref{additional_assumption_J+GCNEU}, and $M=M(K,2)$ the deterministic bound from \eqref{*NEU}. So the sum \eqref{have_to_prove_contibution2NEU} satisfies 
$$
E_{\si^2}\left(\, \left[\, \cals_2(\Delta)\,\right]^2 \,\right) \;=\; \calo\left(\frac{1}{nh}\right) \qquad\mbox{as $\;\Delta\downarrow 0$}\;. 
$$

2) To deal with $\left[\, \cals_1(\Delta)\,\right]^2$ from  \eqref{have_to_prove_contibution3NEU}, we start with the sum of squared diagonal terms  
$$
\cals_3(\Delta) \;:=\; 
\sum\limits_{\al\in\calj(\Delta)} {\rm length}^2(A_\al)\;  1_{{\tt G}(\Delta)}\!\!\left[\, \caly_\al - [\si(\calx_\al)\,\calu_\al(1)]^2 \right]^2  K_h^2(\calx_\al-a)   
$$
to which lemma \ref{3.1.5NEU} applies (with $g(z)=z^2$ since $\caly_\al=\calz^2_\al$): placing conditional expectations \eqref{GJ-93_proposition_6NEU} inside $E_{\si^2}\left(\ldots\right)$ and using bounds  $M=M(K,2)$ from \eqref{*NEU} together with \eqref{n_delta_lambda_d1NEU} and \eqref{bandwidthNEU}, we obtain
$$
E_{\si^2}\left(\, \cals_3(\Delta)\, \right) \;\;\le\;\; C\Delta\;\; M\;  \frac{1}{nh}  
\;=\; \calo\left(\frac{\Delta}{nh}\right)   \;=\; \calo\left(\frac{1}{n^3h}\right)
$$
which is negligible as $\Delta\downarrow 0$ in comparison to the bound $\calo\left(\frac{1}{nh}\right)$ in step 1). 

3) To deal with the sum of non-diagonal contributions 
$$
\cals_4(\Delta):= \left[\, \cals_1(\Delta)\,\right]^2 - \cals_3(\Delta)  
$$
to $\left[\, \cals_1(\Delta)\,\right]^2$ we introduce short notations  
$$
\calv_{\al'} \;:=\; 1_{{\tt G}(\Delta)}\!\!\left[\, \caly_{\al'} - [\si(\calx_{\al'})\,\calu_{\al'}(1)]^2 \right]  \;\;,\;\; \al'\in\calj(\Delta) \quad,\quad  \calg \;:=\; \si\left(\, \calx_\al : \al\in\calj(\Delta) \,\right) 
$$
and write 
$$
\cals_4(\Delta) \;=\; \sum\limits_{\al'\neq\al''} {\rm length}(A_{\al'}) {\rm length}(A_{\al''})\; 
\calv_{\al'}\; \calv_{\al''}\; 
K_h(\calx_{\al'}-a)  K_h(\calx_{\al''}-a) \;.
$$
Using a regular version $K(\om,\cdot)$ of the conditional law   
$$
\call\left(\,
\left(\,  1_{ {\tt G}(\Delta) }\, \left( \caly_\al \atop \calu_\al \right)  \,\right)_{\al\in\calj(\Delta)}  
\,\bigg|\,  \calg  \,\right)   
$$
and Cauchy-Schwarz with respect to $K(\om,\cdot)$, lemma \ref{3.1.5NEU} again applies (with $g(z)=z^2$) and yields deterministic bounds 
$$
E\left(\, | \calv_{\al'} \calv_{\al''} | \mid \calg \right) \;\le\; \sqrt{ E\left(\, \calv_{\al'}^2\!  \mid \calg \right)  E\left(\, \calv^2_{\al''}\! \mid \calg \right)\;  } \;\;\le\; \; C\,\Delta
$$
almost surely, whence  
\beao
E_{\si^2}\left(\, |\cals_4(\Delta)| \,\right) \;\;\le\;\;
C\, \Delta\;\; E_{\si^2}\left[ \sum\limits_{\al\in\calj(\Delta)} {\rm length}(A_\al)\; \frac 1h |K|(\tfrac{\calx_\al-a}{h}) \right]^2 
\;\;\le\;\; CM^2\; \Delta 
\;\;=\;\; \calo\left(\frac{1}{n^2}\right)
\eeao
using again the bounds $M=M(K,1)$ from \eqref{*NEU}. Thus also 
$E_{\si^2}\left(\, |\cals_4(\Delta)| \,\right)$
is negligible in comparison to the bound $\calo\left(\frac{1}{nh}\right)$ obtained in step 1) as $\Delta\downarrow 0$. 

4) As a consequence of steps 1) to 3), we have 
$$
E_{\si^2}\left(\, \left[\, \cals_1(\Delta)\,\right]^2 \;+\; \left[\, \cals_2(\Delta)\,\right]^2 \,\right) 
\;\;=\;\; \calo\left(\frac{1}{nh}\right) 
$$
which finishes the proof of lemma \ref{3.3.yyNEU}. \halmos\\

\begin{lemma}\label{3.3.zzNEU}
In general regression schemes \ref{3.1.1NEU} where $0<\la<\frac 12$,  exceptional events $\,{\tt F}(\Delta)\,$ in theorem \ref{3.1.4NEU} are such that 
$$
\sup_{a \in {\rm int}(A)}\;\;  \limsup_{\Delta\downarrow 0}\quad
E_{\si^2}\left(\, n^{8\la-3} \left[\, \sum\limits_{\al\in\calj(\Delta)} {\rm length}(A_\al)\;  1_{{\tt F}(\Delta)}\!\!\left(\, \caly_\al - \si^2(\calx_\al) \right)  K_h(\calx_\al-a) \,\right]^2  \,\right) \quad<\quad \infty \;.   
$$ 
Using $\la=\la_0(\beta)$, the critical value \eqref{critical_valueNEU} associated to class $\calh(\beta,L)$, we have  
$$
E_{\si^2}\left(\, \left[\, \sum\limits_{\al\in\calj(\Delta)} {\rm length}(A_\al)\;  1_{{\tt F}(\Delta)}\!\!\left(\, \caly_\al - \si^2(\calx_\al) \right)  K_h(\calx_\al-a) \,\right]^2  \,\right) \;\;=\;\; \calo\left( \frac{1}{nh} \right)  
$$
as $\Delta\downarrow 0$. The left hand side is negligible in comparison to  $\calo\left( \frac{1}{nh} \right)$ when $\,\la_0(\beta)<\la<\frac 12\,$. 
\end{lemma}

\vskip0.5cm
{\bf Proof: } Write $\cals_5(\Delta)$ for the sum in square brackets in the assertion. 
Then  $\caly_\al=\calz^2_\al$ combined with the deterministic bounds on $\calz_\al$ from  part iv) of theorem  \ref{3.1.4NEU} give 
$$
\left| \caly_{\al'} - \si^2(\calx_{\al'}) \right|   \left| \caly_{\al''} - \si^2(\calx_{\al''}) \right|  
\;\;\le\;\;   \calo( \Delta^{4(\la -\frac 12)} ) \quad,\quad \al',\al''\in\calj(\Delta)
$$
when $\Delta$ is small enough (recall that $\si^2$ is bounded on the interval $A$). Using $\Delta=\calo(n^{-2})$ from \eqref{n_delta_lambda_d1NEU}, this bound is $\calo\left( n^{4-8\la} \right)$ as $\Delta\downarrow 0$.  From part iii) of theorem  \ref{3.1.4NEU} and \eqref{n_delta_lambda_d1NEU} we have 
$$
Q_{\si^2}\left(\, {\tt F}(\Delta) \,\right) \;\;\le\;\; \calo\left(\Delta^{\frac 12}\right) \;=\; \calo\left( \frac 1n \right)
$$ 
as $\Delta\downarrow 0$. Proceeding as in the proof of lemma \ref{3.3.yyNEU}, using the  constants $M=M(K,1)$ of \eqref{*NEU}, we end up with the bound  
$$
E_{\si^2}\left(\, [\cals_5(\Delta)]^2 \,\right) \;\;\le\;\; \calo\left( n^{3-8\la} \right)   \quad\mbox{as $\Delta\downarrow 0$}   
$$
which proves the first assertion. The second assertion follows since $\la=\la_0(\beta)$ in \eqref{critical_valueNEU} is such that 
$$
\calo\left( n^{3-8\la_0(\beta)} \right) \;=\; \calo\left(\; \frac 1n\; n^{ \frac{1}{2\beta+1} } \right)  
\;=\; \calo\left( \frac{1}{nh} \right)  
$$
by definition of the bandwidth in \eqref{bandwidthNEU}.  \halmos\\

{\bf Proof of theorem \ref{3.2.2'NEU}: } Note first that $\, n^{\frac{2\beta}{2\beta+1}} = nh\,$ by \eqref{n_delta_lambda_d1NEU} and \eqref{bandwidthNEU}. Then the estimation error 
$$
\wh{\si^2_\Delta}(a) \;-\; \si^2(a)  
$$
is decomposed into several terms. First,  the difference 
$$
\si^2(a)  \;-\; \sum\limits_{\al\in\calj(\Delta)} {\rm length}(A_\al)\; \si^2(a)\;   K_h(\calx_\al-a)   
$$
is $\calo(\frac{1}{nh})$ by lemma \ref{3.3.2NEU}. Second, by lemma \ref{3.3.4NEU}, 
$$
\sum\limits_{\al\in\calj(\Delta)} {\rm length}(A_\al)\; (\, \si^2(\calx_\al) - \si^2(a) \,)\;   K_h(\calx_\al-a)   
$$
has squared $L^2(Q_{\si^2})$-norm of order $\calo(\frac{1}{nh})$. Third, on the good events ${\tt G}(\Delta)$ of theorem \ref{3.1.4NEU}, 
$$
\sum\limits_{\al\in\calj(\Delta)} {\rm length}(A_\al)\;  1_{{\tt G}(\Delta)}\!\!\left(\, \caly_\al - \si^2(\calx_\al) \right)  K_h(\calx_\al-a)
$$
has squared $L^2(Q_{\si^2})$-norm of order $\calo(\frac{1}{nh})$ by lemma \ref{3.3.yyNEU}.  So far, we could work with arbitrary $0<\la<\frac 12$ fixed. This situation changes drastically with the final contribution 
$$
\sum\limits_{\al\in\calj(\Delta)} {\rm length}(A_\al)\;  1_{{\tt F}(\Delta)}\!\!\left(\, \caly_\al - \si^2(\calx_\al) \right)  K_h(\calx_\al-a)
$$
of the exceptional events ${\tt F}(\Delta)$: here we have not more than the trivial bounds from theorem \ref{3.1.4NEU}~iv). 
By lemma \ref{3.3.zzNEU}, squared $L^2(Q_{\si^2})$-norms are of order $\calo(n^{3-8\lambda})$ which obliges us to work with 
$$
\la \;\ge\; \la_0(\beta) \;=\; \frac 12 - \frac{1}{8(2\beta+1)} \;, 
$$
the condition introduced in \eqref{critical_valueNEU} and \eqref{critical_value_bedingungNEU}, to get the contribution from exceptional events balanced under a common $\le \calo(\frac{1}{nh})$ for all contributions. The proof of theorem \ref{3.2.2'NEU} is finished. \halmos\\

%
%



\begin{thebibliography}{99}

\bibitem{ADR-69}
Az\'{e}ma, J., Duflo, M., Revuz, D.: Mesure invariante des processus de Markov r\'{e}currents. Sem. Prob. Strasbourg \textbf{3}, 24--33 (1969).

\bibitem{BDMT-11}
Bansaye, V., Delmas, J.-F., Marsalle, L., Tran, V.: Limit theorems for Markov processes indexed by continuous time Galton-Watson trees. The Annals of Applied Probability \textbf{21(6)}, 2263–-2314 (2011).

\bibitem{Be-15}
Berg, T.: Nonparametric estimation of the diffusion coefficient of a branching diffusion with immigration. PhD thesis, Mainz 2015.\\ 
{\tt https://publications.ub.uni-mainz.de/theses/volltexte/2015/4096/pdf/4096.pdf} 

\bibitem{BS-95}
Bibby, B., S\o rensen, M.: Martingale estimation functions for discretely observed diffusion processes. Bernoulli {\bf 1(1/2)}, 17--339 (1995). 



\bibitem{Br-05}
Brandt, C.: Partial reconstruction of the trajectories of a discretely observed branching diffusion with immigration and an application to inference. PhD thesis, Mainz 2005.\\ 
{\tt https://publications.ub.uni-mainz.de/theses/volltexte/2005/756/pdf/756.pdf}

\bibitem{Ca-90} 
Cattiaux, P.: Calcul stochastique et op\'erateurs d\'eg\'en\'er\'es du second ordre I. R\'esolvantes, th\'eor\`eme de H\"ormander et applications. 
Bull.\ Sci.\ Math.\ {\bf 114}, 421-462 (1990).  

\bibitem{Dy-65}
Dynkin, E.: Markov processes Vol.\ II. Springer 1965. 

\bibitem{Els-11}
Elstrodt, J.: Ma{\ss}- und Integrationstheorie. 7th Edition, Berlin: Springer 2011.

\bibitem{EN-00}
Engel, K., Nagel, R.: One-Parameter Semigroups for Linear Evolution Equations. New York: Springer 2000.


\bibitem{Fo-92}
Folland, G.: Fourier analysis and its applications. Wadsworth 1992. 

\bibitem{Fr-75}
Friedman, A.: Stochastic differential equations and applications Vol.\ I. Academic Press 1975. 

\bibitem{GJ-93}
Genon-Catalot, V., Jacod, J.: On the estimation of the diffusion coefficient for multi-dimensional diffusion processes. Ann.\ Inst.\ H.\ Poincar\'e (Proba.\ Stat.) {\bf 29(1)}, 119--151 (1993). 

\bibitem{Go-02}
Gobet, E.: LAN property for ergodic diffusions with discrete observations. Ann.\ Inst.\ H.\ Poincar\'e {\bf 38}, 711--733 (2002). 

\bibitem{Ha-12}
Hammer, M.: Ergodicity and Regularity of Invariant Measure for Branching Markov Processes with Immigration. PhD thesis, Mainz 2012.\\ 
{\tt https://publications.ub.uni-mainz.de/theses/volltexte/2012/3306/pdf/3306.pdf} 

\bibitem{HH-09} 
Hardy, R., Harris, S.: A spine approach to branching diffusions with applications to $\mathcal{L}^p$-convergence of martingales. In S\'{e}minaire de Probabilit\'{e}s XLII, pages 281–330. Springer (2009).

\bibitem{HS-65}
Hewitt, E., Stromberg, K.: Real and abstract analysis. Springer 1965. 

\bibitem{HHL-02} 
H\"opfner, R., Hoffmann, M., L\"ocherbach, E.: Non-parametric estimation of the death rate in branching diffusions. Scand.\ J.\ Statist.\ {\bf 29}, 665--692 (2002). 

\bibitem{HL-05} 
H\"opfner, R., L\"ocherbach, E.: Remarks on ergodicity and invariant occupation measure in branching diffusions with immigration. Ann.\ Instit.\ H.\ Poincar\'e (Proba.\ Statist.) {\bf 41(6)}, 1025--1047 (2005).

\bibitem{INW-66a} 
Ikeda, N., Nagasawa, M., Watanabe, S.: A construction of Markov processes by piecing out. Proc.\ Japan Acad.\ {\bf 42}, 370-375 (1966). 

\bibitem{INW-66b} 
Ikeda, N., Nagasawa, M., Watanabe, S.: A construction of Markov branching processes. Proc.\ Japan Acad.\ {\bf 42}, 380-384 (1966). 

\bibitem{INW-68} 
Ikeda, N., Nagasawa, M., Watanabe, S.: Branching Markov Processes II.  J.\ Math.\ Kyoto Univ.\ {\bf 8}, 365--410 (1968).

\bibitem{INW-69} 
Ikeda, N., Nagasawa, M., Watanabe, S.: Branching Markov Processes III. J.\ Math.\ Kyoto Univ.\ {\bf 9}, 95--160 (1969).

\bibitem{JP-12} 
Jacod, J., Protter, P.: Discretization of processes. Springer 2012. 

\bibitem{Ke-97}
Kessler, M.: Estimation of an ergodic diffusion from discrete observations. Scand.\ J.\ Statist.\ {\bf 24(2)}, 211--229 (1997).



\bibitem{LL-01} 
Lieb, E. H., Loss, M.: Analysis. Second Edition, Providence: American Mathematical Society 2001. 

\bibitem{LKS-12} 
Lindner, A., Kessler, M., S\o rensen, M.: Statistical Methods for Stochastic Differential Equations. CRC Press 2012. 

\bibitem{Lo-02} 
L\"ocherbach, E.: Likelihood ratio processes for Markovian particle systems with killing and jumps. Statist.\ Inference Stoch.\ Proc.\ {\bf 5(1)}, 153--177, 2002. 

\bibitem{Lo-02a} 
L\"ocherbach, E.: LAN and LAMN for systems of interacting diffusions with branching and immigration. Ann.\ I.\ H.\ Poincar\'e {\bf 38(1)}, 59--90, 2002. 

\bibitem{Lo-04} 
L\"ocherbach, E.: Smoothness of the intensity measure density for interacting branching diffusions with immigrations. J.\ Funct.\ Analysis {\bf 215(1)}, 130--177 (2004). 

\bibitem{Ma-18} 
Marguet, A.: Uniform sampling in a structured branching population. \\
https://arxiv.org/abs/1609.05678 (2018).

\bibitem{Na-77} 
Nagasawa, M.: Basic models of branching processes. Bull.\ Inst.\ Intern.\ Statist.\ {\bf 27(2)}, 423--445 (1977). 


\bibitem{Nu-78}
Nummelin, E.: A splitting technique for Harris recurrent Markov chains. Z.\ Wahrscheinlichkeitsth.\ Verw.\ Geb.\ {\bf 43}, 309--318 (1978). 

\bibitem{Nu-85}
Nummelin, E.: General irreducible Markov chains and non-negative operators. Cambridge University Press 1985. 


\bibitem{PV-10} 
Podolskij, M., Vetter, M.: Understanding limit theorems for semimartingales: a short survey. Statistica Neerlandica {\bf 64}, 329--351 (2010). 

\bibitem{Re-84}
Revuz, D.: Markov chains. Rev.\ Ed.\ Springer 1984. 

\bibitem{Sa-99}
Sato, K.: L\'evy processes and infinitely divisible distributions. Cambridge University Press 1999. 

\bibitem{Sc-05}
Schilling, R.: Measures, integrals and martingales. Cambridge University Press 2005. 

\bibitem{St-08} 
Stroock, D.: Partial differential equations for probabilists. 
Cambridge University Press 2008. 


\bibitem{Tsy-08}
Tsybakov, A.: Introduction to nonparametric estimation. Springer 2008. 

\bibitem{Wa-67} 
Watanabe, S.: Limit Theorem for a Class of Branching Processes. In: J. Chover (ed.), Markov processes and potential theory. Proc. Sympos. Math. Res. Center Madison Wis. 1967, pp. 205-232. New York: Wiley 1967.


\bibitem{Yo-92}
Yoshida, N.: Estimation for diffusion processes from discrete observations. J.\ Multivariate Anal.\  {\bf 41(2)}, 220–-242 (1992).

\end{thebibliography}
\end{document}